\documentclass[article]{siamart190516}

\usepackage{lipsum}
\usepackage{amsfonts}
\usepackage{graphicx}
\usepackage{epstopdf}
\usepackage{algorithmic}
\ifpdf
  \DeclareGraphicsExtensions{.eps,.pdf,.png,.jpg}
\else
  \DeclareGraphicsExtensions{.eps}
\fi
\usepackage{amsmath} 
\usepackage{amssymb}  
\usepackage{mathrsfs}
\usepackage{algorithm}
\usepackage{subfigure}
\usepackage{multirow}
\usepackage{appendix}
\usepackage{enumerate}
\usepackage{mathtools}
\usepackage{color}
\usepackage{bm}

\def\opi{{\mathrm{i}\mkern1mu}}

\newcommand{\bC}{\mathbb{C}}
\newcommand{\bR}{\mathbb{R}}
\newcommand{\real}{\operatorname*{Re}}
\newcommand{\imag}{\operatorname*{Im}}

\newcommand{\Id}{\mathrm{I}}
\newcommand{\bmu}{\boldsymbol{\mu}}
\newcommand{\bff}{\boldsymbol}
\newcommand{\norm}[1]{\left\lVert#1\right\rVert}
\DeclareMathOperator*{\argmax}{arg\,max}
\newcommand{\e}{\mathrm{e}}
\DeclareMathOperator*{\argmin}{arg\,min}
\renewcommand{\hat}{\widehat}

\newcommand\cD{\mathcal{D}}

\newcommand\cL{\mathcal{L}}
\newcommand\cM{\mathcal{M}}

\newcommand\C{\mathbb{C}}
\newcommand\F{\mathbb{F}}

\newcommand\R{\mathbb{R}}

\newcommand\bA{\bm{A}}

\newcommand\bb{\bm{b}}
\newcommand\bbf{\bm{f}}

\newcommand\bu{\bm{u}}

\newsiamremark{remark}{Remark}
\newsiamremark{hypothesis}{Hypothesis}
\crefname{hypothesis}{Hypothesis}{Hypotheses}
\newsiamthm{claim}{Claim}

\headers
{Model order reduction contour integral methods}{N. Guglielmi, and M. Manucci}

\title{Model Order Reduction in contour integral methods for parametric PDEs}

\author{Nicola Guglielmi\thanks{Gran Sasso Science Institute, L'Aquila, IT 
  (\email{nicola.guglielmi@gssi.it}).}
\and Mattia Manucci\thanks{Gran Sasso Science Institute, L'Aquila, IT 
  (\email{mattia.manucci@gssi.it}).}
}

\usepackage{amsopn}

\ifpdf
\hypersetup{
	pdftitle={Model order reduction in contour integral methods for parametric PDEs},
	pdfauthor={N. Guglielmi \and M. Manucci}
}
\fi

\begin{document}
	\setlength{\abovedisplayskip}{3pt}
	\setlength{\belowdisplayskip}{3pt}
	\maketitle
	
	\begin{abstract}
		In this paper we discuss a {projection model order reduction (MOR)} method for a class of parametric linear evolution PDEs, which is based on the application of the Laplace transform. The main advantage
		of this approach consists in the fact that, differently from time stepping methods, like Runge-Kutta integrators, the Laplace transform allows to compute the solution directly at a given instant, which can be done by approximating the contour integral associated to the inverse Laplace transform by a suitable quadrature formula. 
		In terms of {some classical MOR} methodology, this determines a significant improvement in the reduction phase - like the one based on the classical proper orthogonal decomposition (POD) - since the number of vectors to which the decomposition applies is drastically reduced as it does not contain all intermediate solutions generated along an integration grid by a time stepping method.
		We show the effectiveness of the method by some illustrative parabolic PDEs arising from finance and also provide some evidence that the method
		we propose, when applied to a linear advection equation, does not suffer the
		problem of slow decay of singular values which instead affects time stepping methods for the numerical approximation of the Cauchy problem arising from space discretization. 
	\end{abstract}
	
	\vspace{-0.2cm}
	
	\vspace{-0.3cm}
	\vspace{-0.4cm}
	\section{Introduction}
	Parametrized time dependent partial differential equations (PPDEs) occur in many contexts such as industrial or financial applications. They are solved by means of high-fidelity (or full-order) approximation techniques such as finite elements, finite volumes, finite differences or spectral methods coupled with an integrator for the {arising} Cauchy problem. However, the use of high-fidelity approximation techniques may become prohibitive when they are required to deal quickly and efficiently with a repetitive solution of PPDEs with different parameters.
	
	For this class of problems, \textit{reduced-order modeling} – also named \textit{model order reduction} – is a generic expression used to identify any approach aimed at replacing the high-fidelity problem by one featuring with a much lower numerical complexity, relying on a controllable error with respect to the high-fidelity solution. The underlying assumption on which the \textit{reduced-order model} (ROM) techniques are based is that the behaviour of a possibly complex system may be described by a small number of dominant modes.
	
	In particular, among the reduced-order modeling techniques, a remarkable instance is represented by \textit{reduced basis} (RB) methods. The strategy is to solve, during a computationally onerous \textit{offline} phase, the high-fidelity problem only for few instances of the input parameters, aiming to construct a set of base solutions, the so called \textit{reduced base}. This \textit{reduced space} is usually of dimension much smaller than the space associated to the high fidelity problem. After that, for every new instance of the input parameters, the associated solution is approximated by a suitable linear combination of the base elements of the reduced space. The unknown coefficients of this combination are obtained during the \textit{online} phase by solving a simpler problem generated through a Galerkin projection of the high fidelity problem onto the reduced space. For a detailed survey on reduced basis methods starting from the framework of the continuous problems, we remand to \cite{Hesthaven2016}.
	
	In this work we are interested in time-dependent problems and their treatment in the {MOR} context. A survey on this topic{, from a reduced basis prospective,} is given in \cite{Glas2017} where two methodologies are described and compared. The first one, and probably the one considered ``standard", is based upon a time stepping solver like a Runge-Kutta method in the offline phase. The reduced basis is then usually formed thanks to a Greedy algorithm which is an iterative procedure where at each iteration one new basis function is added in such a way that the overall approximation capability of the basis set is improved. The greedy strategy is then combined with a compression of the computed numerical solution along the grid through a Proper Orthogonal Decomposition (POD) - that is an SVD decomposition of the array containing the solution vectors -  \cite{Drohmann2012,Haasdonk2013}. The drawback of these time-stepping schemes is that, in order to approximate the solution at a certain time $T=t_n$, one needs to compute an approximation of the solution, for both the full and the reduced problem, at grid points $0<t_1<t_2<\ldots<t_n$, which would be particularly demanding if $T$ is large. The second approach consists in treating time as an additional variable, which results in a problem of dimension $d+1$ where $d$ denotes the spatial dimension. The reduced basis is formed by a standard Greedy algorithm and then the full problem is projected (in the sense of Galerkin or Petrov-Galerkin) onto the reduced space \cite{Urban2012}. It is well known that the size of the discrete problem grows exponentially with respect to the numbers of variables when keeping the same accuracy (the so-called curse of dimensionality), therefore the computational cost for this approach can become prohibitive in the offline phase. {Another survey which discuss about MOR techniques from a system-theoretic
	point of view is \cite{Benner2015}}.   
	
	Here we introduce a new strategy to address time dependence in the {MOR} context. It consists of employing a contour integral method (see e.g. \cite{Guglielmi2018,Guglielmi2020}) to numerically invert the Laplace transform to the the full and reduced problem in order to compute an approximation of the solution at a given time $t$. This method is indicated when one is interested to determine the solution at a specific time $t=T$ or even in certain time window $[t_0,\Lambda t_0]$, with $\Lambda>1$. We will show that the new approach is very efficient when compared to the standard methodologies; moreover it appears to be promising when applied to linear hyperbolic problems which are known to be challenging for the standard {projection reduced order techniques} \cite{Ohlberger2015,Greif2019}.
	
	The article is organized as follows. In Section $2$ we present the mathematical formulation of contour integral methods and projection reduced order methods, we also give a result of Kolmogorov $n$-width decay for discrete parabolic problems. In Section $3$ we describe two types of algorithms for the search of a reduced space based on the new time integrator. In Section $4$ numerical experiments are presented for: 1) the Black-Scholes model, 2) the Heston model and 3) the one dimensional linear advection equation with discontinuous initial data. Finally in Section $5$ we state few conclusions.
	\section{Parameterized evolution problems}
		Let us denote by $\mathcal{A}_{\bmu}$ some elliptic operator (of convection diffusion reaction type) w.r.t. to a spatial variable $\bff x_s\in\Omega\subset\mathbb{R}^d$. For some parameters $\bmu\in\mathcal{D}$ in a compact parameter set $\mathcal{D}\subset \mathbb{R}^p$, we consider the PPDE
	\begin{align}\label{eq:pde}
	    \frac{\partial u}{\partial t}(\bff x_{s},t; \bmu)=& {\cal{A}_{\bmu}} u(\bff x_{s},t; \bmu) +f(\bff x_{s},t; \bmu),\quad \bff x_{s}\in \Omega,\quad t\in(0,T],\nonumber\\
	    &\text{plus boundary condition}\nonumber\\
	    u(0,\bff x_s;\bmu)=&u_0(\bff x_s;\bmu),
	\end{align}
	for some given initial value $u_0$ and external force $f$\footnote{We denote matrices and vectors by boldface characters and function in normal font.}. We assume that the problem is well posed in a strong sense and that it admits a unique solution {$\forall \bmu\in\mathcal{D}$}. 
	\subsection{The detailed discretization: Discrete Cauchy problem}
	
	Next, we perform a sufficiently accurate discretization in space, e.g. by finite elements, finite differences, finite volumes, spectral methods or others. We obtain a Cauchy problem for the coefficients of a discrete expansion in finite dimension, which we call $N_h \in \mathbb{N}$. In particular, we assume that the discrete solution is sufficiently detailed such that we can view it as a “truth” discretization – in the language of the Reduced Basis Method. This means that the number of degrees of freedom $N_h\gg1$ is large, or even huge. The discrete Cauchy problem then reads: find $\mathbf{u}(t; \bmu) \in \mathbb{R}^{N_h}$ such that
	\begin{equation}\label{eq:cauchy}
		\dot{\bff{u}}(t;\bmu) = \bff{A}(\bmu) \bff{u}(t;\bmu) + \bff{b}(t;\bmu), \qquad t\in(0,T], \qquad \bff{u}(0; \bmu)=\bff{u}_0(\bmu),
	\end{equation}
where $\bff{A}(\bmu)\in\mathbb{R}^{N_h\times N_h}$ is the matrix that corresponds to the discretization of $\mathcal{A}_{\bmu}$ and $\bff{b}(t;\bmu)\in\mathbb{R}^{N}$ the discretization of the external force. It is quite common to assume a separation of variables and parameters which is known as \textit{affine decomposition}, i.e.,
\begin{equation}\label{eq:affine2}
		{\bff{A}}({\bmu})=\sum_{q=1}^{Q_{\bff A}}\vartheta_{q}^{\bff A}({\bmu}){\bff{A}}_{q}\quad\text{and} \quad \bff{b}(t;\bmu)=\sum_{q=1}^{{Q}_{\bff b}}\vartheta_{q}^{\bff{b}}({\bmu}){\bff{b}}_q(t),
\end{equation}
where $Q_{\bff A},Q_{\bff b}\in\mathbb{N}$, $\vartheta^{ \bff A}_q$, $\vartheta^{   \bff b}_q:\mathcal{D} \rightarrow \mathbb{R}$ and parameter-independent $\bff{A}_q \in \mathbb{R}^{N_h\times N_h}$, $\bff{b}_q(t)\in \mathbb{R}^{N_h}$.

\subsection{Solution by using the Laplace transform}
The idea of the paper is to use the Laplace transform $\cL[\bff{u}](z;\bmu)\equiv \hat\bu(z;\bmu):=\int_0^\infty \bu(t;\bmu)\, e^{-zt}\, dt$, $z\in\C$, as a backbone of model reduction. Then, the Laplace transform of the solution of \eqref{eq:cauchy} is given by
\begin{align}\label{eq:LT}
	\hat{\bu}(z;\bmu) 
		= \big(z\bff\Id - \bA(\bmu)\big)^{-1}\big(\bu_{0}(\bmu) + \hat\bb(z;\bmu)\big)
			\in\C^{N_h},
			\quad z\in\C,
\end{align}
i.e., we can determine the Laplace transform $\hat\bu(z;\bmu)$ of the desired solution by solving a linear system in $\C^{N_h}$ as long as we can determine the Laplace transform of the force, i.e., $\hat\bb$. 
Once we would know $\hat\bu(z;\bmu)$, this yields the following form of the solution of \eqref{eq:cauchy}, i.e.,

\begin{equation}\label{eq:bromwich}
		\bff u(t;\bmu)=\frac{1}{2\pi \opi}\int_{\gamma-\opi\infty}^{\gamma+\opi\infty}e^{zt}\hat{\bff u}(z; \bmu)\,dz=\frac{1}{2\pi \opi}\int_{\Gamma}e^{zt}\hat{\bff u}(z; \bmu)\,dz,
	\end{equation}
	where $\Gamma$ is an open, piecewise smooth, complex curve\footnote{As opposed to standard notation in PDEs, $\Gamma$ does not denote $\partial\Omega$, the boundary of the domain $\Omega$.} running from $-\opi\infty$ to $+\opi\infty$ and surrounding on its left all the singularities of $\hat{\bff b}(z;\bmu)$ and the eigenvalues of $\bff A(\bmu)$, $\forall \bmu\in\mathcal{D}$. The idea of passing from the integral along the Bromwich contour to the one along $\Gamma$ was first proposed in \cite{Butcher1957}. The purpose of the contour deformation is to exploit the exponential factor. In particular, deforming the path of integration to a Hankel contour, i.e., a contour whose real part begins at negative infinity in the third quadrant, then winds around all singularities and terminates with the real part approaching minus infinity in the second quadrant, we can exploit the  rapid decay due to the exponential factor. This makes the approximation of the integral by the trapezoidal or midpoint rules particularly favourable. The deformation of the contour can be justified by Cauchy’s theorem (see for instance \cite{Butcher1957}), provided the contour is chosen within the domain of analyticity of $\hat{\bff u}(z;\bmu)$. Some mild restrictions on the decay of $\hat{\bff u}(z;\bmu)$ in the left half complex plane are also required (for sufficient conditions, see \cite{TALBOT1979,Trefethen2006}). 

\begin{remark}\label{Rem:AffineLT}
	The affine decomposition \eqref{eq:affine2} can be transferred to the Laplace transform, i.e., 
	\begin{align*}
		\hat\bb(z;\bmu)
		\!=\! \int_0^\infty \bb(t;\bmu)\, e^{-zt}\, dt
		\!=\! \sum_{q=1}^{Q_{\bff b}} \vartheta_q^{\bff b}(\bmu)\, 
				\int_0^\infty\bbf_q(t) \, e^{-zt}\, dt
		= \sum_{q=1}^{Q_{\bff b}} \vartheta_q^{\bff b}(\bmu)\, 
				\hat\bb_q(z),
	\end{align*}
	and similar for $\bu_{0}(\bmu)$. However, for the subsequent investigation, we also need to assume the affine decomposition of $\bff{u}_0(\bmu)$ and a separation w.r.t.\ the variable $z$ for the Laplace transform $\hat\bb(z;\bmu)=\cL[\bb](z;\bmu)$ of the right-hand side of \eqref{eq:LT}, i.e., 
	\begin{align}\label{eq:affinLT}
		\bff{u}_0(\bmu)+ \hat\bb(z;\bmu) = \sum_{q=1}^{Q_{\hat{\bff b}}} 
				\vartheta_{q}^{{\hat{\bff b}}}(\bmu) \,
				\theta_{q}^{{\hat{\bff b}}}(z)\,
				\hat\bb_{q}.
	\end{align}
	Note, that $z\bff\Id - \bA(\bmu)$ admits such a separation with $Q_{\bff A}+1$ terms.
\end{remark}
    \vspace{-0.3cm}
	\subsection{Approximation of the contour integral via  via contour integral methods}
	
	To approximate the integral \eqref{eq:bromwich}, we fix a time $t=T$ and then we parametrize the integration contour $\Gamma$ by a suitable map $z=z(x)\in \mathbb{C}$, $x\in \mathbb{R}$, so that
	\begin{equation}\label{eq:solu}
	\bff u(T;\bmu)=\frac{1}{2\pi \opi}\int_{\Gamma}e^{z{T}}\hat{\bff u}(z; \bmu)\,dz = \frac{1}{2\pi \opi}\int\limits_{\mathbb{R}} \bff G(x; \bmu) dx,
	\end{equation}
	with 
	\begin{equation}\label{eq:int}
		\bff G(x; \bmu)=e^{z(x){T}}\hat{\bff u}(z(x); \bmu)z'(x).
	\end{equation} 
	Since we are interested in approximating $\bff u({T};\bmu)$ with prescribed precision $tol$, we will only consider the portion of the Bromwich integral parameterized in a closed interval,
	\begin{equation}\label{eq:par}
		\bff u(T;\bmu) = \frac{1}{2\pi \opi}\int\limits_{\bR} \bff G(x; \bmu) dx \approx \frac{1}{2\pi \opi}\int\limits_{-c \pi}^{c \pi} \bff G(x;\bmu) dx,
	\end{equation}
	for a suitable positive truncation parameter $c \le c_{\text{max}}$, which we determine a priori by imposing the condition $\| \bff G(c \pi;\bmu )\| = tol$. The application of a quadrature formula to approximate \eqref{eq:par} provides a numerical approximation of $\bff u$, at the given time ${T}$ without need of
	computing it at intermediate time instants. For instance, an application of the trapezoidal rule with constant stepsize $\frac{2 c \pi }{N}$ provides the desired approximation
	\begin{align} 
		& \bff u_N(\bmu)=\frac{ c}{iN}\sum_{j=1}^{N-1}\bff G(\xi_j;\bmu),\label{eq:tr}\\
		\mbox{with } \ & \xi_j=-c \pi + j \frac{2 c \pi }{N},\quad j=1,\ldots,N-1.\label{eq:nodes}
	\end{align}
	
	Under standard regularity assumptions on $\bff G(x;\bmu)$ it can be shown, see Theorem $1$ in \cite{Guglielmi2018}, that the trapezoidal quadrature rule has an exponential rate of convergence with respect to the number of quadrature nodes $N$.
	
	The crucial point is the construction of the integration contour $\Gamma$. Assuming that the Laplace transform can be analytically extended to the left half of the complex plane and that its extension is properly bounded with respect to $z$, several authors have proposed different choices of the contour $\Gamma$. In the literature this approach has been widely studied for pure diffusion equations (see e.g. \cite{Gavrilyuk2001, LopezFernandez2004, LopezFernandez2006, Sheen2003}), and more recently for convection diffusion {reaction} equations in \cite{Guglielmi2018, Guglielmi2020} and for fractional time PDEs \cite{Colbrook2021, Colbrook2021a}. We refer the reader to \cite{Guglielmi2018} for a detailed review of the literature concerning the choice of the profile $\Gamma$.
	
	The magnitude of the resolvent norm $\|\left(z\bff \Id-\bff A(\bmu)\right)^{-1}\|$ (with $\|\cdot\|$ being the spectral norm when not indicated differently) has a crucial role in the rate of convergence of any contour integral method based on Laplace transform. Due to this, the choice and parametrization of the integration contour requires some information on the behaviour of the resolvent norm. In the recent work \cite{Guglielmi2020} elliptic, parabolic and hyperbolic profiles have been proposed, in connection to the knowledge of the so-called weighted $\varepsilon$-pseudospectrum of $\bff A$ (see \cite{Trefethen2005}), which is defined by
	\begin{equation*}  
		\sigma_{\varepsilon,t}(\bff A)  =  \left\{ z \in \bC \ : \  \e^{-\real{(z)}t}\sigma_{\min}(\bff A-z \bff \Id) \le \varepsilon\right\} 
	\end{equation*}
	for suitable $\varepsilon > 0$, with $\sigma_{\min}$ denoting the smallest singular value. Since $\bff A$ is in general non-normal the pseudospectrum may rapidly increase around the spectrum of $\bff A$, making the problem quite challenging.
    \vspace{-0.3cm}
	\subsection{Model reduction -- The Kolmogorov n-width}
	\label{sec:KnW}
	The aim of the present model reduction is to reduce the size $N_h$ of the discrete system to some $n\ll N_h$. This will be done by determining a suitable $n$-dimensional linear subspace in $\F^{N_h}$, $\F\in\{ \R,\C\}$,  formed by snapshots.\footnote{We leave the freedom for real and complex values in order to allow for reduced spaces in the original and the Laplace space.} In order to study the error of such a reduction, it is quite standard to consider the \emph{Kolmogorov} $n-$\emph{width} of some subset $\cM\subset\F^{N_h}$, namely
\begin{align*}
	d_n(\cM) &:= \inf_{Y_n\in\F^{N_h}\!;\, \dim(Y_n)=n}\,
		\sup_{\bff y\in\cM}
		\inf_{\bff y_n\in Y_n} \| \bff y-\bff y_n\|,
\end{align*} 
where the norm $\|\cdot\|$ is the Euclidean norm in $\F^{N_h}$. The specific focus here lies on the approximation of PPDEs, i.e., the subset $\cM$ is chosen as the set of solutions of the Cauchy problem or its Laplace transform. We will focus on the solution $\bu(t;\bmu)$ at a fixed time $t=T$ or a (small) time window, i.e., the \emph{solution set}
\begin{align}\label{def:SolSet}
	\cM(I;\cD) \!:=\! \{ \bu(t;\bmu)\in\R^{N_h}:\,
		t\in I, \bmu\in\cD
		\text{ and } \bu(t;\bmu) \text{ solves } \eqref{eq:cauchy}\,
		 \forall t\in I\},
\end{align}
where $I\subseteq (0,T]$ is some time window, which might also consist of only one instant, i.e., $I=\{T\}$, in which case we shall simply write $\cM(T;\cD)$.

 If the $n$-width decays rapidly as $n$ increases, for instance $d_n(\mathcal{M}) = {\mathcal{O}}(e^{-n})$, it means that the solution set can be well approximated by a space of {relatively small} dimension $n$, that is called \textit{reduced space}. The fact that the Kolmogorov $n$-width decays fast is problem-specific, i.e., not all problems exhibit this behaviour. Up to our knowledge there does not exist a general theory which, for a given a PPDE and the associated solution set, is able to determine the behaviour of its Kolmogorov $n$-width. However, there are rigorous results for certain classes of continuous problems. In \cite{Cohen2015} it is shown that when the Kolmogorov $n$-width of a compact set $K$ of $X$ is $\mathcal{O}(n^{-r})$ for some $r>1$, then those of $\mathcal{U}(K)$ are $\mathcal{O}(n^{-s})$ for any $s<r-1$ under the hypothesis of holomorphicity of $\mathcal{U}$. This result is used in \cite{Ohlberger2015} to prove exponential decay of the Kolmogorov $n$-width for elliptic operators that admits an affine decomposition. In the same work a polynomial rate of decay is proved for the one dimensional advection equation with discontinuous initial data. In \cite{Greif2019} it is shown that $d_n(\mathcal{M})\ge\mathcal{O}(n^{-1/2})$ for $1$D-hyperbolic wave equation with discontinuous initial data.
	
One contribution of this article is a result on the Kolmogorov $n$-width behaviour for problems of type (\ref{eq:cauchy}), which is given {in Theorem \ref{teo1}. Before stating the theorem we provide a lemma which directly comes from the result in \cite[Theorem $3.1$]{Ohlberger2015}}. Let us define the map
	\begin{equation}\label{eq:map_Lap2}
		\Phi_{\mathcal{L}}:\;\;{\eta=(z,\bmu)}\in\Gamma_r\times\mathcal{D}\longrightarrow \hat{\bff u}(\eta),
	\end{equation}
	where 
	$\hat{\bff u}(\eta)$ is solution of
	\begin{equation}\label{eq:pde_l2}
		P(\hat{\bff u}(\eta), \eta)=(z\bff \Id-\bff A (\bmu))\hat{\bff u}(\eta)-\bff u_0(\bmu)-\hat{\bff b}(\eta)=0,
	\end{equation}
	and the solution set
	\begin{equation}\label{eq:solmanLap}
		\hat{\mathcal{M}}(\mathcal{D})=\Big\{\hat{\bff u}(z;\bmu)|\;\bmu\in\mathcal{D},\; {z\in \Gamma_r} \text{ and } \hat{\bff u} {(\eta)}\text{ solves } \eqref{eq:pde_l2} \text{ for } \eta\in \Gamma_r\times\mathcal{D} \Big\}.
	\end{equation}
	Now we can provide the following result.
    \begin{lemma}\label{lemma1}
    If the affine decomposition \eqref{eq:affine2}, \eqref{eq:affinLT}  hold, then
    \begin{equation*}
    		d_n(\hat{\mathcal{M}})\le k_1 e^{-k_2n^{1/Q}}, \quad Q=Q_{\bff A}+{Q}_{\hat{\bff b}}+1.
    \end{equation*}
    with $k_1>0$ and $k_2>0$ constants independent from $n$.
    \end{lemma}
\begin{proof}
	The proof is analogous to the one of \cite[Theorem $3.1$]{Ohlberger2015} and is based on an application of \cite[Theorem $4.1$]{Cohen2015}. The only difference consists on the fact that we also consider a parametric dependence of the forcing term $\hat{b}$ which, under our hypotheses, does not affect the validity of the argument used in \cite[Theorem $3.1$]{Ohlberger2015}.
\end{proof}
 We denote by $V^{*}_n$ the $n$-dimensional subspace such that
	\begin{equation*}
	    d_n(\hat{\mathcal{M}})\le \sup_{(z,\bmu)\in\Gamma_{r}\times\mathcal{D}}\inf_{y_n\in V^{*}_n}\|\hat{\bff u}(z;\bmu)-\bff y_n\|\le k_1 e^{-k_2n^{1/Q}},
	\end{equation*}
	that in general is not the optimal $n$-dimensional subspace but that it is suitable to show the exponential rate of decay of the Kolmogorov $n$-width $ d_n(\hat{\mathcal{M}})$.
	\begin{samepage}
	\begin{theorem} \label{teo1}
	Let the affine decomposition \eqref{eq:affine2}, \eqref{eq:affinLT} hold and assume:
            {\begin{enumerate}[(i)]
	        \item $\Gamma$ is an analytic curve that keeps on its left all the singularities of $\hat{\bff u}(z;\bmu)$, $\forall \bmu\in\mathcal{D}$;
			\item the Laplace transform of the source term for $\real(z)\rightarrow -\infty$ and $z\in \Gamma$ grows at most linearly, i.e.
	\begin{equation}\label{eq:boundsLP}
	    \lim_{\real(z)\rightarrow -\infty}\|\hat{\bff b}(z;\bmu)\|\le M_1 z,\quad\forall\bmu\in\mathcal{D}.
	\end{equation}
		\end{enumerate}}
		Then it holds that the Kolmogorov $n$-width associated to the solution set $\mathcal{M}(T;\mathcal{D})$ has a exponential decay that is, for suitable constants $\tilde{k}_1$ and $\tilde{k}_2$ independent of
		$n$, it holds
		\begin{equation*} 
            d_n(\mathcal{M}(T;\mathcal{D}))\le 
            \tilde{k}_1 e^{-\tilde{k}_2 n^{1/Q}},\quad Q=Q_{\bff A}+Q_{\hat{\bff b}}+1.
		\end{equation*}
	\end{theorem}
	\end{samepage}
Before proceeding with the proof we comment the result. The required hypothesis holds for discretized linear parabolic problems where the exponential decay of the Kolmogorov $n$-width is notoriously expected. Moreover, the well-known fact that the velocity of decay is influenced by the number of considered parameters in the problem and is expressed by the power $n^{1/Q}$ in the exponential term. What about those problems where the Kolmogorov $n$-width has an algebraic decay? 
Let us consider the $1$D linear advection problem discussed in \cite{Ohlberger2015}. The solution of this problem is $H(\bmu t-x)$ where $\bmu\in[0,1]$ is the parametric velocity, $t\in[0,1]$ stands for time and $x\in[0,1]$ is the space variable. The Laplace transform of this function is $\mathcal{L}(H(\bmu t-x))=(1/(z\bmu))(e^{-xz/\bmu})$ and if we insert it in \eqref{eq:boundsLP} we find that the condition is not satisfied. The Laplace transform of the Heaviside function, in the discrete problem, appears in the known term due to the boundary conditions on the extrema of the space domain.
	
	\begin{proof}
		Applying the Laplace transform to (\ref{eq:cauchy}) we get
		\eqref{eq:LT}
{ and by assumptions $(i)$ and $(ii)$, see \cite{TALBOT1979, Trefethen2006}, we obtain the solution via inverse Laplace transform \eqref{eq:solu}, in particular the first assumption implies that there exists $M_2\in\mathbb{R}^{+}$ such that
\begin{equation}\label{eq:bound A and b}
	\|(z\bff \Id-\bff A (\bmu))^{-1}\|\le \frac{M_2}{z}
\end{equation}
for all $\bmu\in\mathcal{D}$ and $z\in\Gamma$, which implies that $\|\hat{\bff u}(z;\bmu)\|$ is bounded $\forall z\in \Gamma$.} 

We consider the parametrization of the parabolic contour $\Gamma$
\begin{equation}\label{eq:map_p}
 z(s)\;:\;\mathbb{R}\rightarrow\Gamma, \quad z(s)\;:\;-s^2-2i s a_1+a_2    
\end{equation}
with real constants $a_1<0$ and $a_2>0$.

Due to the symmetry of the integration profile and of the integrand function, the integration can be restricted to the upper half complex plane $\{z:\;\imag(z)\ge 0\}$. The use of the trapezoidal rule with constant stepsize $\delta$ gives the approximation of \eqref{eq:solu}:
\begin{equation*}
	\bff u_\delta(T;\bmu)=2\delta\imag\left(\sum_{j=0}^{\infty}e^{z(s_j) T}\hat{\bff u}(z(s_j);\bmu)z'(s_j)\right),\;s_j=j\delta\pi
\end{equation*}
where, by assumption ($ii$), that is \eqref{eq:boundsLP}, making use of \cite[Theorem $5.1$]{Trefethen2014TheEC}, we obtain
$\forall\bmu\in\mathcal{D}$:
\begin{equation}\label{eq:errquad}
	\|\bff u(T;\bmu)-\bff u_\delta(T;\bmu)\|\le C(\Omega)\frac{2{Q}}{\e^{2\pi {v}/\delta}-1} 
	\le c_1 e^{-c_2/\delta},
\end{equation}
where $Q$ and $v$ are two constants, $\delta$ can be made arbitrarily small, $ C({\Omega})$ is a constant that depends from the measure of the domain $\Omega$ and $c_1$, $c_2$ are
constants independent of $\delta$.  

{
Let us define $\bff u_{\delta,M}(T;\bmu)$ as 
\begin{equation}\label{eq:udeltaM}
    \bff u_{\delta,M}(T;\bmu)=2\delta\imag\left(\sum_{j=0}^{M}e^{z(s_j)T}\hat{\bff u}(z(s_j);\bmu)z'(s_j)\right),
\end{equation}
By assumptions $(i)-(ii)$ we have that 
\begin{equation}\label{eq:truncation}
    \|\bff u_{\delta}(T;\bmu)-\bff u_{\delta,M}(T;\bmu)\|
    \le c_0 e^{-M^2} \qquad \forall\; \bmu\in\mathcal{D},
\end{equation}
with the constant $c_0$ independent of $M$.}

We let $\real( z(s_M) ) = r$ and define
\begin{equation*}
\Gamma_r = \Gamma \cap \{z \in \mathbb{C}\; |\; \real(z) \ge r \}.
\end{equation*}
Then recalling the definition of Kolmogorov $n$-width we have
\begin{align*}\label{eq:bound}
	&d_n(\mathcal{M}(T;\mathcal{D}))=\inf_{Y_n\subset\mathbb{C}^{N_h};\text{dim}(Y_n)=n}\sup_{\bmu\in\mathcal{D}}\inf_{\bff y_n\in Y_n}\|\bff u(T;\bmu)-\bff y_n\|\nonumber\\
	=&\inf_{Y_n\subset\mathbb{C}^{N_h};\text{dim}(Y_n)=n}\sup_{\bmu\in\mathcal{D}}\inf_{\bff y_n\in Y_n}\Bigg\| \bff u(T;\bmu)- \bff u_{\delta}(T;\bmu)\nonumber\\
	+&\bff u_{\delta}(T;\bmu)-\bff u_{\delta,M}(T;\bmu)+\bff u_{\delta,M}(T;\bmu)-\bff  y_n  \Bigg\|    \nonumber\\
	\le &c_1e^{-c_2/\delta}+c_0 e^{-M^2} +   \inf_{Y_n\subset\mathbb{C}^{N_h};\text{dim}(Y_n)=n}\sup_{\bmu\in\mathcal{D}}\inf_{\bff y_n\in Y_n}\Bigg\|\bff u_{\delta,M}(T;\bmu)-\bff y_n \Bigg\|.
\end{align*}
where we used \eqref{eq:errquad} and \eqref{eq:truncation}. Note that $\bff u_{\delta,M}(T;\bmu)$ as defined in (\ref{eq:udeltaM}) is a linear combination of $\hat{\bff u}(z(s_j);\bmu)$ for $j=0,...,M$ with $z(s_j)\in\Gamma_{r}$, $\forall j$. By Lemma \ref{lemma1} we have that $\hat{\bff u}(z(s_j);\bmu)$ belongs to the manifold (\ref{eq:solmanLap}) that has an exponential decaying Kolmogorov $n$-width, therefore if we choose $Y_n$ as the same subspace chosen for Lemma \ref{lemma1}, that was denoted as $V^{*}_n$, we get
\begin{align*}
	&\inf_{Y_n\subset\mathbb{C}^{N_h};\text{dim}(Y_n)=n}\sup_{\bmu\in\mathcal{D}}\inf_{\bff y_n\in Y_n}\Bigg\| \bff u_{\delta,M}(T;\bmu)-\bff y_n\Bigg\|\\
	\le&\sup_{\bmu\in\mathcal{D}}\inf_{\bff y_n\in V^{*}_n}\Bigg\| \bff u_{\delta,M}(T;\bmu)-\bff y_n\Bigg\|\le \delta 
	k_1c_3e^{\real(z(s_0)) T}(M+1)e^{-k_2n^{1/Q}},
\end{align*}
Now choosing $\delta=c_2(k_2n^{1/Q})^{-1}$ and
$M = (c_2/\delta)^{1/2}$ we get
\begin{align*}
	d_n(\mathcal{M}(T;\mathcal{D}))\le& c_1e^{-k_2n^{1/Q}}+c_0 e^{-k_2n^{1/Q}}+ c_2(k_2n^{1/Q})^{-1} c_3e^{\real (z(s_0))T}(M+1)e^{-k_2n^{1/Q}}\nonumber\\
	\le&\left(c_1+c_0+c_2c_3(k_2n^{1/Q})^{-1/2}e^{c_5 T}\right)  e^{-k_2n^{1/Q}}\\
	\le&\left(c_1+c_0+c_2c_3(k_2)^{-1/2}e^{c_5 T}\right)  e^{-k_2n^{1/Q}}
	\nonumber
\end{align*}
and setting $\tilde{k}_1 = c_0+c_1+c_2 c_3(k_2)^{-1/2}e^{c_5 T}, \tilde{k}_2 = c_4$, the previous finally yields
\[
d_n(\mathcal{M}(T;\mathcal{D}))  \le \tilde{k}_1 e^{-\tilde{k}_2 n^{1/Q}}.
\]
\end{proof}
We remark that in the theorem we have chosen, for convenience, a parabolic profile; in principle we could choose any type of smooth curve satisfying the first hypothesis of the theorem. The outline of the proof would be the same but the final polynomial rate multiplying the exponential would change exponent according to the behaviour of the real part of the chosen profile of integration. {Also the value of the constants would be modified}.

Another remark consists in the fact that our proof holds also if we consider the time window $t\in(0,T]$ rather than the fixed time $t=T$. This is due to the fact that the dependence from time only appears on the scalar quantities $e^{\real(z)t}$ therefore the same subspace $V^{*}_n$ defined in the proof of the theorem can approximate $\bff u_{\delta,M}$ with the an uniform accuracy in all the time window $(0,T]$. 
\begin{corollary}
Under the same hypothesis of Theorem \ref{teo1} the Kolmogorov $n$-width associated to the solution set
\begin{equation*}
    \mathcal{M}(I;\mathcal{D})=\cup_{t\in(0,T]}\mathcal{M}_h(t;\mathcal{D}),\quad I=(0,T]
\end{equation*}
has a exponential decay that is, for suitable constants $\tilde{k}_1$ and $\tilde{k}_2$ independent of $n$, it holds
\begin{equation*}
    d_n(\mathcal{M}(I;\mathcal{D}))\le 
            \tilde{k}_1 e^{-\tilde{k}_2 n^{1/Q}},\quad Q=Q_{\bff A}+{Q}_{\hat{\bff b}}+1.
\end{equation*}
\end{corollary}
Note that, since we work in finite machine precision, the time window that is possible to consider in the applications with only one profile $\Gamma$ is $[t_0,\Lambda t_0]$, $\Lambda>1$, always smaller than $(0,T]$.
	
	\subsection{The reduced problem}
	We aim to apply projection model order reduction to problem (\ref{eq:cauchy}) where the time integration is performed through the contour integral method proposed in \cite{Guglielmi2020}.
	
	First, we observe that the numerical approximation given by (\ref{eq:tr}) is a linear combination of the function $\bff G(x;\bmu)$ evaluated at the quadrature points $\xi_j$ defined in (\ref{eq:nodes}). Therefore, recalling (\ref{eq:int}), the result of Lemma \ref{lemma1} and the fact that $\hat{\bff u}(z(\xi_j);\bmu)\in\mathbb{C}^{N_h}$ for all $j$ and $\bmu$, we define the solution set of our problem as the restriction of \eqref{eq:solmanLap} on the quadrature points, i.e.
	\begin{equation}\label{eq:manLap}
		\hat{\mathcal{M}}_r(\mathcal{D})=\Big\{\hat{\bff u}(z(\xi_j);\bmu)\;|\;\bmu\in\mathcal{D},\;j=1,...,N-1\Big\}\subset\hat{\mathcal{M}}(\mathcal{D}).
	\end{equation}
	In Section \ref{sec:Greedy} we will discuss how to build an appropriate low dimensional subspace of {$\mathbb{C}^{N_h}$ that well approximate (\ref{eq:manLap})} but for the moment we simply assume that such a subspace has been constructed. We can uniquely represent the reduced space by the $N_r$ (usually much smaller than $N_h$) base vectors $\zeta_i\in\mathbb{C}^{N_h}$, with $i=1,...,N_r$. 
	
	The $N-1$ (see \eqref{eq:nodes}) reduced solutions are taken as linear combinations of such base vectors and are determined by the Galerkin projection of (\ref{eq:LT}) onto the reduced space, that means
	\begin{equation*}
	\vspace{-0.5cm}
		\hat{\bff u}_r(z_j; \bmu)=\sum_{i=1}^{N_r}\beta_{j,i}(\bmu)\bff \zeta_i
	\end{equation*}
	satisfying
	\begin{equation}\label{eq:galproj}
		\langle (z_j\bff \Id-\bff A(\bmu))\hat{\bff u}_r(z_j;\bmu),\;\bff \zeta_i\rangle=\langle\bff u_0(\bmu)+\hat{\bff b}(z_j;\bmu),\;\bff \zeta_i\rangle\quad\text{for}\;i=1,...,N_r\;
	\end{equation}
	where $\langle\cdot,\cdot\rangle$ denotes the scalar product on $\mathbb{C}^{N_h}$ and $z_j=z(\xi_j)$. Thanks to Galerkin orthogonality we can make use of Céa's Lemma to connect the best approximation error with the reduced approximation:
	\begin{align}\label{eq:Cea}
		\max_{1\le j \le N-1}\Big\|\hat{\bff u}(z_j;\bmu)-\hat{\bff u}_{r}(z_j;\bmu)\Big\|\le&\max_{1\le j \le N-1}\left(\left(1+\frac{\gamma_j(\bmu)}{\alpha_j(\bmu)}\right)\inf_{\bff v\in\mathbb{C}^{N_r}}\Big\|\hat{\bff u}(z_j;\bmu)-\bff v\Big\|\right)\nonumber\\
		\le&d_{N_r}(\hat{\mathcal{M}}_r)\max_{1\le j \le N-1}\left(1+\frac{\gamma_j(\bmu)}{\alpha_j(\bmu)}\right);
	\end{align}
	where $\alpha_j(\bmu)$ and $\gamma_j(\bmu)$ are the coercvity and continuity constants of $z_j \bff \Id-\bff A(\bmu)$, that means
	\begin{align*}
		\alpha_j(\bmu)&=\inf\{|\lambda|:\lambda\in W(z_j\bff \Id-\bff A)\}=\inf_{\bff w \in \mathbb{C}^{N_h}, \norm{\bff w}=1}{|\bff w^*(z_j\bff \Id-\bff A(\bmu))\bff w|};\\
		\gamma_j(\bmu)&=\sup_{\bff w \in \mathbb{C}^{N_h}, \norm{\bff w}=1}\sup_{\bff y \in \mathbb{C}^{N_h}, \norm{\bff y}=1}{\bff y^*(z_j\bff \Id-\bff A(\bmu))\bff w}=\norm{z_j\bff \Id-\bff A(\bmu)}_{2},
	\end{align*}
	where $W(z_j\bff \Id-\bff A)$ is the so called numerical range or field of values of the complex matrix $z_j\bff \Id-\bff A$. Thus, the quality of the reduced order approximation, based on a Galerkin projection, depends, as is well known, on the coercivity and continuity constants of the operators $z_j\bff \Id-\bff A(\bmu)$, which are problem-dependent. Furthermore, it also depends on the accuracy of the reduced space in approximating the entire solution set, which is related to the Kolmogorov $n$-width. 
	Note that, for those operators which do not satisfy the coercivity assumption (i.e. $\alpha_j(\bmu)>0,\;\forall \bmu\in\mathcal{D}$), the {projection model order reduction} framework can still be effectively applied by making use of the inf-sup stability constant which replaces the coercivity constant (see \cite[Chapter $6$]{Hesthaven2016}).
	
	Once $\hat{\bff u}_r(z_j)$ is computed for all $j=1,...,N-1$ we can approximate the solution of the reduced problem {at time $t=T$} as
	\begin{equation}\label{eq:redsol}
		\bff u_r({T};\bmu)=\frac{ c}{iN}\sum_{j=1}^{N-1}e^{z_j {T}}\hat{\bff u}_{r}(z_j; \bmu)z'_j\;,
	\end{equation}
	where $z'_j=z'(\xi_j)$. Its reliability with respect to the high fidelity solution $\bff u({T};\bmu)$ is determined by (\ref{eq:Cea}), in fact \vspace{-0.2cm}
	\begin{align*}
		\Big\|\bff u({T};\bmu)-\bff u_{r}({T};\bmu)\Big\|=&\Bigg\|\frac{ c}{iN}\sum_{j=1}^{N-1}e^{z_j{T}}\Big(\hat{\bff u}(z_j; \bmu)-\hat{\bff u}_{r}(z_j; \bmu)\Big)z'_j\Bigg\|\nonumber\\
		\le&\frac{c}{N}\sum_{j=1}^{N-1}\Bigg(e^{\real{             ({z_j})}{T}}\Big|z'_j\Big|\Big\|\hat{\bff u}(z_j;\bmu)-\hat{\bff u}_{r}(z_j;\bmu)\Big\|\Bigg)\nonumber\\
		\le&d_{N_r}(\mathcal{M}(T))\frac{c}{N}\sum_{j=1}^{N-1}\Bigg(e^{\real{             ({z_j})}{T}}\Big|z'_j\Big|\left(1+\frac{\gamma_j(\bmu)}{\alpha_j(\bmu)}\right)\Bigg).
	\end{align*}
	Note that the reduced space lies in $\mathbb{C}^{N_r}$ whose dimension is $2N_r$. However, whenever we compare our reduced order method with the ones constructed with real snapshots, we directly confront it with spaces that lie in $\mathbb{R}^{N_r}$ rather than $\mathbb{R}^{2N_r}$. Having larger spaces reflects in a better accuracy, but it increases the computational cost. In our comparisons with real reduced spaces we have both these features.
	
	{For simplicity we have presented our model order reduction contour integral method formulation for a fixed time instant $t=T$. However, everything can be extended to a time window of (usually) moderate size $[t_0, \Lambda t_0]$, with $\Lambda>1$, without increasing the computational cost, see \cite{Guglielmi2018, Guglielmi2020} for details. The main idea relies on the fact that the dependence from time of the high fidelity solution appears only in the scalar quantity $e^{\real(z)t}$, therefore the reduced space constructed for time $t=T$ can also approximate the solution in a suitable  time window.}
	
	\subsection{Offline-Online phases} \label{sec:On-Off}
	
	According to the literature, the offline phase is the one where the reduced space is constructed. This phase is the one computationally demanding since building the reduced space requires the solution of several full-size problems. The online phase instead is the one where the reduced solution $\bff u_r(t;\bmu)$, (see (\ref{eq:redsol})), is computed. In an ideal setting, the computational cost of the online phase should be independent of the complexity of the full problem (which depends on $N_h$) and should depend only on the size $N_r$ of the reduced problem. To compute $\bff u_r(t;\bmu)$ we need to solve (\ref{eq:galproj}), that is to calculate 
	\begin{equation*}
		{\bff \beta}_j(\bmu)=(z_j\bff \Id_{N_r}-\bff A_r(\bmu))^{-1}(\bff u^r_{0}(\bmu)+\hat{\bff b}_r(z_j;\bmu)),\quad \text{for}\;j=1,...,N-1,
	\end{equation*}
	being $\bff \Id_{N_r}$ the $N_r\times N_r$ identity matrix,
	\begin{align}
		&\bff A_r(\bmu)=\bff B^H\bff A(\bmu)\bff B,\quad\hat{\bff b}_r(z_j;\bmu)=\bff B^H\hat{\bff b}(z_j;\bmu),\quad \bff u^r_{0}(\bmu)=\bff B^H\bff u_0(\bmu). \label{eq:redA}
	\end{align}
	with $\bff B$ the $N_h\times N_r$ matrix whose columns are given by the reduced space basis vectors $\bff \zeta_j$, $j=1,...,N_r$. This requires to form (\ref{eq:redA}) for each instance of $\bmu$, whose cost still depends on $N_h$ (since $\bff A(\bmu)\in\mathbb{R}^{N_h\times N_h}$, $\bff B\in\mathbb{C}^{N_h\times N_r}$ and $\hat{\bff b}(z;\bmu)\in\mathbb{C}^{N_h}$). However, this is simplified by employing assumption (\ref{eq:affine2}) and \eqref{eq:affinLT}. Therefore the $Q_{\bff A}$ square matrices $\bff A^r_{q}$ of dimension $N_r$ (associated to $\bff A_q$), can be precomputed in the offline phase once the reduced space is known. This computation is done only once and may be run in parallel since $\bff A_q$ are independent from $\bmu$. Then, during the online phase, for a new instance of the parameter $\bmu$, the solution matrix is built as 
	\begin{equation*}
		\bff A_r(\bmu)=\sum_{q=1}^{Q_{\bff A}}\vartheta_{q}^{A}({\bmu}){\bff A}_q^{r}.
	\end{equation*}
	This operation is independent of $N_h$ and scales proportionally to $Q_{\bff A}\cdot N_r$. The treatment of the known term $\bff u_0(\bmu)+\hat{\bff b}(z,\bmu)$ is analogous. {Algorithm \ref{al3} resumes the online phase procedure of our method. First, for a given instance of $\bmu\in\mathcal{D}$ and of $t$ belonging to a suitable time window $[t_0,\Lambda t_0]$, the reduced matrix $\bff A_r(\bmu)$ and the reduced known term $\bff u^r_0(\bmu)+\hat{\bff b}_r(z,\bmu)$ are assembled employing the affine decomposition property, than the vectors $\hat{\bff u}_r(z_j;\bmu)$ are computed on the quadrature points and, finally, the reduced solution $\bff u_r(t;\bmu)$ is assembled as linear combination of $\{\hat{\bff u}_r(z_j;\bmu)\}_{j=1}^{N-1}$. Note that the number of quadrature points can be halved due to symmetry of the integrand function and of the integration profile, moreover the evaluations of $\hat{\bff u}_r(z_j;\bmu)$ for $j=1,...,N-1$ can be done in parallel.} 
	
	For cases where an affine decomposition is not available, one can often find an approximate form that satisfies this property through empirical interpolation (see e.g. \cite{Barrault2004}) or empirical cubature methods (see e.g. \cite{Farhat2015,Yano2019,Manucci2022}). However, for the purposes of this work {that concerns linear parabolic problems}, we assume that the affine decomposition is available.
	\begin{small}
	\begin{algorithm}[t]
		\caption{{Online phase MOR contour integral method}}\label{al3}
		\hspace*{\algorithmicindent} \textbf{Input:} $\bmu\in\mathcal{D}$, $t\in[t_0,\Lambda t_0]$
		
		\hspace*{\algorithmicindent} \textbf{Output:} The reduced solution $\bff u_r(t;\bmu)$
		\begin{algorithmic}[1]
		    \STATE Compute $\bff A_r(\bmu)=\sum_{q=1}^{Q_{\bff A}}\vartheta_{q}^{\bff A}({\bmu}){\bff A}_q^{r}$,  $\bff u^r_{0}(\bmu)+\hat{\bff b}_r(z,\bmu)=\sum_{q=1}^{Q_{\bff b}}\vartheta^{ \bff b}_{q}({\bmu}){ \hat{\bff b}}_q^{r}(z)$
		    \FOR {$j=1,...,N-1$} 
		    \STATE Compute ${\bff \beta}_{j}(\bmu)=\left(z_j\bff \Id_{N_r}-\bff A_r(\bmu)\right) \setminus (\bff u^r_{0}(\bmu)+\hat{\bff b}_r(z,\bmu)))$
		    \ENDFOR
			\STATE Compute $\bff u_r(t;\bmu)=\frac{ c}{iN}\bff B\left(\sum_{j=1}^{N-1}e^{z_j t}{\bff \beta}_{j}(\bmu)z'_j\right)$ 
		\end{algorithmic}
	\end{algorithm}
	\end{small}
	\subsection{Sensitivity of the projection reduced order methods}
	As we have seen the construction of the reduced space starts from the solution of full-size problems; this process is affected by errors and therefore one should control that the resulting reduced space (and thus the perturbed reduced problem) still guarantees a reliable approximation of the full problem.
	
	Given an integration contour and given a set of quadrature points we adopt the following notation:
	\begin{align*}
		&\bff A_j(\bmu)=z_j\bff \Id-\bff A(\bmu), \quad \bff g_j(\bmu)=\bff u_0(\bmu)+\hat{\bff b}(z_j;\bmu),\nonumber\\
		&\bff x_j(\bmu)=\bff A_j^{-1} \bff g_j(\bmu), \quad \tilde{\bff x}_j=\bff A_j^{-1}\tilde{\bff g}_j(\bmu);  
	\end{align*}
	where $\tilde{\bff g}_j(\bmu)$ is a perturbation of $\bff g_j(\bmu)$. We construct the matrix $\tilde{\bff B}$ of dimension $N_r$ collecting snapshots of the perturbed full problems
	\begin{equation*}
		\bff A_j(\bmu)\tilde{\bff x}_j=\tilde{\bff g}_j(\bmu),\quad\text{for } j=1,...,N-1,\quad\text{for some }\bmu\in\mathcal{D};
	\end{equation*} 
	we then have the perturbed reduced systems\vspace{-0.2cm}
	\begin{equation*}
		\tilde{\bff B}^H\bff A_j(\bmu)\tilde{\bff x}^r_j(\bmu)=\tilde{\bff B}^H\tilde{\bff g}_j(\bmu),\quad\forall j=1,...,N-1,
	\end{equation*} 
	where $\tilde{\bff x}^r_j(\bmu)$ is the solution of the perturbed reduced system. Then we define the perturbed solution set as
	\begin{equation*}
		\tilde{{\mathcal{M}}}(\mathcal{D}):=\Big\{\tilde{\bff x}_j(\bmu)\;|\;\bmu\in\mathcal{D},\;j\in\{1,...,N-1\}\Big\},
	\end{equation*}
	and the associated Kolmogorov $n$-width $d_n(\tilde{\mathcal{M}}(\mathcal{D}))$.
	
	The next Proposition \ref{Prop:1} provides an upper bound for the error
	\begin{equation*}
		e(t,\bmu)=\|\bff u(t,\bmu)-\tilde{\bff u}_r(t,\bmu)\|,
	\end{equation*}
	where $\bff u(t,\bmu)$ is given by (\ref{eq:tr}) and $\tilde{\bff u}_r(t,\bmu)$ is similarly given by
	\begin{equation*}
		\tilde{\bff u}_r(t;\bmu)=\frac{ c}{iN}\sum_{j=1}^{N-1}e^{z_jt}\tilde{\bff x}^r_j(\bmu)z'_j, \quad \mbox{for } \ j=1,...,N-1.
	\end{equation*}
	\vspace{-0.5cm}
	\begin{proposition}\label{Prop:1}
		For any $\bmu\in\mathcal{D}$, $t>0$ and $N_r\in\mathbb{N}$ it holds
		\begin{align*}
			e(t,\bmu)\le\frac{c}{N}\sum_{j=1}^{N-1}\Bigg(\e^{\real({z_j})t}\Big|z'_j\Big|\Bigg(\omega_j(\bmu)+\left(1+\frac{\gamma_j(\bmu)}{\alpha_j(\bmu)}\right)\Big(q+d_{N_r}(\hat{\mathcal{M}}(\mathcal{D}))\Big)\Bigg)\Bigg)
		\end{align*} 
		where
		\begin{equation*}
		\omega_j(\bmu)=\norm{\bff A_j^{-1}(\bmu)}\norm{\bff g_j(\bmu)-\tilde{\bff g}_j(\bmu)}\quad\text{and}\quad q=\sup_{\bmu\in\mathcal{D}}\max_{1\le j\le N-1}\omega_j(\bmu).
		\end{equation*}

	\end{proposition}
	\begin{proof}
	We first make use of the triangular inequality and then Céa's Lemma for the perturbed problem
		\begin{align*}
			e(t,\bmu)=&\norm{\bff u(t,\bmu)-\tilde{\bff u}_r(t,\bmu)}\le\frac{c}{N}\sum_{j=1}^{N-1}\Bigg(e^{\real({z_j})t}\Big|z'_j\Big|\norm{{\bff x}_j(\bmu)-\tilde{\bff x}^{r}_{j}(\bmu)}\Bigg)\nonumber\\
			\le&\frac{c}{N}\sum_{j=1}^{N-1}\Bigg(e^{\real             ({z_j})t}\Big|z'_j\Big|\Big(\norm{{\bff x}_j(\bmu)-\tilde{\bff x}_{j}(\bmu)}+\norm{\tilde{\bff x}_{j}(\bmu)-\tilde{\bff x}^{r}_{j}(\bmu)}\Big)\Bigg)\nonumber\\
			\le&\frac{c}{N}\sum_{j=1}^{N-1}\Bigg(e^{\real             ({z_j})t}\Big|z'_j\Big|\left(\omega_j(\bmu)+\left(1+\frac{\gamma_j(\bmu)}{\alpha_j(\bmu)}\right)\inf_{\bff v\in \text{Range}(\tilde{\bff B})}\Big\|\tilde{\bff x}_j(\bmu)-\bff v\Big\|\right)\Bigg)\nonumber\\
			\le&\frac{c}{N}\sum_{j=1}^{N-1}\Bigg(e^{\real            ({z_j})t}\Big|z'_j\Big|\left(\omega_j(\bmu)+\left(1+\frac{\gamma_j(\bmu)}{\alpha_j(\bmu)}\right)d_{N_r}(\tilde{\mathcal{M}}(\mathcal{D}))\right)\Bigg),
		\end{align*}
		where we recall that $\tilde{B}$ is the matrix whose columns form a basis for the reduced space generated from the perturbed problem. We now need to bound the Kolmogorov $N_r$-width of the perturbed problem. We have
		\begin{align*}
			d_{N_r}(\tilde{\mathcal{M}}(\mathcal{D}))=&\inf_{\text{dim}(Y)=N_r}\sup_{\tilde{\bff x}\in\tilde{\mathcal{M}}}\inf_{{\bff v}\in Y}\norm{\tilde{\bff x}-\bff v} \le \inf_{\text{dim}(Y)=N_r}\sup_{\tilde{\bff x}\in\tilde{\mathcal{M}}}\inf_{{\bff v}\in Y}\norm{\tilde{\bff x}-\bff x+\bff x-\bff v}\nonumber\\
			\le&\inf_{\text{dim}(Y)=N_r}\sup_{\tilde{\bff x}\in\tilde{\mathcal{M}}}\inf_{{\bff v}\in Y}\left(\norm{\tilde{\bff x}-\bff x}+\norm{\bff x-\bff v}\right)\nonumber\\
			\le&\sup_{\tilde{\bff x}\in\tilde{\mathcal{M}}}\norm{\tilde{\bff x}-\bff x}+d_{N_r}(\hat{\mathcal{M}}(\mathcal{D}))\nonumber\\
			\le&\sup_{\bmu\in\mathcal{D}}\max_{1\le j\le N-1}\norm{\bff A_j^{-1}(\bmu)}\norm{\bff g_j(\bmu)-\tilde{\bff g}_j(\bmu)}+d_{N_r}(\hat{\mathcal{M}}(\mathcal{D})).
		\end{align*}
	\end{proof}
	As a consequence, up to 
	\begin{equation*}
		\sup_{\bmu\in\mathcal{D}}\max_{1\le j\le N-1}\norm{\bff A_j^{-1}(\bmu)}\norm{\bff g_j(\bmu)-\tilde{\bff g}_j(\bmu)},
	\end{equation*}
	the decay of $d_{N_r}(\tilde{\mathcal{M}}(\mathcal{D}))$ with respect to $N_r$ is at least as fast as the one of $d_{N_r}(\hat{\mathcal{M}}(\mathcal{D}))$.

	Proposition \ref{Prop:1} states that the difference between the solution of the full problem and the solution of the perturbed reduced problem depends on two factors: the resolvent norm, which we control in the construction of the integration profile, and the Kolmogorov $n$-width associated to the full problem, which we assume to know a priori.
	
	\section{Reduced space generation} \label{sec:Greedy}
	The most popular strategies to generate reduced spaces we consider are the proper orthogonal decomposition (POD) \cite{Sirovich1987}, essentially a truncated singular value decomposition of the snapshots matrix, and the greedy algorithm \cite{Prudhomme2001}. These two are also combined when dealing with time dependent problems.
	
	In this section we first introduce an error estimator based on the discrete formulation of the problem we consider; then we use it to define two novel greedy algorithms which provide suitable reduced {spaces} to solve the time integration with contour integral methods. An essential ingredient of these methods is a reliable error estimator.\vspace{-0.3cm}
	\subsection{Greedy basis generation}
	The standard greedy generation of the reduced space is an iterative procedure where at each iteration one new basis function is added and the overall approximation capability of the basis set is possibly improved. It requires the computation of one high fidelity solution per iteration and a total of $N_r$ linearly independent full-size solution vectors to generate the $N_r$-dimensional reduced space. One key requirement of the greedy algorithm is the availability of an error estimate $\Delta(\bmu)$ which predicts the error due to the model order reduction, that is the error induced by replacing the full space of dimension $N_h$ by the reduced space of dimension $N_r$. 
	
	This approach is also denoted as weak-greedy algorithm, to distinguish it from the strong-greedy algorithm where the direct measure of the error is used. A strong-greedy algorithm requires the evaluation of the high fidelity solution for each instance of the parameters in the parametric domain, which is against the principle of model order reduction. For this reason the error estimator $\Delta(\bmu)$ has to be computationally cheap to evaluate, for every $\bmu\in\mathcal{D}$. A key feature of this strategy is the ability to construct subspaces which preserve the type of decay of the Kolmogorov $n$-width inherited by the underlying problem. This aspect is investigated and illustrated in \cite{Binev2011}.
	
	Let us define the residual at the point $z_j$ as the vector generated by replacing $\hat{\bff u}(z_j;\bmu)$ by $\hat{\bff u}_r(z_j;\bmu)$ in (\ref{eq:LT}), i.e.
	\begin{equation*}
		\bff r(z_j;\bmu)=\Big(z_j\bff \Id-\bff A(\bmu)\Big)\hat{\bff u}_r(z_j;\bmu)-\bff u_0(\bmu)-\hat{\bff b}(z_j;\bmu),
	\end{equation*}
	so that
	\begin{equation} \label{eq:red}
		\hat{\bff u}_r(z_j;\bmu)=\Big(z_j\bff \Id-\bff A(\bmu)\Big)^{-1}\left(\bff r(z_j;\bmu)+\bff u_0(\bmu)+\hat{\bff b}(z_j;\bmu)\right).
	\end{equation}

	In order to exploit the method in a time window instead of using it uniquely at a given instant, we construct the error estimator for a given time window $[t_0,\Lambda t_0]$, with $\Lambda>1$.
	\vspace{-0.3cm}
	\begin{proposition}
		For a prescribed time window $[t_0,\Lambda t_0]\subset\mathbb{R}^{+}$, $\Lambda>1$, and for all $\bmu\in \mathcal{D}$ it holds that
		\begin{equation*}
			\max_{t\in[t_0,\Lambda t_0]}\Big\|\bff u(t,\bmu)-\bff u_r(t,\bmu)\Big\|\le\Delta(\bmu),
		\end{equation*}
		where the error estimator $\Delta(\bmu)$ is given by
		\begin{equation}\label{eq:est}
			\Delta(\bmu):=\left(\frac{c}{N}\sum_{j=1}^{N-1}e^{\real({z_j}) t_j}\Big|z'_j\Big|\Big\|\Big(z_j\bff \Id-\bff A(\bmu)\Big)^{-1}\Big\|\Big\|\hat{\bff r}(z_j;\bmu)\Big\|\right).
		\end{equation}
		with
		\begin{equation*}
		    t_j=\begin{cases}
		    \Lambda t_0, & \mbox{if }\real(z_j)\ge 0\\ t_0, & \mbox{if }\real(z_j)<0
             \end{cases}
		\end{equation*}
	\end{proposition}
	\begin{proof}
		We have
			\begin{small}
		\begin{align*}
			\max_{t\in[t_0,\Lambda t_0]}&\Big\|\bff u(t,\bmu)-\bff u_r(t,\bmu)\Big\|=\left(\frac{c}{N}\Bigg\|\sum_{j=1}^{N-1} \e^{z_jt_j}z'_j\Big(\hat{\bff u}(z_j;\bmu)-\hat{\bff u}_{r}(z_j;\bmu)\Big)\Bigg\|\right)\nonumber\\
			\le&\left(\frac{c}{N}\sum_{j=1}^{N-1}\e^{\real({z_j})t_j}\Big|z'_j\Big|\Big\|\hat{\bff u}(z_j;\bmu)-\hat{\bff u}_{r}(z_j;\bmu)\Big\|\right)\nonumber\\
			\le&\left(\frac{c}{N}\sum_{j=1}^{N-1}\e^{\real             ({z_j})t_j}\Big|z'_j\Big|\Big\|\Big(z_j\bff \Id-\bff A(\bmu)\Big)^{-1}\Big\|\Big\|\hat{\bff r}(z_j;\bmu)\Big\|\right)=\Delta(\bmu)\label{eq:prop3},
		\end{align*}
	\end{small}
		where we have made use of (\ref{eq:LT}), (\ref{eq:redsol}) and (\ref{eq:red}).
	\end{proof}
	As mentioned, a good error estimator needs to be cheaply computed for all $\bmu\in\mathcal{D}$. The maximum over the time window is fast to compute since the dependence from $t$ only appears in the scalar exponential term, therefore the maximum value is attained at one of the two extremes of the interval according to the sign of $\real(z)$. The dependence of $\Delta(\bmu)$ on $\bmu$ is in the norm of the residual and in the norm of the resolvent. The first one {can be computed in an efficient way} for each value of $\bmu$ (see \cite[Section 4.2.5]{Hesthaven2016}). The second one involves the computation of the smallest singular value of $z_j\bff \Id-\bff A(\bmu)$ whose evaluation is computationally expensive if it has to be performed for many values of $\bmu\in\mathcal{D}$. Therefore the strategy employed is to look for a cheap upper bound for $\|(z_j\bff \Id-\bff A(\bmu))^{-1}\|$. Note that the integration contour is constructed in such a way that $\|(z_j\bff \Id-\bff A(\bmu))^{-1}\|$ is bounded, however the use of this bound in (\ref{eq:est}) is not sharp. In Section \ref{sec:lb} we shall review the methods proposed in literature to deal with this problem and we will describe a new approach based on a singular value optimization procedure. 
	\vspace{-0.2cm}
	\subsection{The inverse Laplace transform Greedy-POD algorithm}
	The Greedy -POD algorithm is a widely used sampling strategy to construct reduced spaces for time dependent problems, it is based on the combined use of POD in time with greedy sampling in the parameter space $\mathcal{D}_{\nu}$ where $\mathcal{D}_{\nu}$ is a finite subset of $\mathcal{D}$. {For instance} it can consist of a regular lattice of size $\nu$, a randomly generated point-set belonging to $\mathcal{D}$ {or it can be chosen as the set of the so called magic points \cite{Maday2008}}. The fact that such strategy is able to construct reduced spaces which preserve the exponential or algebraic convergence rates of the Kolmogorov $n$-width is discussed in \cite{Haasdonk2013}.
	
	Our approach, {that we name Laplace POD-Greedy algorithm,} consists of replacing the POD in time with the one in the Laplace transform domain. The steps for the construction of the reduced spaces are reported in Algorithm \ref{al1}. At each iteration $m$, for the selected vector of parameters $\bmu_m$, the algorithm first computes the Laplace transform of the solution on the quadrature nodes (lines 2-4), then it compresses the set of collected snapshots through a truncated singular value decomposition where the truncation is performed according to a tolerance $tol_{POD}$ (line 6). Finally it computes the vector of parameters which maximize the error estimator (line 7). 
	
	The procedure stops when the error estimator $\Delta(\bmu)$ is smaller then the prescribed accuracy $tol$ for each $\bmu\in\mathcal{D}_{\nu}$.
	\begin{algorithm}[t]
		\begin{small}
		\caption{{{Laplace} POD-Greedy construction of the reduced space}}\label{al1}
		\hspace*{\algorithmicindent} \textbf{Input:} $tol$, $tol_{POD}$, $\bmu_1$, $m=1$, $\mathbb{B}= \emptyset$
		
		\hspace*{\algorithmicindent} \textbf{Output:} The reduced space: $\mathbb{V}_{N_r}=\text{span}\{\bff \zeta_1,...,\bff \zeta_{N_r}\}$
		\begin{algorithmic}[1]
			\WHILE{$\Delta(\bmu_{m})>tol$}
			\FOR{$j=1,...,N-1$}
			\STATE Compute $\hat{\bff u}(z_j;\bmu_m)$ \label{line1}
			\ENDFOR
			\STATE Set $\mathbb{B}=\{\mathbb{B},\{\hat{\bff u}(z_1;\bmu_m),...,\hat{\bff u}(z_{N-1};\bmu_m)\}$
			\STATE  Compute $\mathbb{V}_{N_r}=[\bff \zeta_1,...,\bff \zeta_{N_r}]=\text{POD}(\mathbb{B},tol_{POD})$
			\STATE Compute $\bmu_{m+1}=\argmax_{\bmu\in\mathcal{D}_{\nu}}\Delta(\bmu)$
			\STATE Set $m=m+1$
			\ENDWHILE
		\end{algorithmic}
		\end{small}
	\end{algorithm}
	\subsection{An alternative approach: local reduced spaces on the quadrature nodes}
	Once a parameter is selected, Algorithm \ref{al1} computes the high fidelity solution related to that parameter on each quadrature node and adds those solutions to the reduced space. Therefore, the constructed reduced space turns out to be an approximation of the whole solution set (\ref{eq:manLap}). 
	
	A different strategy consists of approximating the solutions manifold associated to each quadrature node, i.e. to consider \vspace{-0.4cm}
	\begin{equation*}
		\hat{\mathcal{M}}_{j}(\mathcal{D}_{\nu})=\{\hat{\bff u}(z_j; \bmu)\;|\; \bmu\in \mathcal{D}_{\nu}\},\quad 1\le j \le N-1;
	\end{equation*}
	and approximating separately each of these solutions sets.
	
	Since we are interested in the solution $\bff u(t,\bmu)$ we still rely on the error estimate (\ref{eq:est}), but this time the stopping criterion for the greedy algorithm is that each term of the sum in (\ref{eq:est}) is required to be smaller than the given tolerance, divided by $N-1$. The strategy is summarized in Algorithm \ref{al2}.
	\begin{algorithm}[t]
	\begin{small}
		\caption{{Greedy construction of the reduced space at the $j$-th quadrature node}}\label{al2}
		\hspace*{\algorithmicindent} \textbf{Input:} $tol_j=tol/(N-1)$, $\bmu_1$, $j$
		
		\hspace*{\algorithmicindent} \textbf{Output:} The reduced space: $\mathbb{V}_{N_r}^{j}$
		\begin{algorithmic}[1]
			
			\STATE Set $\mathbb{B}=[\;],\;m=1,\;\Delta_j(\bmu)=(\ref{eqA2})$
			\label{line2}
			\WHILE{$\Delta_j(\bmu_{m})>tol_j$}
			\STATE Compute $\hat{\bff u}(z_j;\bmu_m)$
			\STATE Set $\mathbb{B}=\{\mathbb{B},\hat{\bff u}(z_j;\bmu_m)\}$
			\STATE $\mathbb{V}^{j}_{N_r}={\text{orth}}\{\mathbb{B}\}$
			\STATE Compute $\bmu_{m+1}=\argmax_{\bmu\in\mathcal{D}_{\nu}}\Delta_j(\bmu)$
			\STATE Set $m=m+1$
			\ENDWHILE
			
		\end{algorithmic}
		\end{small}
	\end{algorithm}
	In general, the reduced spaces constructed at each quadrature point by Algorithm \ref{al2} have different sizes. Indeed, due to the presence of the exponential term in 
	\begin{equation}
	\Delta_j(\bmu)=\frac{c}{N}\max_{t\in[t_0,\Lambda t_0]}\left(\e^{\real({z_j}) t}\right)\Big|z'_j\Big|\Big\|\Big(z_j\bff \Id-\bff A(\bmu)\Big)^{-1}\Big\|\Big\|\hat{\bff r}(z_j;\bmu)\Big\|,\label{eqA2}
	\end{equation}
	the size of the spaces is expected to be larger for those nodes with larger positive real part and smaller for the nodes with smaller real part. 
	
	The main difference between Algorithm \ref{al1} and Algorithm \ref{al2}, is that the first collects all the basis vectors together and then it constructs a unique reduced space, while the second one builds a specific reduced space for each quadrature node. As a result, the size of the reduced space in the first approach turns out to be significantly larger than any of the sets constructed with the second one.
	
	From a computational point of view, with Algorithm \ref{al2}, we require a greedy type search on each quadrature point while Algorithm \ref{al1} has only one greedy type search. For a fixed accuracy $tol$, since the dimension of the local reduced spaces is smaller than the dimension of the space constructed with Algorithm \ref{al1}, the computation of $\Delta(\bmu)\;\forall\bmu\in\mathcal{D}_{\nu}$ is faster with the local approach, although it needs to be performed at each quadrature point. Note that the construction of the reduced spaces on the quadrature nodes in Algorithm \ref{al2} can be parallelized in a natural way. From a storage point of view the cost is higher for Algorithm \ref{al2} since, in Algorithm \ref{al1}, redundant information across the nodes are compressed thanks to the POD.
	
	We will discuss in more detail the computational and storage costs associated to the two algorithms in Section $4$, which is dedicated to numerical experiments.
	\vspace{-0.2cm}
    \subsection{Optimizing the resolvent norm}\label{sec:lb}
	The error estimate (\ref{eq:est}) depends on the resolvent norm. Indicating ${\sigma_{min}}(\bmu)$ the smallest singular value of $z\bff \Id-\bff A(\bmu)$ it holds \vspace{-0.4cm}
	\begin{equation*}
		\|(z\bff \Id-\bff A(\bmu))^{-1}\|=\frac{1}{{\sigma_{min}}(\bmu)}.
	\end{equation*}
    The direct computation of ${\sigma_{min}}(\bmu)$ $\forall\bmu \in\mathcal{D}_{\nu}$ is an expensive procedure that is better to avoid in the framework of reduced order methods. Few techniques have been developed to perform a fast and cheap estimation of ${\sigma_{min}}(\bmu)$; in particular we mention the Successive Constrain Minimization method (SCM) \cite{Huynh2007,Chen2009,Hesthaven2012} and the heuristic strategy based on radial basis functions interpolation (RBFI) \cite{Manzoni2015}.
	
	The SCM is an offline/online procedure in which generalized eigenvalue problems of size $N_h$ need to be solved during the offline phase. The online part is then restricted to provide a lower bound $\tilde{\sigma}_{LB}(\bmu)$ of the smallest singular value ${\sigma_{min}}(\bmu)$ for each new parameter $\bmu\in\mathcal{D}_{\nu}$ with a computational complexity that is independent of the dimension $N_h$. Nevertheless, although this procedure enables a very rapid online evaluation of $\tilde{\sigma}_{LB}(\bmu)$, it still requires a quite expensive offline stage (see the numerical experiments of Section \ref{sec:Num}), which may compromise the efficiency of the whole reduction process. The heuristic strategy proposed in \cite{Manzoni2015} combines a radial basis interpolant to the smallest singular value, incorporating suitable criteria to ensure its positivity, with an adaptive choice of interpolation points through a greedy procedure. In this way it is possible to obtain a reliable approximation of ${\sigma_{min}}(\bmu)$, whose offline construction and online evaluation are fast. This strategy turns out to be effective in our applications but it lacks rigorous theoretical foundations. In fact, for multivariate interpolation the error introduced by the approximation is not quantifiable as in the one-dimensional case. We propose instead an approach aiming to approximate the lower bound
	\begin{equation} \label{eq:sigmaLB}
		{\sigma}_{LB}=\inf_{\bmu\in\mathcal{D}}{\sigma_{min}}(\bmu).  	
	\end{equation}
	The optimization is performed by a gradient-descent iteration method, i.e. given the parameter $\bmu_i$ at iteration $i$ we compute the update at iteration $i+1$ by a suitable step descent gradient method\vspace{-0.5cm}
	\begin{equation*}
		\bmu_{i+1}=\bmu_i-\delta_i\nabla{\sigma_{min}}(\bmu_i)\quad i\ge 0
	\end{equation*}
	where $\bmu_0$ is suitably chosen and $\delta_i$ is a variable step size which may be determined through a line search. Since the optimization problem (\ref{eq:sigmaLB}) is not convex, this method converges in general to a local minimum. Therefore, in order to be effective, one should have some knowledge about the behaviour of ${\sigma_{min}}(\bmu)$ with respect to the parameters and should combine the local optimization with some global strategies. For instance one could set different starting point $\bmu_0$ or partition $\mathcal{D}$ into several sub-domains $\mathcal{D}_j$ such that $\mathcal{D}=\bigcup_{j=1}^{J}\mathcal{D}_j$ and then minimize ${\sigma_{min}}(\bmu)$ on each sub-domain. Another approach could be 
	that of combining the proposed gradient strategy with the eigenvalue optimization developed by Mengi et al. \cite{Mengi2014}, which addresses the global minimization of a prescribed eigenvalue of a Hermitian matrix-valued function depending analytically on its parameters in a box. Moreover this approach may be accelerated in terms of computational speed by employing the subspace projection techniques discussed in \cite{Kangal2016,Kressner2014}.
	
	The expression of the partial derivatives associated to $\nabla{\sigma_{min}}(\bmu)$, with respect to $\bmu$, can be computed analytically. In fact we have for a simple non zero singular value (see e.g. \cite{Kato2013,Horn1985}),
	\begin{equation*}
		\frac{\partial {\sigma_{min}}(\bmu)}{\partial \bmu_i}=\real\left(\hat{\bff u}^*\left(\frac{\partial \bff A({\bmu})}{\partial \bmu_i}\right)\hat{\bff v}\right)\quad\text{for}\quad i=1,...,d;
	\end{equation*}
	where $\hat{\bff u}$ and $\hat{\bff v}$ are the associated left and right singular vectors, respectively, to ${\sigma_{min}}(\bmu)$ (see \cite[Section 3]{Guglielmi2020} for more details). Due to assumption (\ref{eq:affine2}) we have
	\begin{equation*}
		\frac{\partial \bff A({\bmu})}{\partial \bmu_i}=\sum_{q=1}^{Q_{\bff A}}\frac{\partial \Theta^{A}_{q}({\bmu})}{\partial \bmu_i}{\bff A}_{q}\quad\text{for}\quad i=1,...,d.
	\end{equation*}
	An advantage of this optimization method is that we do not need to prescribe a discrete counterpart of $\mathcal{D}$, while for SCM and interpolation, we need to use a discrete set as training sample.
	\subsection{Choice of the integration profile}
	The integration profile is chosen according to the information provided by the approximation of the weighted $\varepsilon$ - pseudospectrum of $\bff A(\bmu)$ (see \cite{Guglielmi2018}). This set is crucial in order to control the numerical error (due to the fact that we operate in finite arithmetic precision). This quantity was estimated in \cite{Guglielmi2020} as \vspace{-0.5cm}
	\begin{equation*}
		{err}_{N}^{num}\le \frac{c}{N}\sum_{j=1}^{N-1}\e^{\real{(z(x_i))}t}\|(z_j\bff \Id-\bff A(\bmu))^{-1}\|\|\bff r_j(\bmu)\||z'(x_j)|;\label{Sec3.5-eq1}
	\end{equation*}
	where $\bff r_j(\bmu)$ is the residual originated by the numerical solution of
	\begin{equation*}
	(z_j\bff \Id-\bff A(\bmu))\hat{\bff u}=\bff u_0(\bmu)+\hat{\bff b}(z_j;\bmu).
	\end{equation*} 
	Therefore, in order to construct a unique integration profile, one needs information about the $\varepsilon$-pseudospectrum of $\bff A(\bmu)$ for each $\bmu\in\mathcal{D}_{\nu}$. The approximation of weighted $\varepsilon$-pseudospectrum level sets is computationally expensive and despite new methods have been recently proposed to reduce this cost (see \cite{Guglielmi2020}), it is not reasonable to compute this set for many instances of $\bmu\in\mathcal{D}$. To overcome this problem we adopted the following strategy. Starting from an initial profile computed for a certain parameter $\bmu^*\in\mathcal{D}_{\nu}$:
	\begin{enumerate}
		\item for each quadrature node $z_i$ we compute 
		\begin{equation*}
			{\sigma}^i_{LB}=\min_{\bmu\in\mathcal{D}_{\nu}}{\sigma^i_{min}}(\bmu),\quad\bmu^i=\argmin_{\bmu\in\mathcal{D}_{\nu}}{\sigma^i_{min}}(\bmu),
		\end{equation*}
		according the procedure described in Section \ref{sec:lb};
		\item for each $\bmu^i$ and for $t\in[t_0,\Lambda t_0]$, $\Lambda>1$, we estimate ${err}_N^{num}$ (see \eqref{eq:est})  as
		\begin{align}
			{err}_N^{num}\le&\frac{c}{N}\sum_{j=1}^{N-1}\e^{\real{(z(x_j))}t}\|(z_j\bff \Id-\bff A(\bmu^i))^{-1}\|\|\bff r_j(\bmu^i)\||z'(x_j)|\nonumber\\
			\le&\frac{c}{N}\sum_{j=1}^{N-1}\e^{\real{(z(x_i))}t}\frac{1}{{\sigma}^j_{LB}}\|\bff r_j(\bmu^i)\||z'(x_j)|;\label{Sec3.5-eq2}
		\end{align}
		\item if $err_N^{num}$ is smaller than the prescribed tolerance for each $i$ we keep our initial guess as integration profile, otherwise we construct the integration profile associated to one of the $\bff A(\bmu^i)$ for which (\ref{Sec3.5-eq2}) is maximal and we iterate steps $1$-$2$-$3$ until a suitable integration profile is determined.
	\end{enumerate}
	This procedure turns out to be effective in all our test problems, which we are going to present in the next section. 
	\vspace{-0.2cm}
	\section{Numerical tests}\label{sec:Num}
	In this section we show few numerical results to validate our method. We first test our approach on the Black–Scholes and Heston equations. The Black–Scholes model here is the same as the one considered in \cite{Hout2009}, while for the Heston model, following \cite{Hout2008}, we consider a slightly different boundary condition with respect to that taken in \cite{Hout2009}. Then we consider a linear advection equation, which is known to be challenging for reduced order methods based on time step discretization. With these test problems we demonstrate the effectivity of our methodology; moreover we show that our approach does not generate a large size reduced spaces when dealing with linear hyperbolic problems. We refer to the reduction associated to the time step discretization as ``classical reduced order" or classical reduction while we call ``Laplace reduced order" the strategy we have presented in this article. In the following numerical examples we measure the relative reduction error over the parametric domain $\mathcal{D}$ and the time window $[t_0,\;\Lambda t_0]$, $\Lambda>1$, as
	\begin{equation}\label{eq:errBS}
		E_r(\mathcal{D},t_0,\Lambda t_0)=\max_{\bmu\in\mathcal{D}}\max_{t\in[t_0,\;\Lambda t_0]}\frac{\|\bff u(t;\bmu)-\bff u_r(t;\bmu)\|}{\|\bff u(t;\bmu)\|}.
	\end{equation}
	Regarding the classical reduction strategy, for the Black-Scholes test problem we have constructed the classical reduced spaces with the greedy-POD algorithm while, for the Heston and advection equations we used POD to compress the full order solutions computed at each time step on a prescribed train parameter sample. {These comparisons are done in terms of reduction error \eqref{eq:errBS} versus the reduced spaces size and computational time versus reduction error.} All the computations were performed using \verb|Matlab 2021a| on a laptop with 2.60 GHz Intel Core i7 processor.
	\vspace{-0.2cm}
	\subsection{Black–Scholes equation}
	\label{subsec:BS}
	The well known (deterministic) Black-Scholes equation (see e.g. \cite{Black1973}) has the following form:
	\begin{equation}\label{eq:BS}
		\frac{\partial u}{\partial \tau}=\frac{1}{2}\sigma^2s^2\frac{\partial^2u}{\partial s^2}+rs\frac{\partial u}{\partial s}-ru,\;\;s>L,\;\;0<\tau\le t,
	\end{equation}
	for $L$, $t$ given, where the unknown function $u(s,\tau)$ stands for the fair price of the option when the corresponding asset price at time $t-\tau$ is $s$, and $t$ is the maturity time of the option. Moreover $r\ge0$, $\sigma>0$ are given constants (representing respectively the interest rate and the volatility). In practise we consider a bounded spatial domain, setting $L<s<S$ for a sufficiently large $S$. We study \eqref{eq:BS} with the following conditions, typical of the European option call (see \cite{Hout2007}):
	\begin{equation}\label{ibcBS}
		u(s,0)=\max(0,s-K),\quad u(L,\tau)=0, \quad u(S,\tau)=S-e^{-r\tau}K,\quad0\le \tau \le t,
	\end{equation}
	being $K$ the reference strike price. Following the same strategy adopted in \cite{Hout2009}, we discretize in space on a uniform grid of $N_h =1000$ points in the interval $[0, 200]$, using the classical centered finite difference scheme.
	We consider the time window $t\in[1, 10]$ and the parametric domain
	\[
	\bmu=(\sigma,\;r)\in\mathcal{D}:=[0.05,\;0.25]\times [0.001,\;0.02],
	\]
	which we restrict to the finite subset $\Xi$ made of $20\times20$ points uniformly distributed in $\mathcal{D}$. Therefore we measure the relative error (\ref{eq:errBS}) for $\mathcal{D}\equiv\Xi$, $t_0=1$ and $\Lambda=10$.
	\begin{figure}[t]
		\centering{
			\subfigure{
				\includegraphics[width=0.4\textwidth]{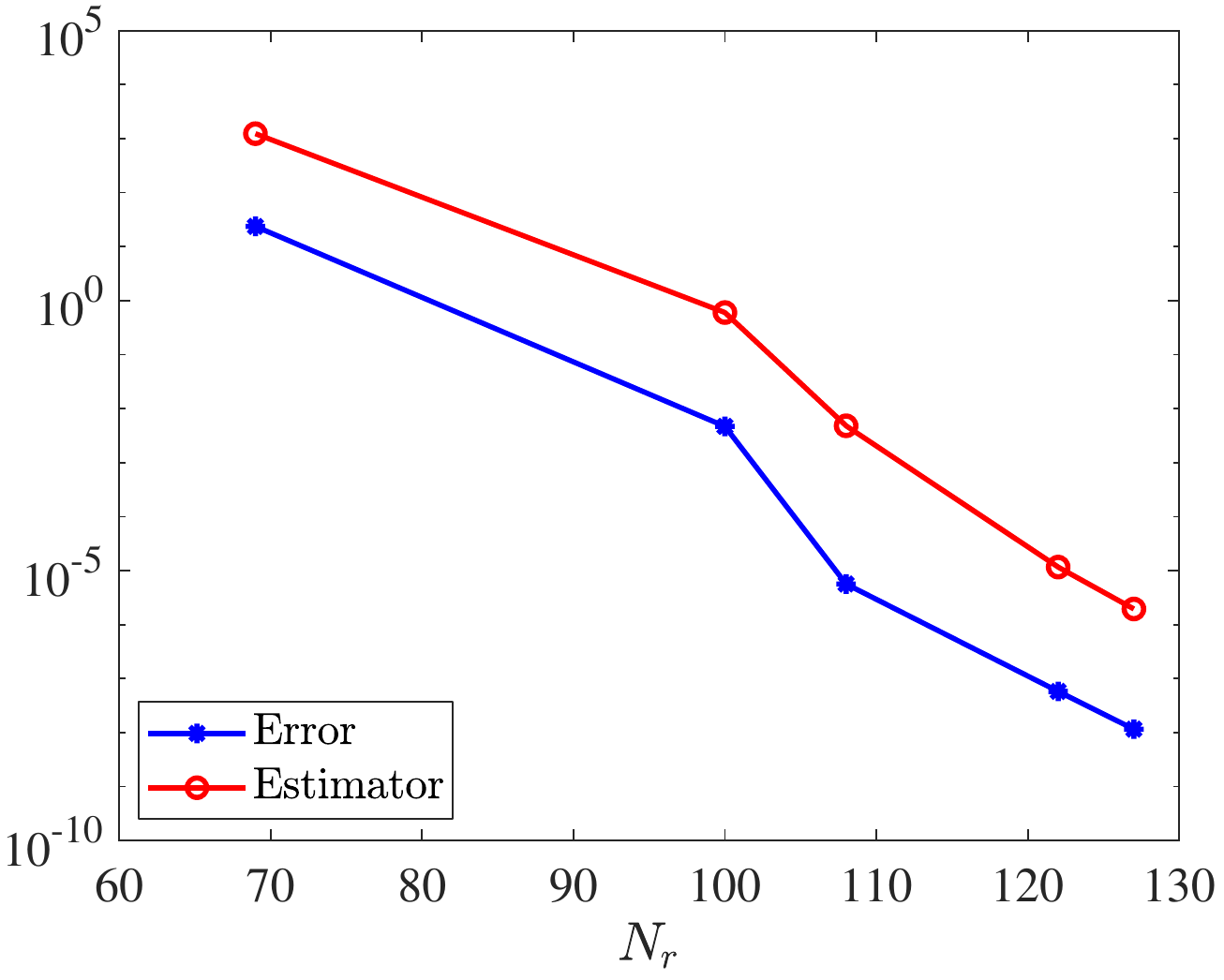}}
		\hspace{0.05cm}
			\subfigure{
				\includegraphics[width=0.4\textwidth]{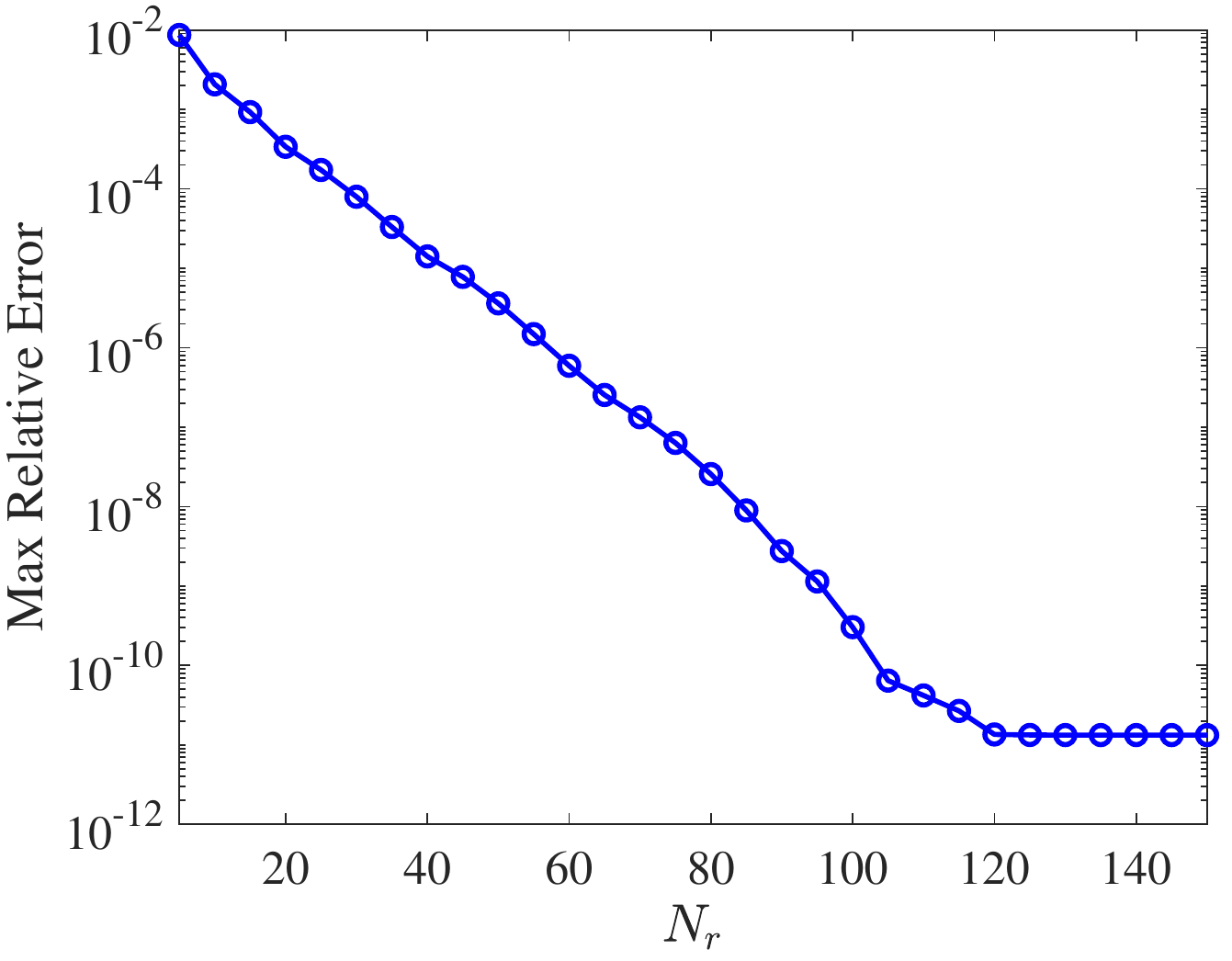}}
		}	
	    \vspace{-0.3cm}
		\caption{Black-Scholes test problem. Behaviour of the Greedy-POD selection algorithm (left); and decay of the relative reduction error (\normalfont{\ref{eq:errBS}}) with respect to the reduced spaces size (right), for a time window $[1, 10]$ and $|\Xi|=400$ uniformly distributed samples.}
		\label{fig2}
		\vspace{-0.2cm}
	\end{figure}

	Figure \ref{fig2} shows the behaviour of Algorithm \ref{al1}. Indeed, on the left picture, we see that the error estimator (\ref{eq:est}) has the same behaviour of the absolute error, while, on the right plot, we recover the desired exponential decay of the relative error with respect to the size of the reduced space.
	\begin{table}[t]  
	\begin{small}
		\centering
		\begin{tabular}{|c|cc|}
			\cline{2-3}
			\multicolumn{1}{c|}{} & CPU time (s) & \# Snap.\\ \hline
			Algorithm \ref{al1} & 19   & 115    \\ 
			Algorithm \ref{al2} & 34   & 1143    \\ \hline
		\end{tabular}
		\caption{CPU time and number of stored snapshots (Snap.) with the greedy-POD strategy of Algorithm \normalfont{\ref{al1}} and the local greedy strategy of Algorithm \normalfont{\ref{al2}}.}
		\vspace{-0.7cm}
		\label{T4}
		\end{small}
	\end{table}
	\begin{figure}
		\centering{
			
			\includegraphics[width=0.4\textwidth]{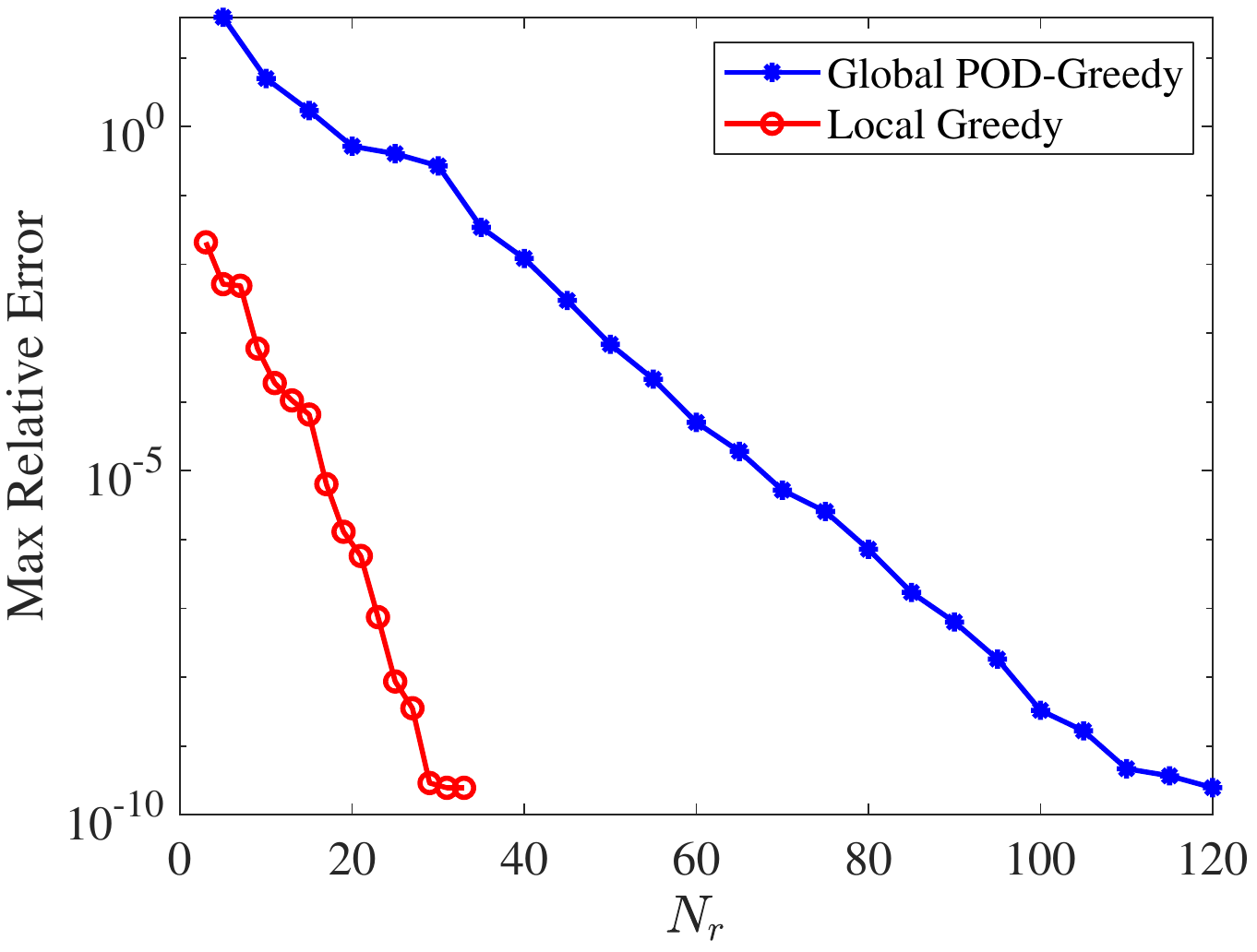}
		}	
	\vspace{-0.4cm}
		\caption{Black-Scholes test problem. Decay of the relative reduction error (\normalfont{\ref{eq:errBS}}), with respect to the size of reduced spaces constructed with the POD-Greedy strategy (Algorithm \normalfont{\ref{al1}}) and the local Greedy one (Algorithm \normalfont{\ref{al2}}). Time window $[1, 10]$ and $|\Xi|=400$ uniformly distributed samples.}
		\vspace{-0.5cm}
		\label{fig11}
	\end{figure}

	Next we illustrate in Table \ref{T4} and Figure \ref{fig11} the performances of Algorithm \ref{al1} and Algorithm \ref{al2}. It turns out that Algorithm \ref{al2} is more expensive in terms of computational time and number of stored snapshots; at the same time the decay of error with respect to the reduced space size is clearly better. Based on these results we consider Algorithm \ref{al2} more suitable to this problem since the gain in the online phase is remarkable and the additional cost in the offline phase is moderate and may be reduced using multi core processors.
	
\begin{figure}[t]
		\centering{
			\subfigure{
				\includegraphics[width=0.4\textwidth]{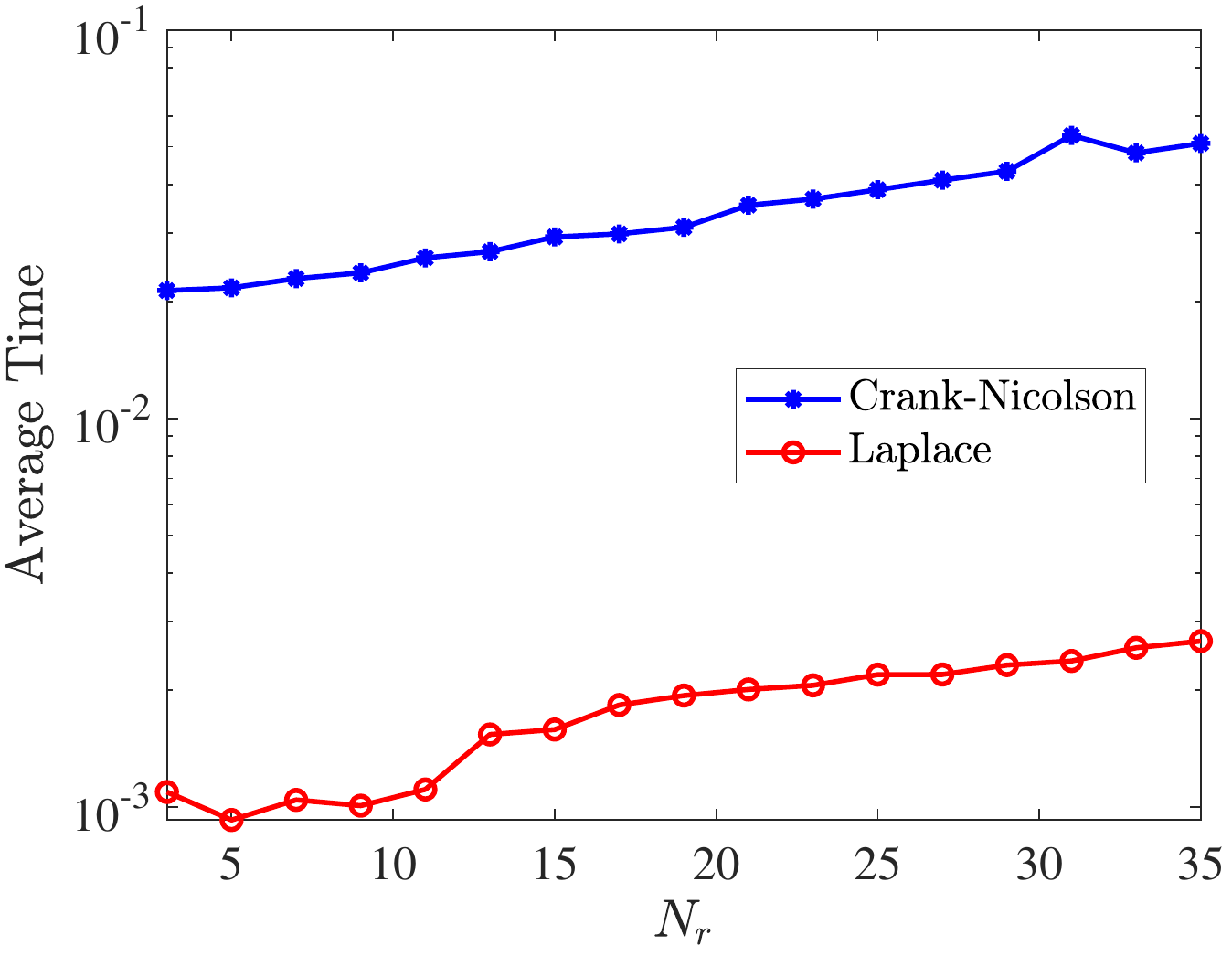}}
			\hspace{0.05cm}
			\subfigure{
				\includegraphics[width=0.4\textwidth]{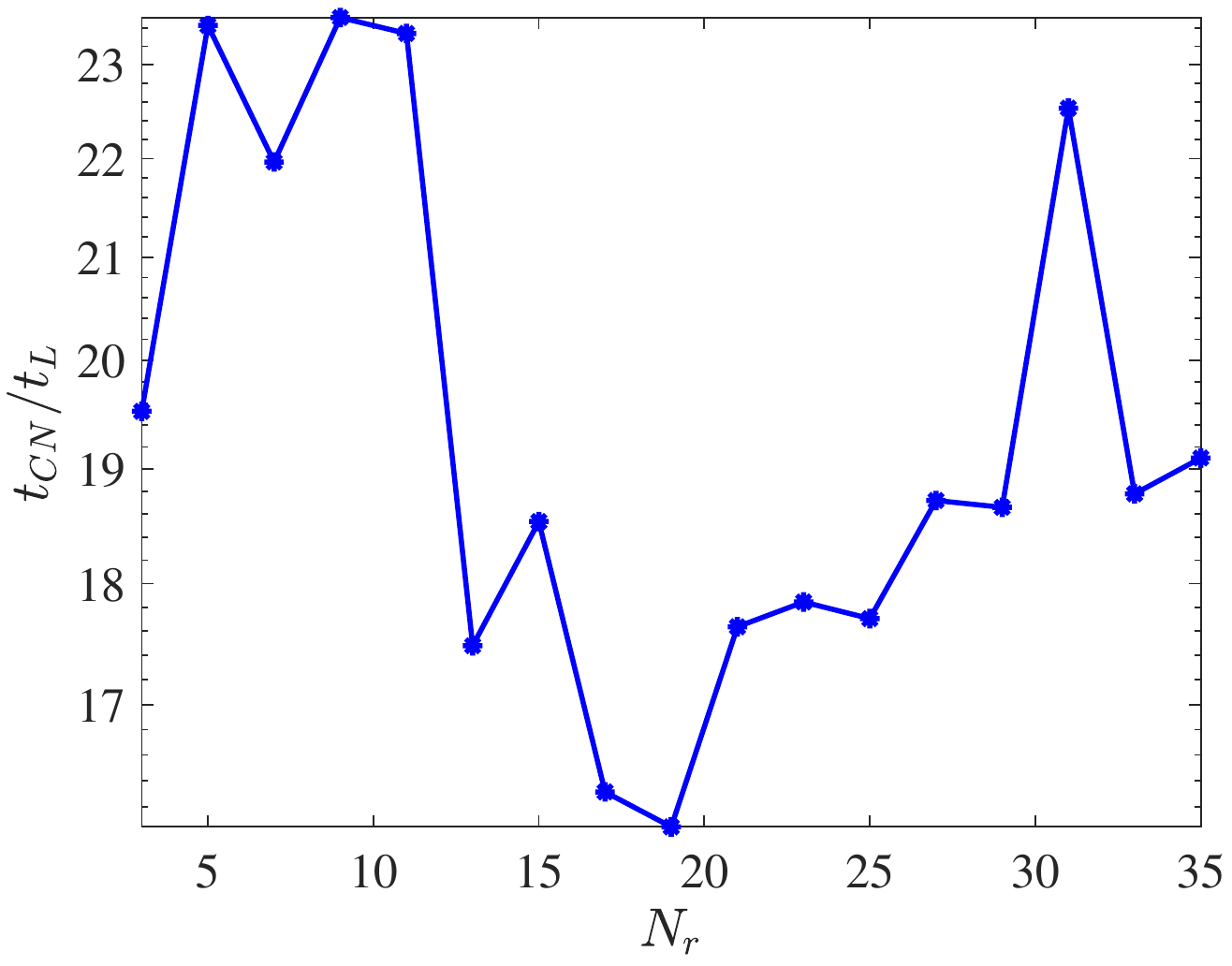}}
		}	
	    \vspace{-0.5cm}
		\caption{Black-Scholes test problem on the time window $[0.1, 1]$. Behaviour of the computational time with respect to size of the reduced spaces for Crank-Nicolson and Laplace. Direct comparison (left) and ratio between the average times (right).}
		\vspace{-0.7cm}
		\label{fig4}
	\end{figure}

\begin{figure}
		
		\centering{
			\subfigure[$\Delta t=10^{-2}$, $tol=10^{-2}$, $t\in {[0.1,1]}$]{
				\includegraphics[width=0.48\textwidth]{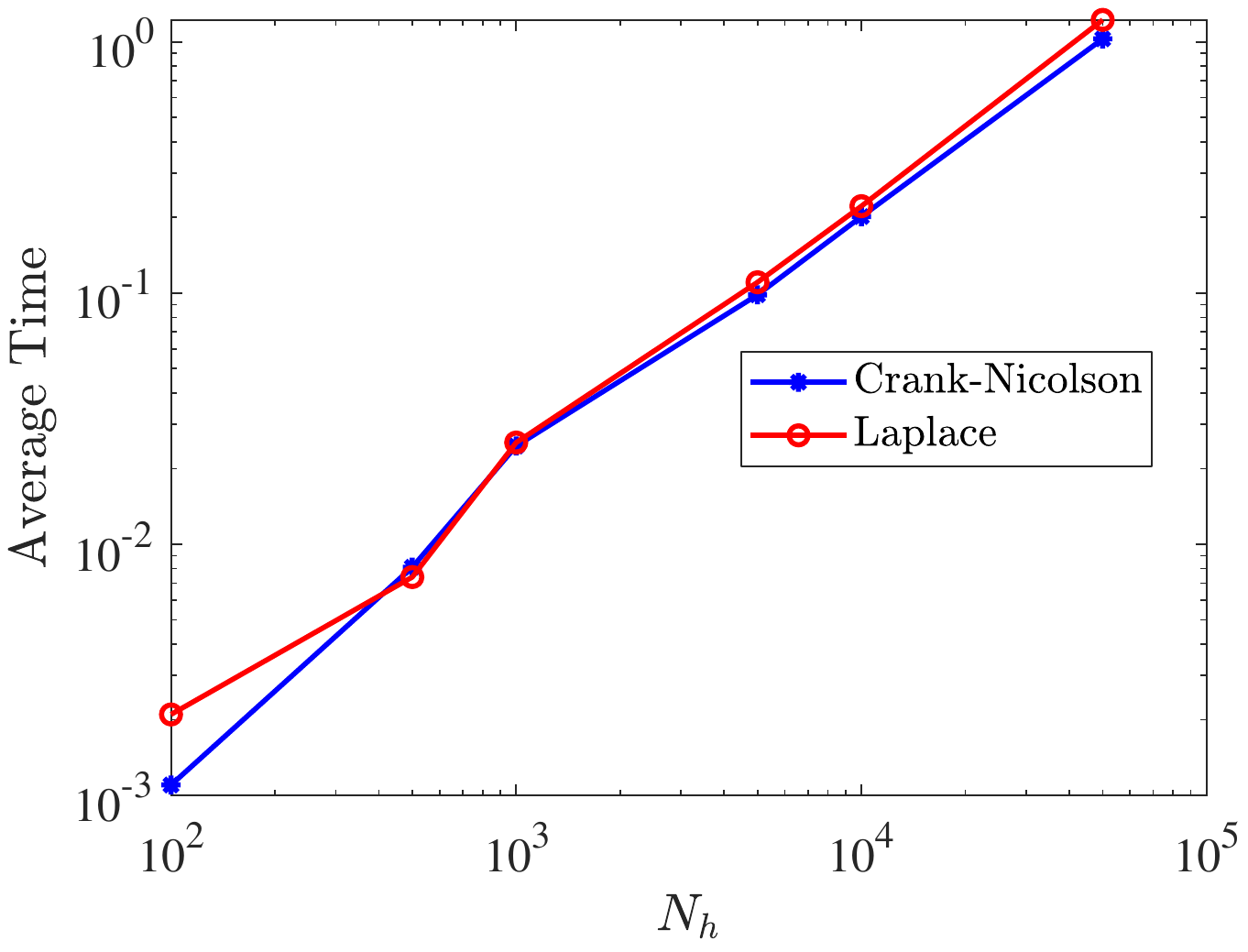}\label{fig8.a}}
			\subfigure[$\Delta t=10^{-2}$, $tol=10^{-2}$, $t\in {[0.1,1]}$]{
				\includegraphics[width=0.48\textwidth]{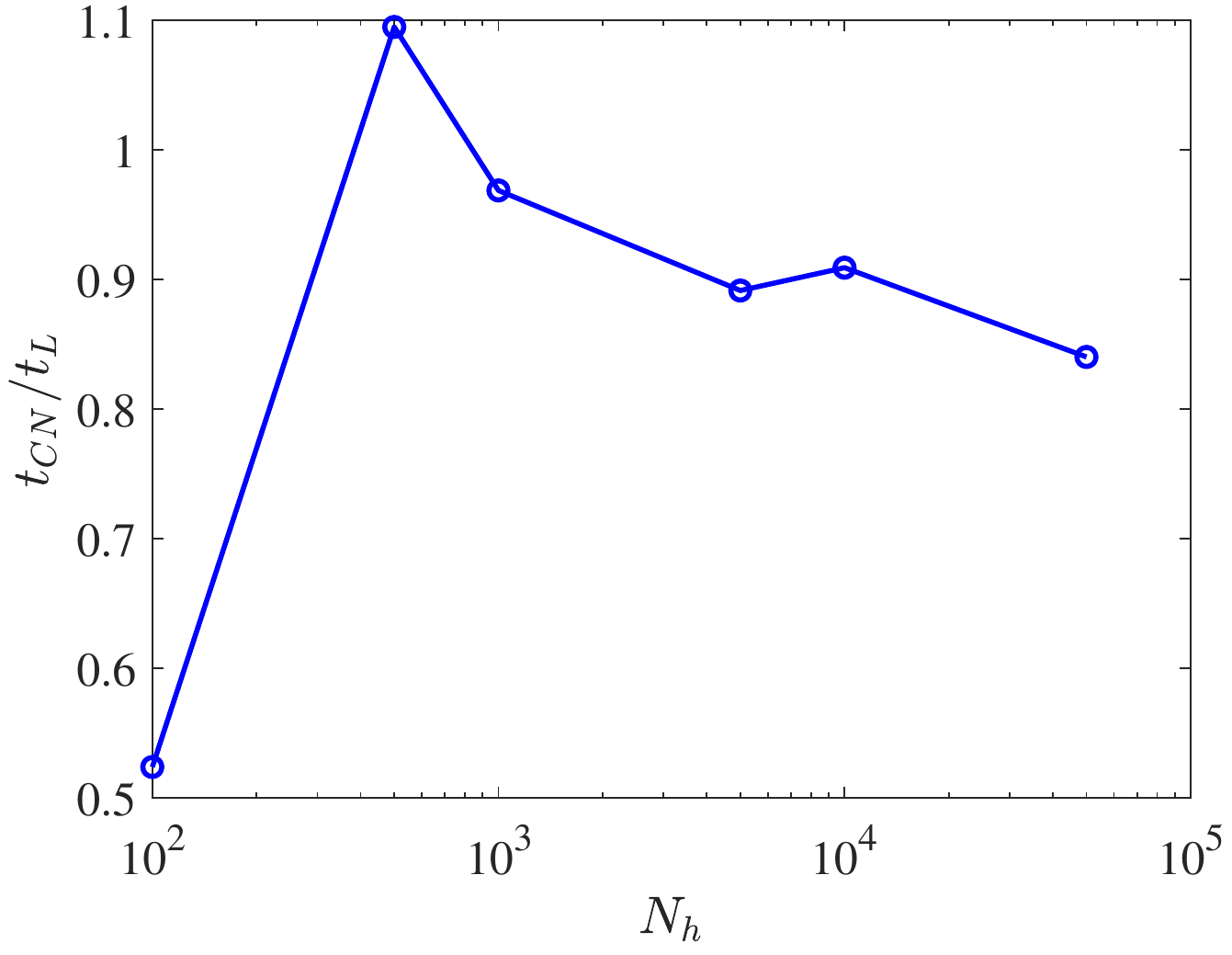}\label{fig8.b}}
			
		}
		
		\centering{
			\subfigure[$\Delta t=10^{-2}$, $tol=10^{-2}$, $t\in {[1,10]}$]{
				\includegraphics[width=0.48\textwidth]{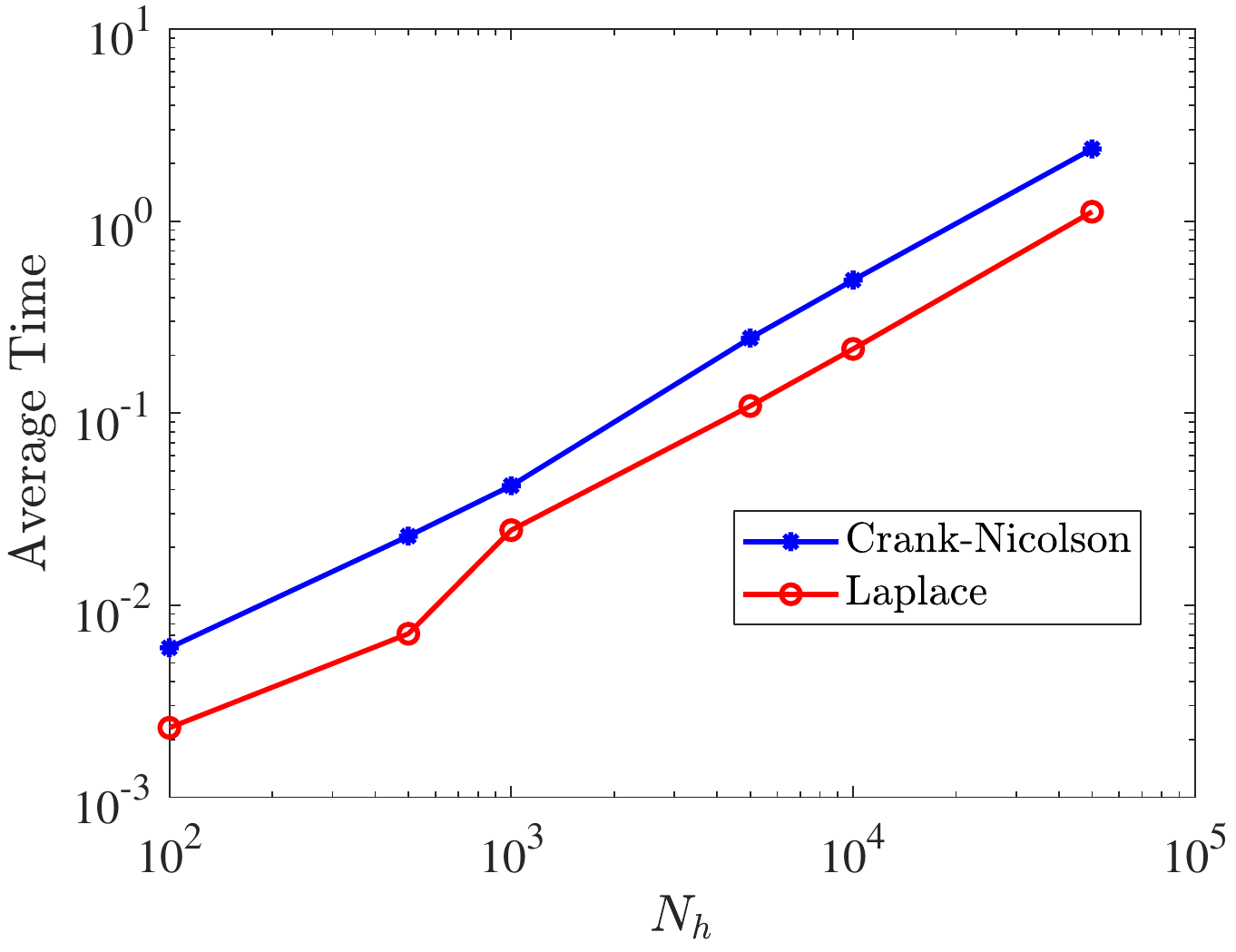}\label{fig8.c}}
			\subfigure[$\Delta t=10^{-2}$, $tol=10^{-2}$, $t\in {[1,10]}$]{
				\includegraphics[width=0.48\textwidth]{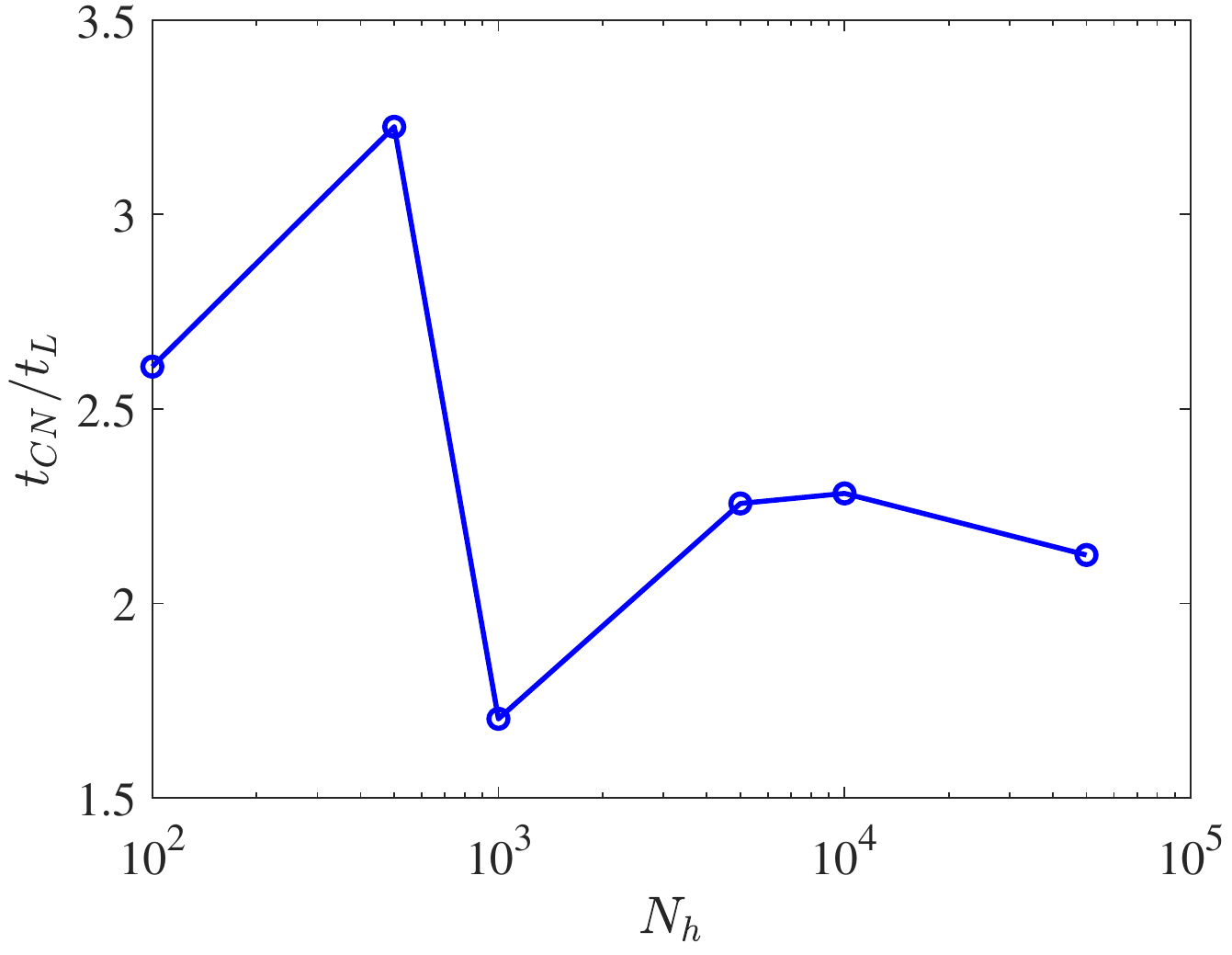}\label{fig8.d}}
		}

		\centering{
			\subfigure[$\Delta t=10^{-4}$, $tol=5\cdot 10^{-6}$, $t\in {[0.1,1]}$]{
				\includegraphics[width=0.48\textwidth]{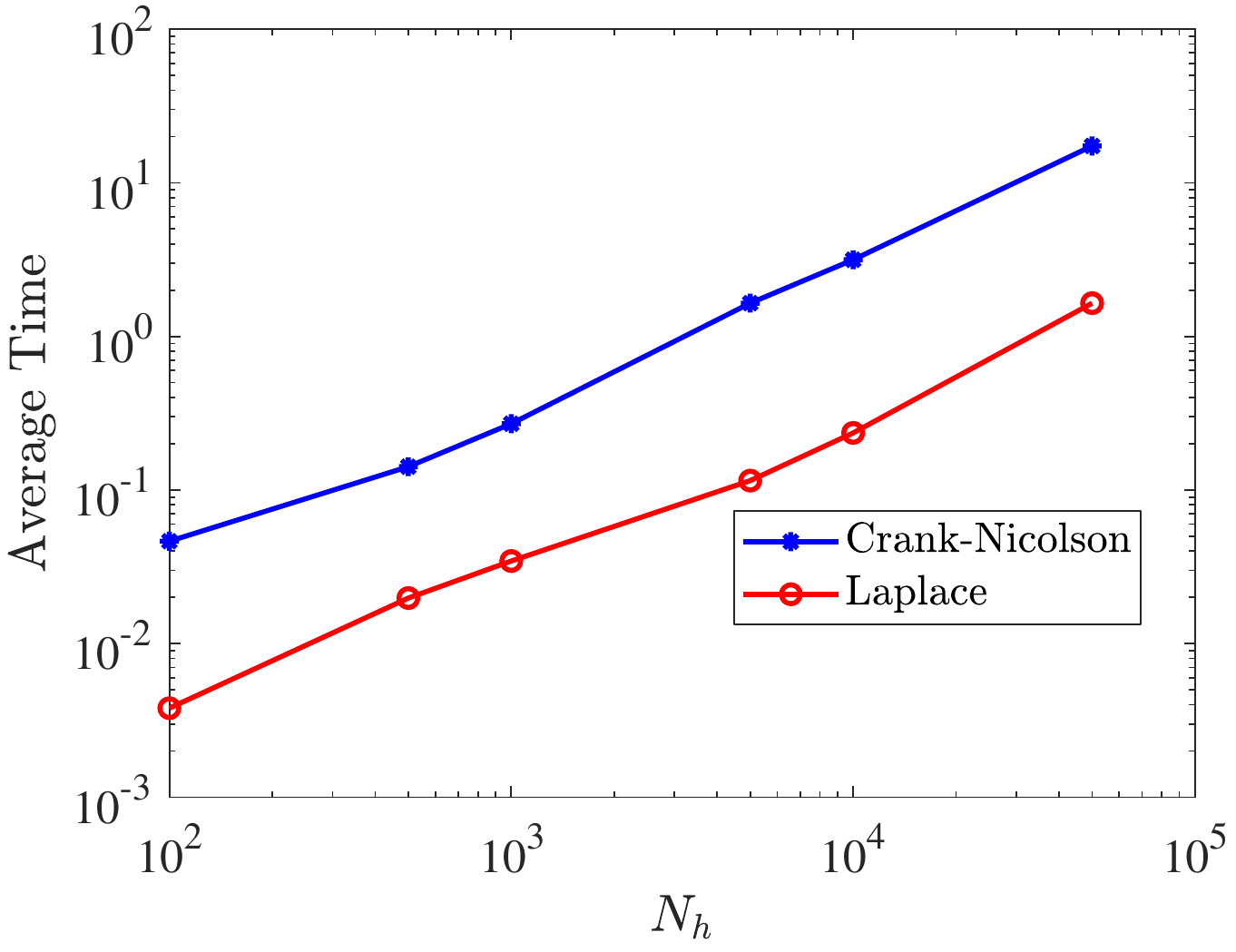}\label{fig8.e}}
			\subfigure[$\Delta t=10^{-4}$, $tol=5\cdot 10^{-6}$, $t\in {[0.1,1]}$]{
				\includegraphics[width=0.48\textwidth]{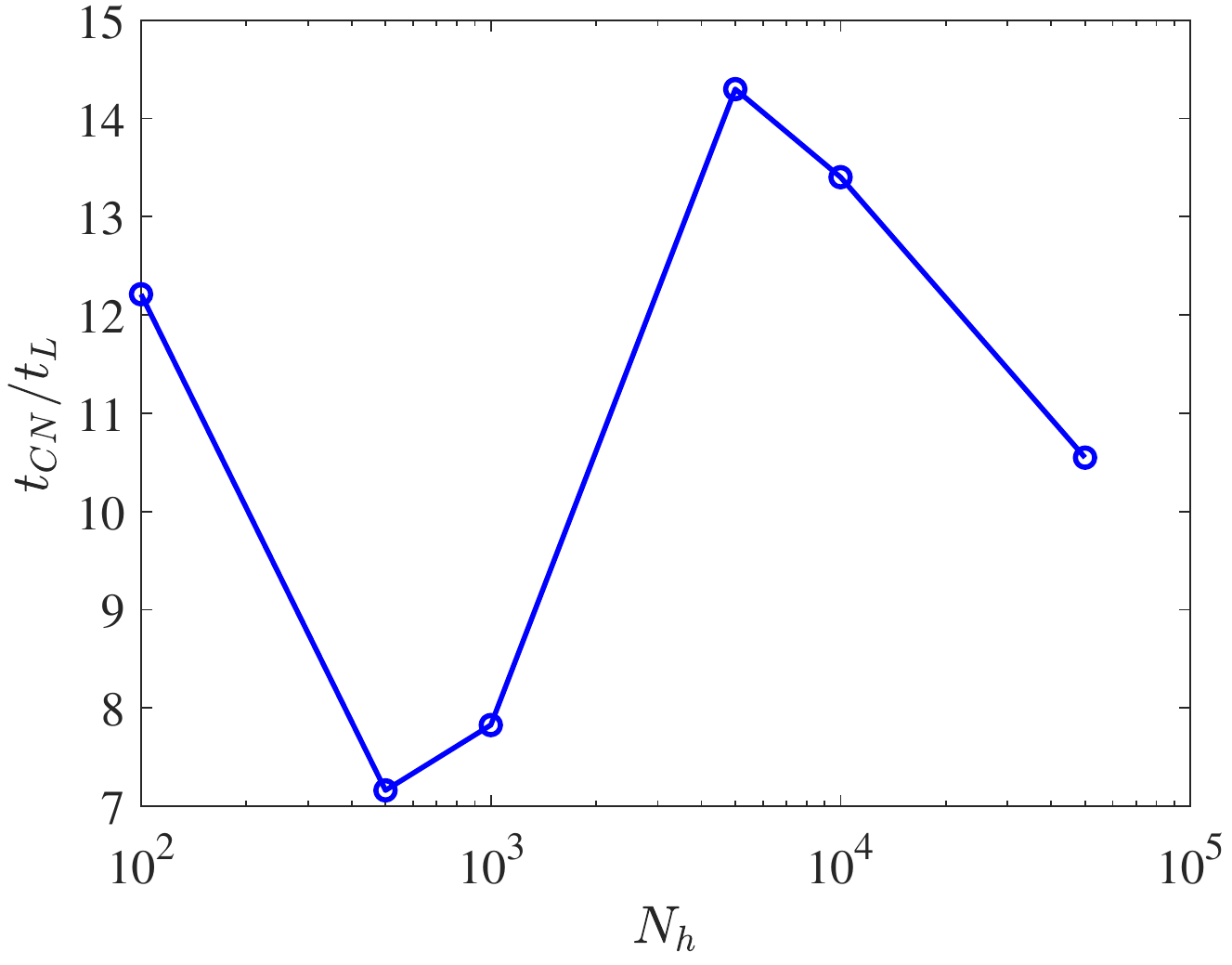}
				
				\label{fig8.f}}
			
		}	
		
		\caption{Black-Scholes equation for different choices of $\Delta t$ (for Crank-Nicolson integrator), $tol$ and time window. On the left: averaged CPU time with respect to the discretization size $N_h$ for Crank-Nicolson and Laplace. On the right: ratio between the CPU time of Crank-Nicolson and Laplace methods.}
		\label{fig8}
	\end{figure}

	We have compared the Laplace reduced {order} method to a classical one based on the Crank-Nicolson method implemented with constant step size $\Delta t=10^{-4}$ on the time window $[0.1, 1]$. The error in the full problem is approximately $ 10^{-3}$ for both Laplace and Crank-Nicolson method. Figure \ref{fig4} displays the CPU time for the solution of the reduced problem, averaged over the parametric domain $\Xi$. We note that the contour integral method is between $16$ and $23$ times faster.
	We also analysed the CPU time in the offline phase, i.e. the time needed to solve the full problem. In particular, we considered the Crank-Nicolson scheme implemented with stepsize: $\Delta t=10^{-4}$ and $\Delta t=10^{-2}$; and the contour integral method with tolerance (for the numerical quadrature): $5\cdot 10^{-6}$ and $10^{-2}$. Results are illustrated in Figure \ref{fig8} were the CPU time is measured as an average of $100$ different computations and is plotted with respect to the size of the discretization $N_h$. When $\Delta t=10^{-2}$, $tol=10^{-2}$ and $T=1$ (Figure \ref{fig8.a} and \ref{fig8.b}) the CPU time is similar for the two solvers. When we set the final time to $T=10$ (Figure \ref{fig8.c} and Figure \ref{fig8.d}) we do not see significant changes in the CPU time of the Laplace method, while the one of Crank-Nicolson clearly increases, making it significantly slower than the contour integral method. Moreover, considering the time window $[0.1,1]$ and increasing the time accuracy, the contour integral method outperforms the time stepping scheme, being $7$ to $14$ times faster (see Figure \ref{fig8.f}). 

\begin{figure}[t]
		\centering{
			\subfigure{
				\includegraphics[width=0.45\textwidth]{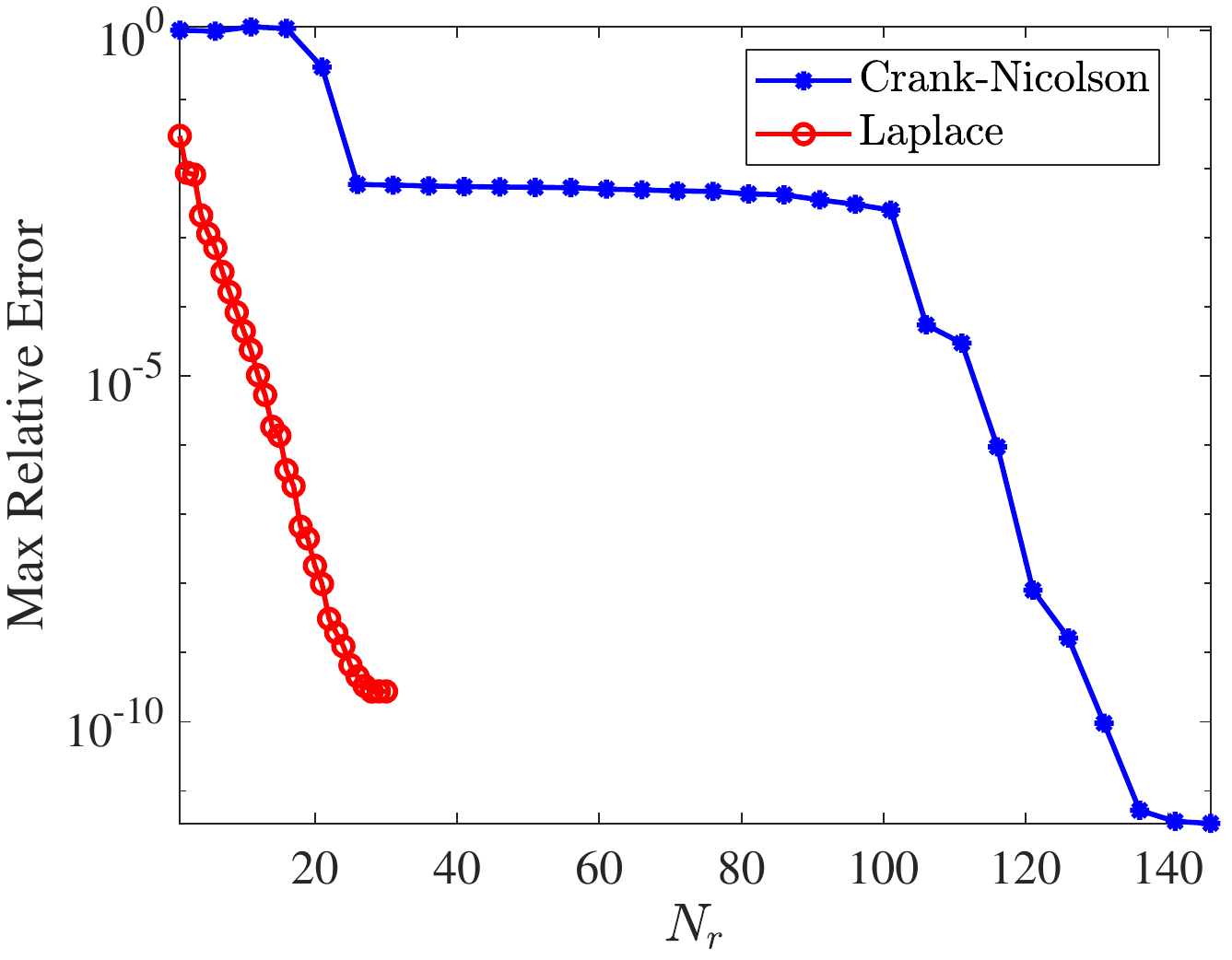}}
			\hspace{0.05cm}
			\subfigure{
				\includegraphics[width=0.45\textwidth]{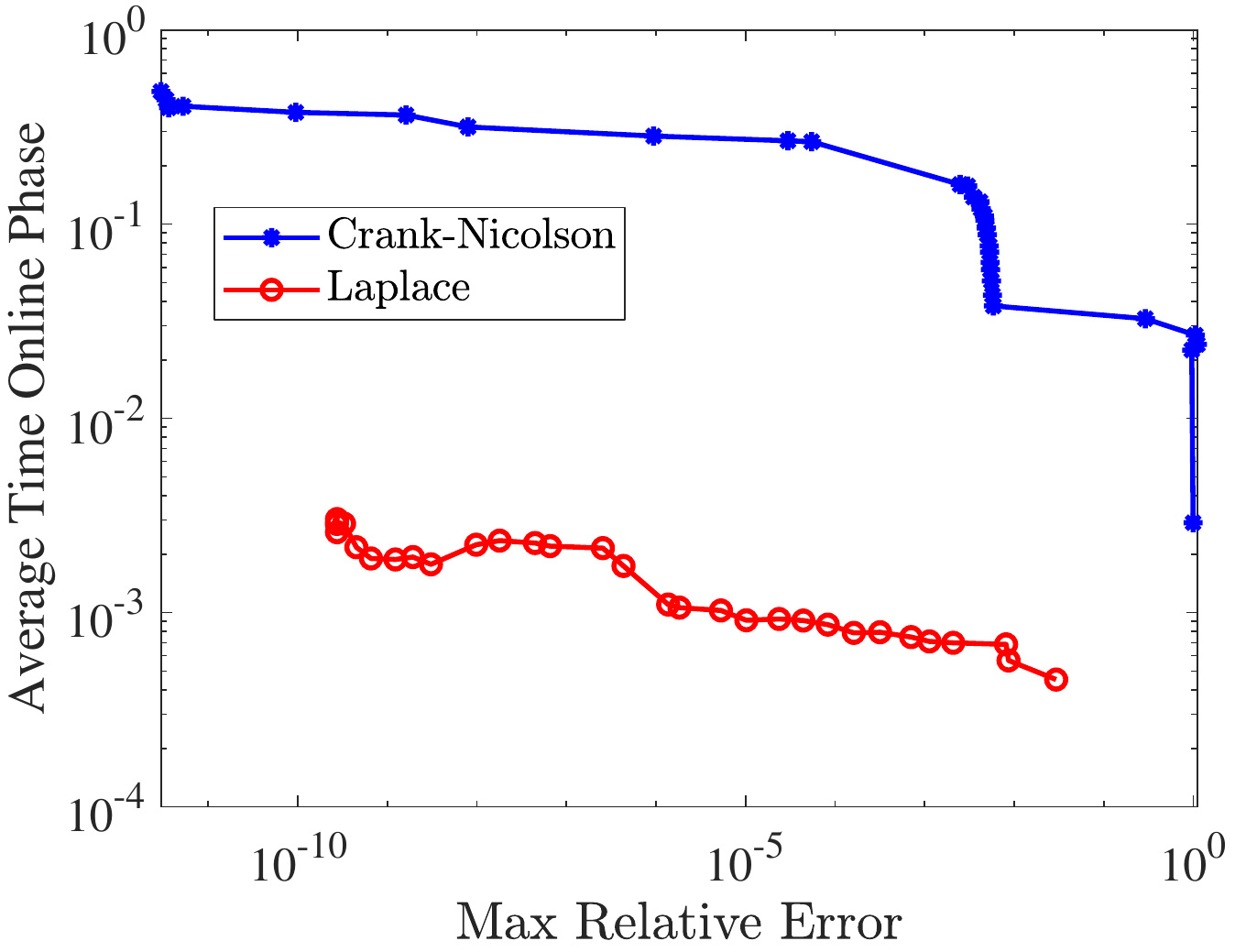}}
		}	
		\caption{Black-Scholes test problem. Decay of the relative reduction error (\normalfont{\ref{eq:errBS}}) with respect to the reduced spaces size (left) and computational time versus the reduction error (right) for the Crank-Nicolson and the contour integral methods. Time window $[0.1, 1]$ and $|\Xi|=400$ uniformly distributed samples.}
		\label{fig3}
	\end{figure}

	Figure \ref{fig3} displays two pictures which compare the two reduction techniques we have proposed. The left one illustrates the behaviour of the maximum relative error with respect to the size of the reduced {space} while the right one is related to the behaviour of the computational time with respect to the associated error. The error originated by the Laplace reduced {order} is always smaller than the one associated to the classical reduced {order} method. Concerning the classical method we observe an initial stagnation which precedes a rapid decrease of the error. This feature can be explained by the fact that the initial datum is only continuous, a feature which, combined with convection, can originate a slow decay of the reduction error for classical reduced order methods (see \cite{Greif2019}). In the right plot it can be seen that the Laplace reduced {order} method is almost two orders of magnitude faster than the classical reduced {order} method.

\begin{figure}
		
		\centering{
			\subfigure[$z_1=0.7747 + 0.0850i$]{
				\includegraphics[width=0.48\textwidth]{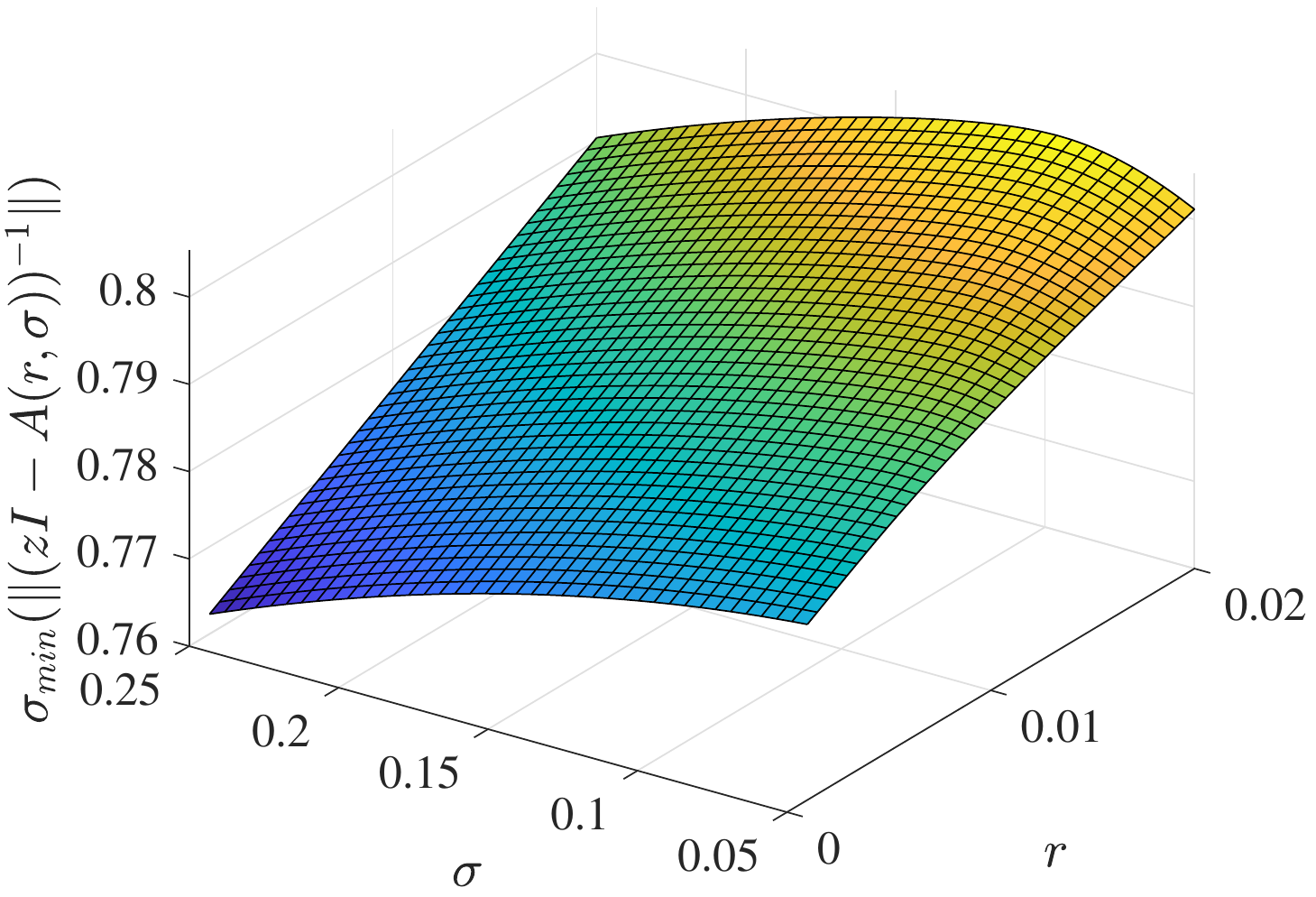}}
			\subfigure[$z_2=-17.5797 + 6.4018i$]{
				\includegraphics[width=0.48\textwidth]{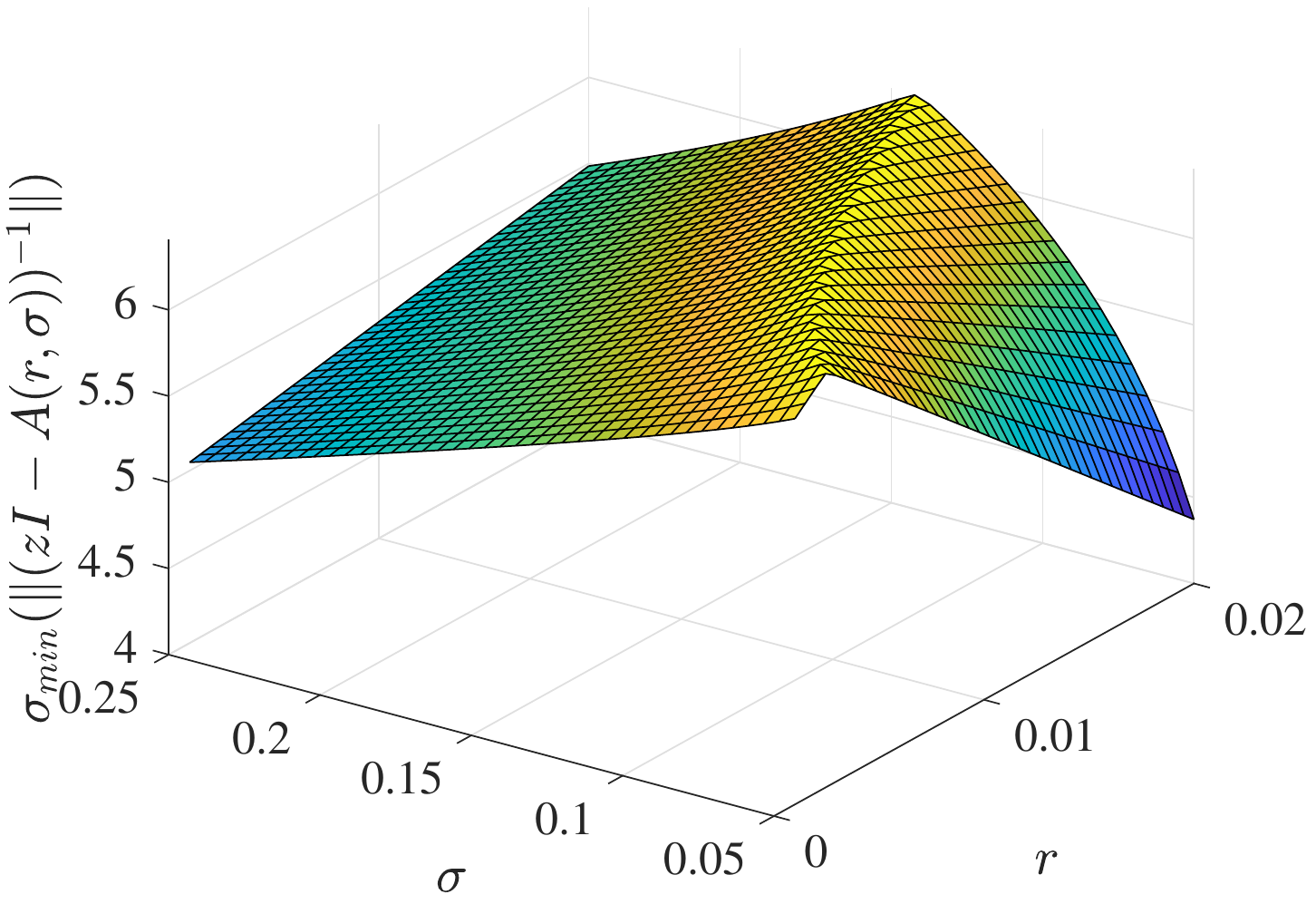}}
			
		}
		
		\caption{Black-Scholes test problem. Plot of the computed smallest singular value for $\bmu\in\Xi$, in two different quadrature points $z_1$ and $z_2$.}
		\label{fig10}
	\end{figure}

	Finally we illustrate the numerical results obtained through the gradient based strategy for the computation of ${\sigma}_{LB}$ (see (\ref{eq:sigmaLB})). As first result we show that the Black-Scholes problem is well suited to be treated with a gradient type method. Figure \ref{fig10} shows the computed smallest singular value for each parameter in $\Xi$ in correspondence of $2$ arbitrarily chosen quadrature points. We see that the behaviour depicted in the second plot is common to most of the quadrature points. It can be seen that the smallest singular value is quite smooth with respect to the parameters in all the considered cases. In Table \ref{T1} we report the computed approximation (which we denote $\sigma_{GS}$) of ${\sigma}_{LB}$ for $4$ quadrature points, together with the relative error and the number of eigenvalue problems we have solved to obtain the final approximation. The value of ${\sigma}_{LB}$ has been estimated as
		\begin{equation}\label{eq:sigLBdisc}
		{\sigma}_{LB}=\inf_{\bmu\in\Xi}{\sigma_{min}}(\bmu).	
	\end{equation}
	We have considered four different starting points, according to the criteria described in Section \ref{sec:lb}. With a maximum of $22$ eigenvalue computations the method is able to return the exact lower bound on the parametric domain. The ability to compute the exact value is due to the fact that the lower bound is attended on a vertex of the parametric domain. In Table \ref{T3} we compare the gradient strategy with both the SCM method \cite{Huynh2007,Chen2009,Hesthaven2012} and the strategy based on RBFI \cite{Manzoni2015} for the Black-Scholes operator\vspace{-0.5cm}
	\begin{equation*}
		\mathcal{A}u=\frac{1}{2}\sigma^2s^2\frac{\partial^2u}{\partial s^2}+rs\frac{\partial u}{\partial s}-ru.
	\end{equation*}
	Note that the negative values of $(\sigma -\sigma_{LB})/\sigma_{LB}$ are justified by the fact that the computed $\sigma$ results to be smaller than the $\sigma_{LB}$ defined in (\ref{eq:sigLBdisc}).
		\begin{table}
		\begin{small}
		\centering
		\begin{tabular}{|c|cccc|}
			\hline
			\multirow{2}{*}{$z$}&\multirow{2}{*}{${\sigma}_{GS}$}&\multirow{2}{*}{${\sigma}_{LB}$}&\multirow{2}{*}{$\frac{{\sigma}_{GS}-{\sigma}_{LB}}{{\sigma}_{LB}}$}&\multirow{2}{*}{\#$EP$}\\
			&&&&\\
			\hline
			\hline
			$0.4190 + 0.0803i$   & $0.4093$    &$0.4093$& $0$&$21$ \\
			
			
			
			
			
			$-3.6612 + 2.3961i$   & $1.4558$    &$1.4558$& $0$&$22$ \\
			
			
			
			$-9.4930 + 3.5718i$   & $2.0782$    &$2.0782$& $0$&$22$\\
			
			
			
			$-17.3555 + 4.4742i$  & $2.4755$    &$2.4755$& $0$& $22$\\
			
			\hline
		\end{tabular}
		\caption{Black-Scholes test problem. Approximations of ${\sigma}_{LB}$ through the gradient type method. The table reports the relative error in the approximation and the number of solved eigenvalue problems (EP) to compute the approximation ${\sigma}_{GS}$. The results are shown for $4$ quadrature points $z$.}
	
		\label{T1}
		\end{small}
	\end{table}
	\begin{table}
	\begin{small}
	\vspace{-0.5cm}
		\centering
		\begin{tabular}{|c|cccc|}
			\hline
			\multirow{2}{*}{Methods}&\multirow{2}{*}{${\sigma}$}&\multirow{2}{*}{$\frac{{\sigma}-{\sigma}_{LB}}{{\sigma}_{LB}}$}&\multirow{2}{*}{\#$EP$}&\multirow{2}{*}{CPU Time}\\
			&&&&\\
			\hline
			\hline
			GS   & $3.3\cdot10^{-4}$    &$-2.0\cdot10^{-3}$& 111&$1.29$ \\
			
			SCM    & $3.1\cdot10^{-7}$ &$-9.9\cdot10^{-1}$& $190$&$>400$\\
			
			RBFI  & $3.3\cdot10^{-4}$ &$\quad1.9\cdot10^{-11}$& $332$&$55.8$\\

			\hline
		\end{tabular}
		\caption{Black-Scholes test problem. Approximation of ${\sigma}_{LB}$ with the gradient solver (\normalfont{GS}), the Successive Constrain Minimization method (\normalfont{SCM}) and the radial basis function interpolation (\normalfont{RBFI}). The methods are compered in terms of: error in the approximation of ${\sigma}_{LB}$, number of solved eigenvalue problem (\#$EP$) and CPU time.}
		\vspace{-1cm}
		\label{T3}
		\end{small}
	\end{table}
	\subsection{The Heston equation}
	The Heston equation \cite{Heston1993} is a $2$D convection diffusion equation given by
	\begin{equation}
		\frac{\partial u}{\partial \tau}=\frac{1}{2}s^2v\frac{\partial^2u}{\partial s^2}+\rho\sigma sv\frac{\partial^2 u}{\partial s \partial v}+\frac{1}{2}\sigma^2v\frac{\partial^2u}{\partial v^2}+(r_d-r_f)s\frac{\partial u}{\partial s}+\kappa(\eta-v)\frac{\partial u}{\partial v}-r_du.
		\label{5.2}
	\end{equation}
	The unknown function $u(s, v,\tau)$ represents the price of a European option when at time $t-\tau$ the corresponding asset price is equal to $s$ and has variance $v$. We consider the equation on the unbounded domain $0\le\tau\le t,\;\;s>0,\;v>0$, where the time $t$ is fixed. The parameters $\kappa>0$, $\sigma>0$, and $\rho\in[-1,1]$ are given. Moreover equation (\ref{5.2}) is usually considered under the condition $2\kappa\eta>\sigma^2$ that is known as the Feller condition (see \cite{Janek2011}). We consider equation (\ref{5.2}) together with the initial condition $u(s,v,0)=\max(0,s-K)$ (where $K>0$ is fixed a priori and represents the strike price of the option), and the boundary condition $u(L,v,\tau)=0,\;\;0\le\tau\le t$.

	For the numerical approximation of (\ref{5.2}), we choose a suitable bounded domain of integration and follow \cite{Hout2008}. In particular, we fix two positive sufficiently large constants $S$, $V$ and we let the two variables $s$, $v$ vary in the set $0\le s\le S,\;\; 0\le v\le V$. On the new boundary, we need to add two extra conditions (specific for the European call option),\vspace{-0.4cm}
	\begin{equation*}
		\frac{\partial u}{\partial s}(S,v,\tau)=e^{-r_f\tau}\;\;\text{and}\;\;u(s,V,\tau)=se^{-r_f\tau},\; \;0\le\tau\le t,
	\end{equation*}
	which are treated analogously to the boundary condition in (\ref{ibcBS}). The spatial discretization we adopted is the same introduced in \cite{Hout2008} with $N_h=10^4$. We set $r_f = 0$, $K = 100$, $L = 0$, $S = 8K$, $V = 5$, and consider the time window $t\in[0.5, 1]$ and two distinct parametric domains for $\bmu =(\sigma,\;r_d,\;\kappa,\;\eta,\;\rho)^T$ and $\bmu=(\kappa,\;\eta)^{T}$
	\begin{align*}
	\mathcal{D}_{1}&\equiv[0.18,\;0.4]\times [0.001,\;0.2]\times [1.2,\;3]\times [0.08,\;0.15]\times [0.21,\;0.9]\subset\mathbb{R}^{5};\\
	\mathcal{D}_{2}&\equiv [1.2,\;3]\times [0.08,\;0.15]\subset\mathbb{R}^{2}.
	\end{align*}
	For each choice of $\bmu$, the remaining parameter values are assumed to be fixed and taken from the reference parameter vector $\bmu^{*}=(0.3,0.02,2,0.1,0.21)^T$.  In our first test we consider $\bmu\in\mathcal{D}_2$ and $|\Xi|=15^2=225$ equidistantly distributed points. To quantify the efficiency of our reduced order method, we have measured the error decay when increasing the dimension $N_r$. For each reduced model, we have computed the relative reduction error according to (\ref{eq:errBS}) with $\mathcal{D}_2=\Xi$, $t_0=0.5$ and $\Lambda=2$.
	\begin{figure}[t]
		\centering{
			\subfigure{
				\includegraphics[width=0.45\textwidth]{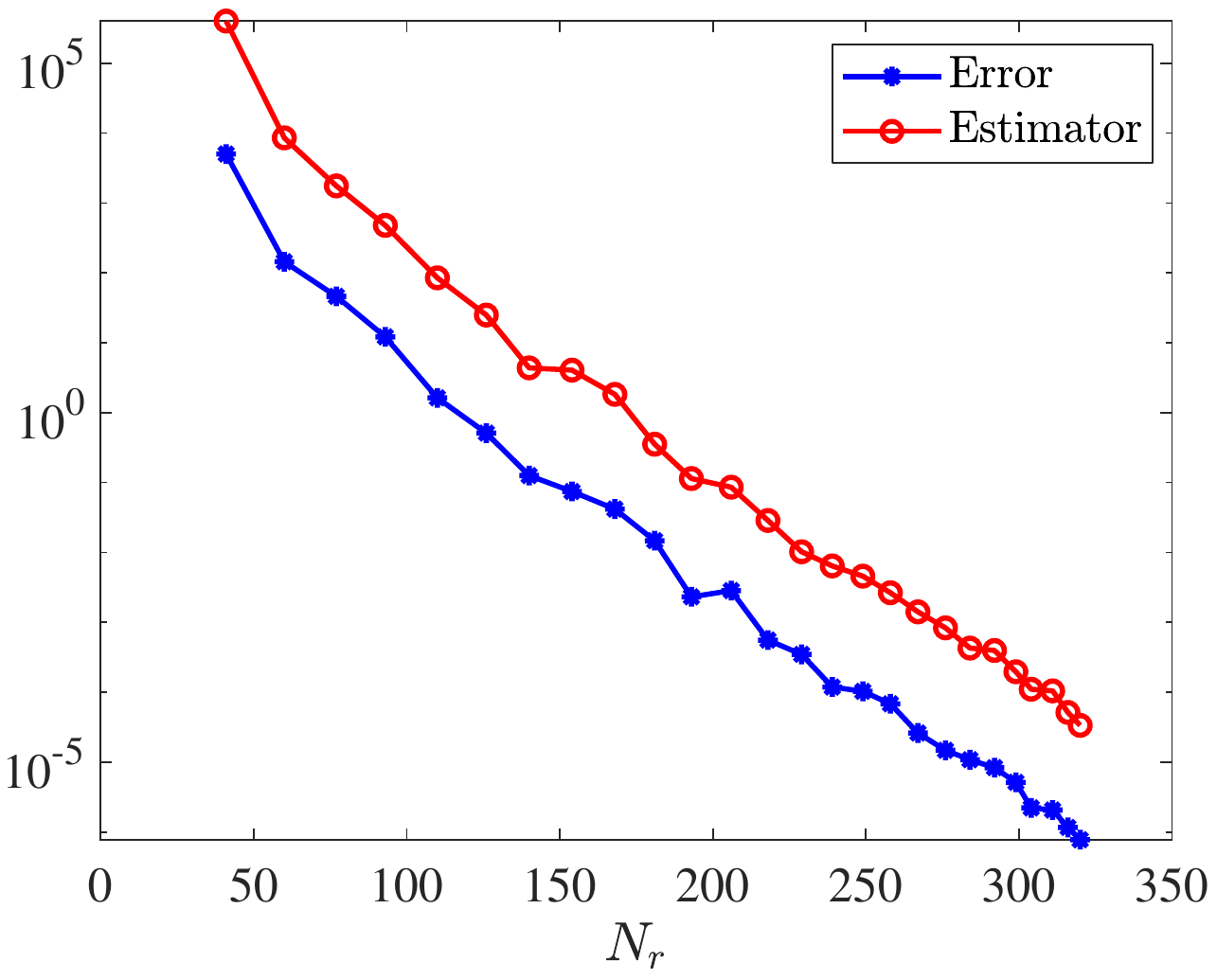}}
			\hspace{0.05cm}
			\subfigure{
				\includegraphics[width=0.45\textwidth]{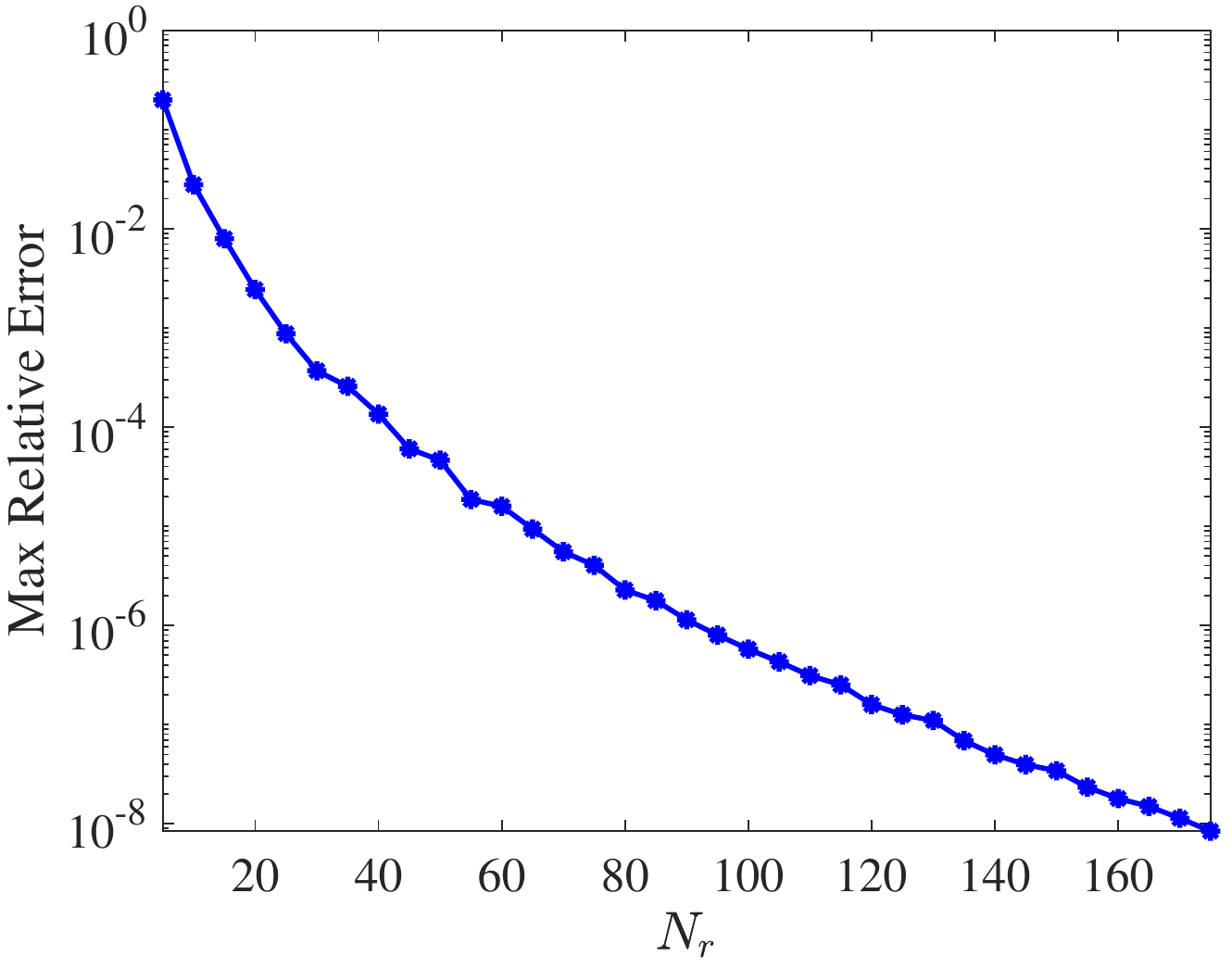}}
		}	
		\caption{Heston test problem with $\bmu=(\kappa, \eta)^T\in\mathbb{R}^2$. Behaviour of the Greedy-POD selection algorithm (left); decay of the relative error (\normalfont{\ref{eq:errBS}}) versus the reduced spaces size (right) for $t_0=0.5$ and $\Lambda=2$. The parametric training set is composed of $|\Xi|=15^2$ equidistantly distributed points.}
		\label{fig9}
	\end{figure}
	\begin{figure}[t]
		\centering{
			\subfigure{
				\includegraphics[width=0.45\textwidth]{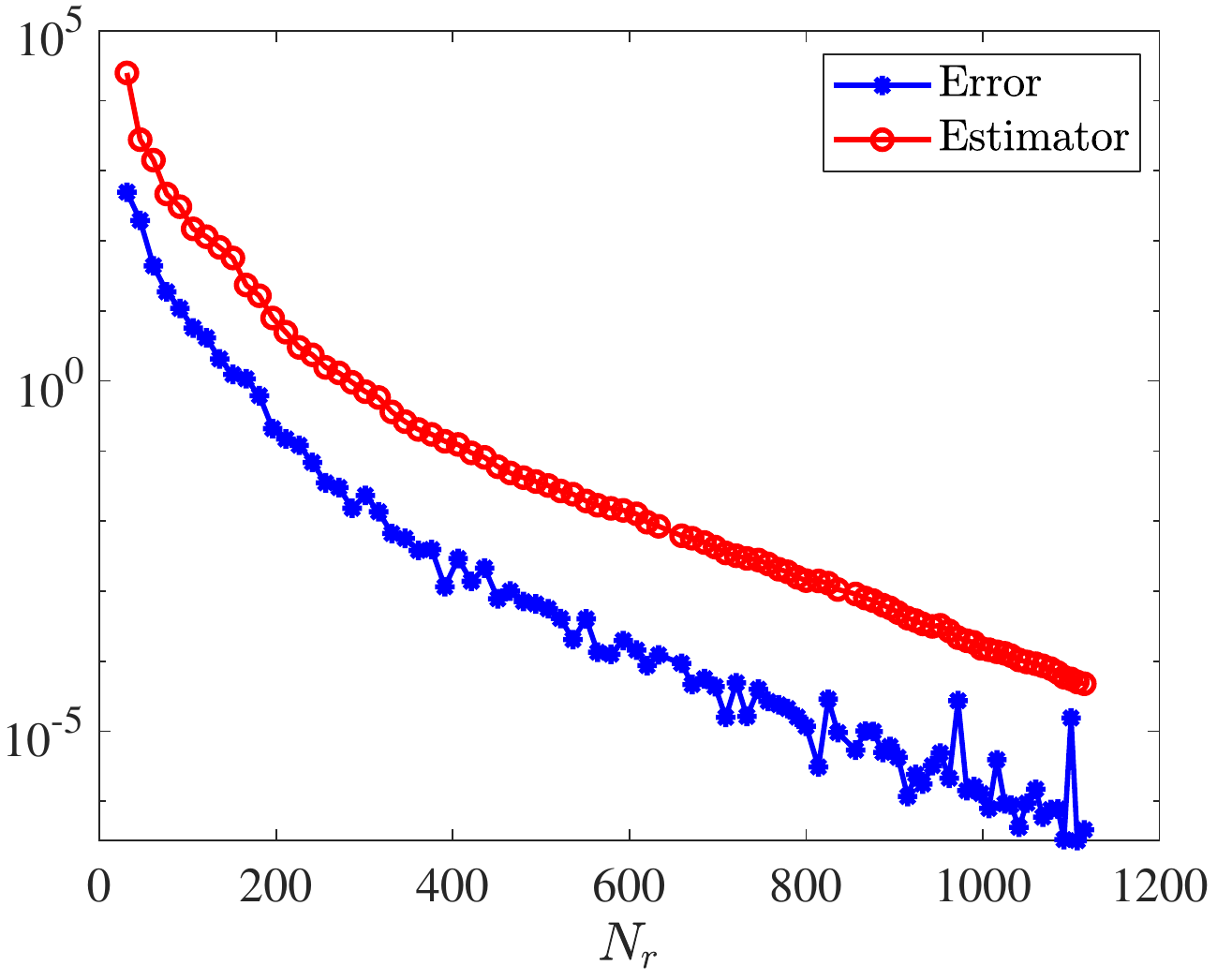}}
			\hspace{0.05cm}
			\subfigure{
				\includegraphics[width=0.45\textwidth]{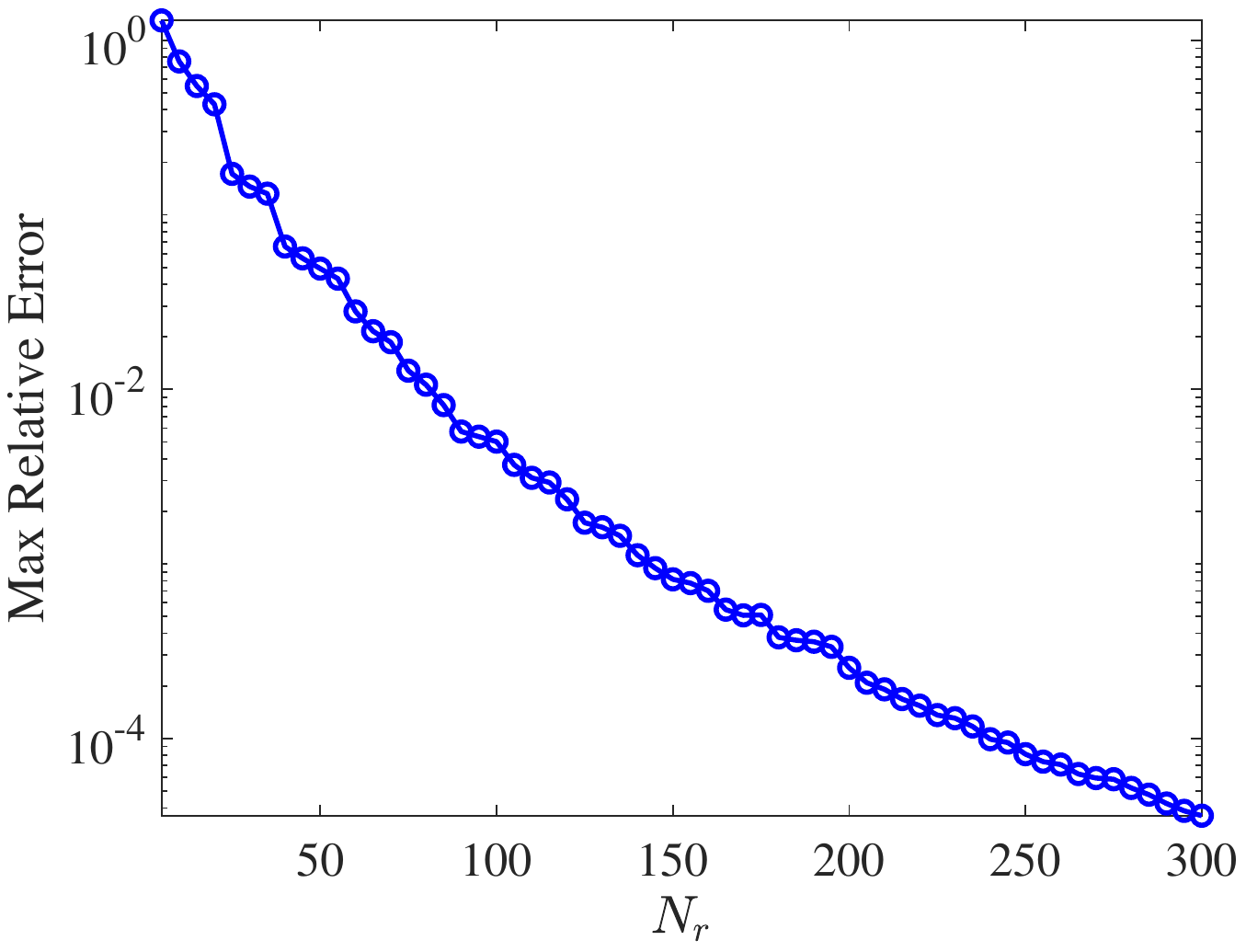}}
		}	
		\caption{Heston test problem with $\bmu=(\sigma,\;r_d,\;\kappa,\;\eta,\;\rho)^T\in\mathbb{R}^5$. Behaviour of the Greedy-POD selection algorithm (left); decay of the relative error (\normalfont{\ref{eq:errBS}}) versus the reduced spaces size (right) for $t_0=0.5$ and $\Lambda=2$. The discrete parametric set is composed of $|\Xi|=5^5$ equidistantly distributed points.}
		\label{fig5}
	\end{figure}

	In Figure \ref{fig9}  we illustrate the effectivity of our Laplace POD-Greedy algorithm based on the error estimator (\ref{eq:est}). Similar to the Black-Scholes problem, in the left picture, we see that the decay of the absolute error is well reproduced by the error estimator, while, in the right picture, we observe the desired exponential decay of the relative error with respect to the size of the reduced spaces. Analogous results are obtained for $\bmu \in \mathcal{D}_1$ and are illustrated in Figure \ref{fig5} where $|\Xi|=5^5=3125$ equidistantly distributed points. We have observed also in this case the exponential decay of the error. However we note that the convergence is slower, which can be explained by the increased complexity in the parameter dependence.
	
	\begin{table}[t]  
	\begin{small}
		\centering
		\begin{tabular}{|c|cc|cc|}
			\cline{2-5}
			\multicolumn{1}{c|}{} & \multicolumn{2}{c}{$d=2$}& \multicolumn{2}{|c|}{$d=5$}\\
			
			\multicolumn{1}{c|}{} & CPU time (s) & \# Snap.& CPU time (s) & \# Snap.\\ 
			\hline
			Algorithm \ref{al1} & 618   & 320&18421&  1115 \\ 
			Algorithm \ref{al2} & 944   & 811&14420&  4848   \\ \hline
		\end{tabular}
		\caption{Heston test problem with $d=2$ and $d=5$ parameters. CPU time and number of stored snapshots (Snap) with the greedy-POD strategy of Algorithm \normalfont{\ref{al1}} and the local greedy strategy of Algorithm \normalfont{\ref{al2}}.}
		\label{T5}
		\end{small}
		\vspace{-0.8cm}
	\end{table}
	\begin{figure}
		\centering{
			\subfigure{
				\includegraphics[width=0.45\textwidth]{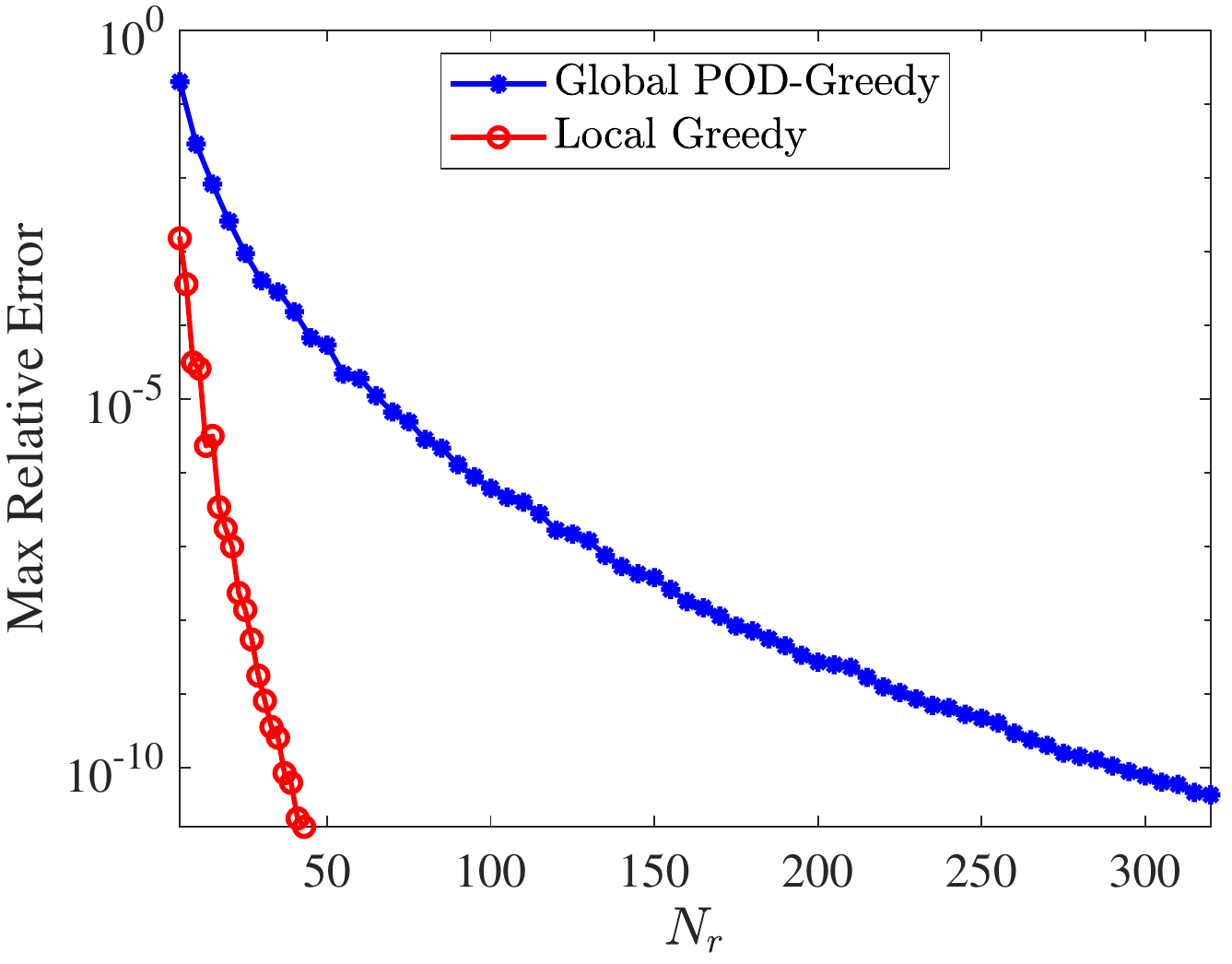}}
			\hspace{0.05cm}
			\subfigure{
				\includegraphics[width=0.45\textwidth]{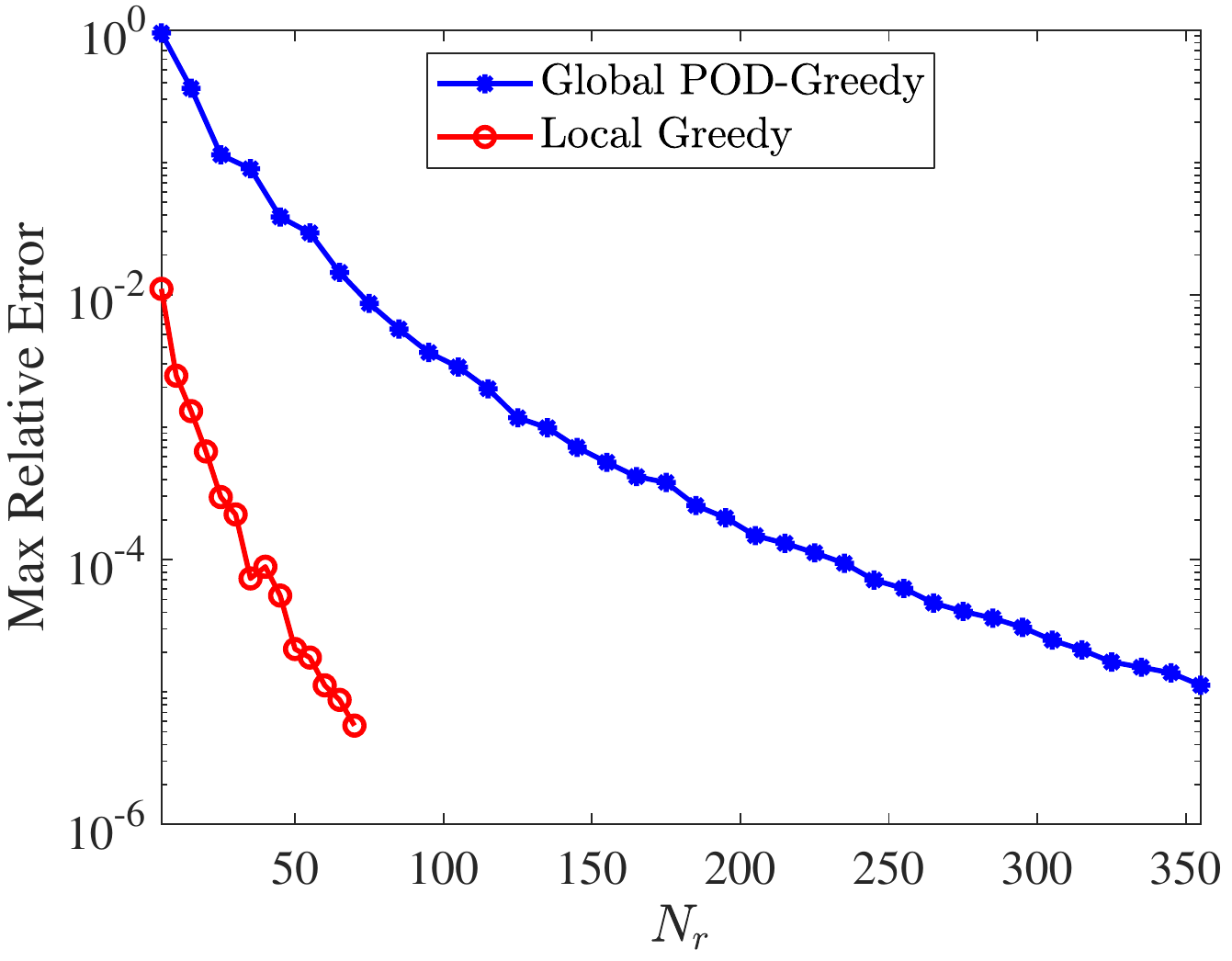}}
		}	
		\caption{Heston test problem with $d=2$ (left) and $d=5$ parameters (right). Decay of the relative reduction error (\normalfont{\ref{eq:errBS}}), versus the reduced spaces size constructed with the POD-Greedy strategy (Algorithm \normalfont{\ref{al1}}) and the local Greedy one (Algorithm \normalfont{\ref{al2}}). Time window $[0.5, 1]$, the number of uniformly distributed samples is $|\Xi|=15^2$ for $d=2$ (left) and $|\Xi|=5^5$ for $d=5$ (right).}
		\label{fig12}
	\end{figure}

	In Table \ref{T5} and Figure \ref{fig12} the performances of Algorithm \ref{al1} and Algorithm \ref{al2} for $d=2$ and $d=5$ are illustrated. The results confirm what has been observed for the Black-Scholes test problem; the local approach is more expensive in terms of computational time and number of stored snapshots but the decay of error with respect to the reduced spaces size is strongly improved. Again, we can conclude that the local approach is better suited for this problem due to the significant improvement in the accuracy reached in the online phase with respect to the reduced spaces size. Moreover for $d=5$ Algorithm \ref{al2} is faster than Algorithm \ref{al1}. This is due to the fact that for $d=5$ the error estimator has a slower decay than $d=2$ and therefore the reduced spaces will result to be larger. Since Algorithm \ref{al2} constructs smaller reduced spaces for each quadrature node with respect to Algorithm \ref{al1} the outcome is that it also becomes more convenient to employ, form a computational point of view, in the offline phase when the dimension of the reduced spaces becomes larger than a few dozen.
	
\begin{figure}[t]
		\centering{
			\subfigure{
				\includegraphics[width=0.45\textwidth]{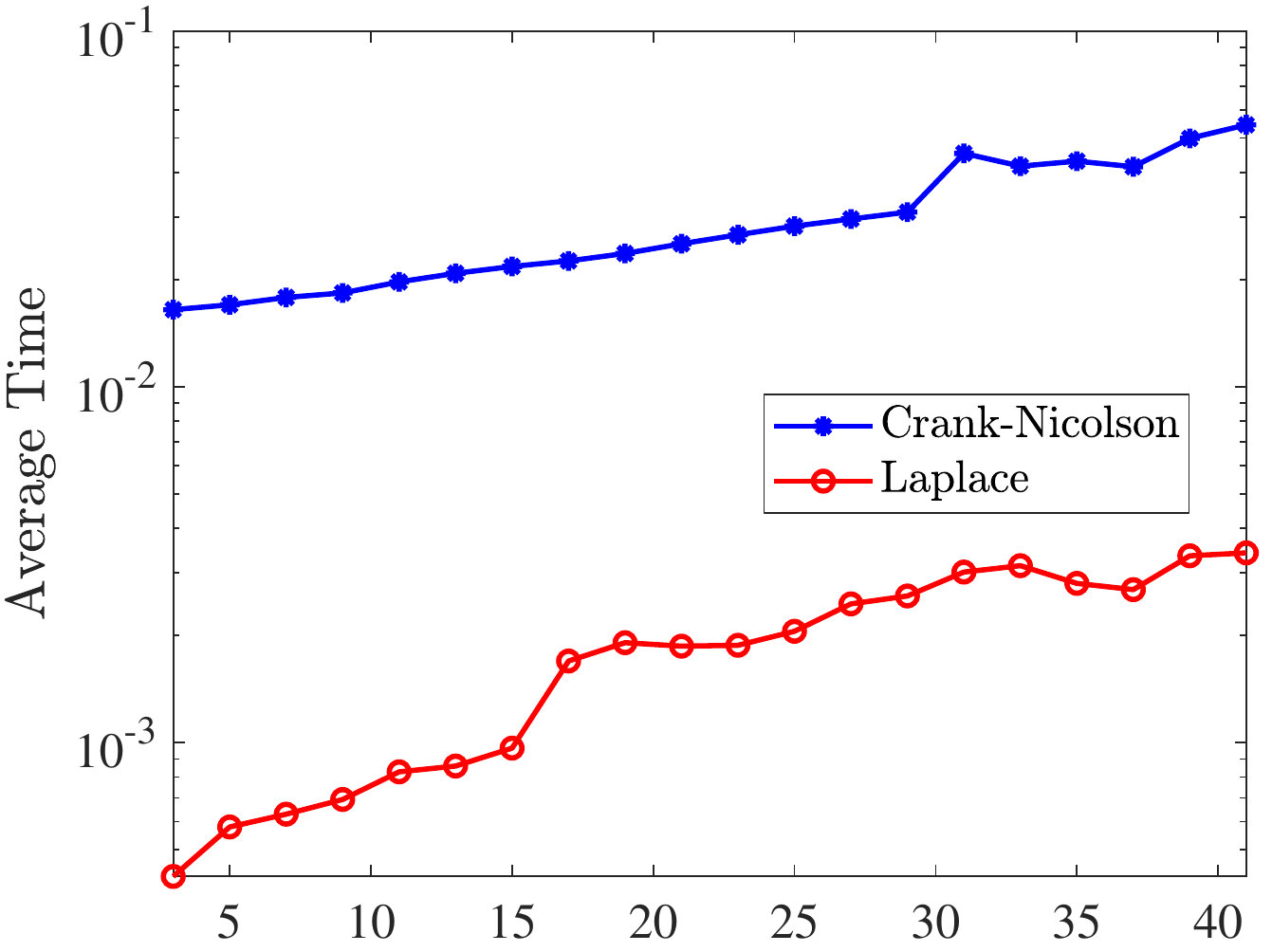}}
			\hspace{0.05cm}
			\subfigure{
				\includegraphics[width=0.45\textwidth]{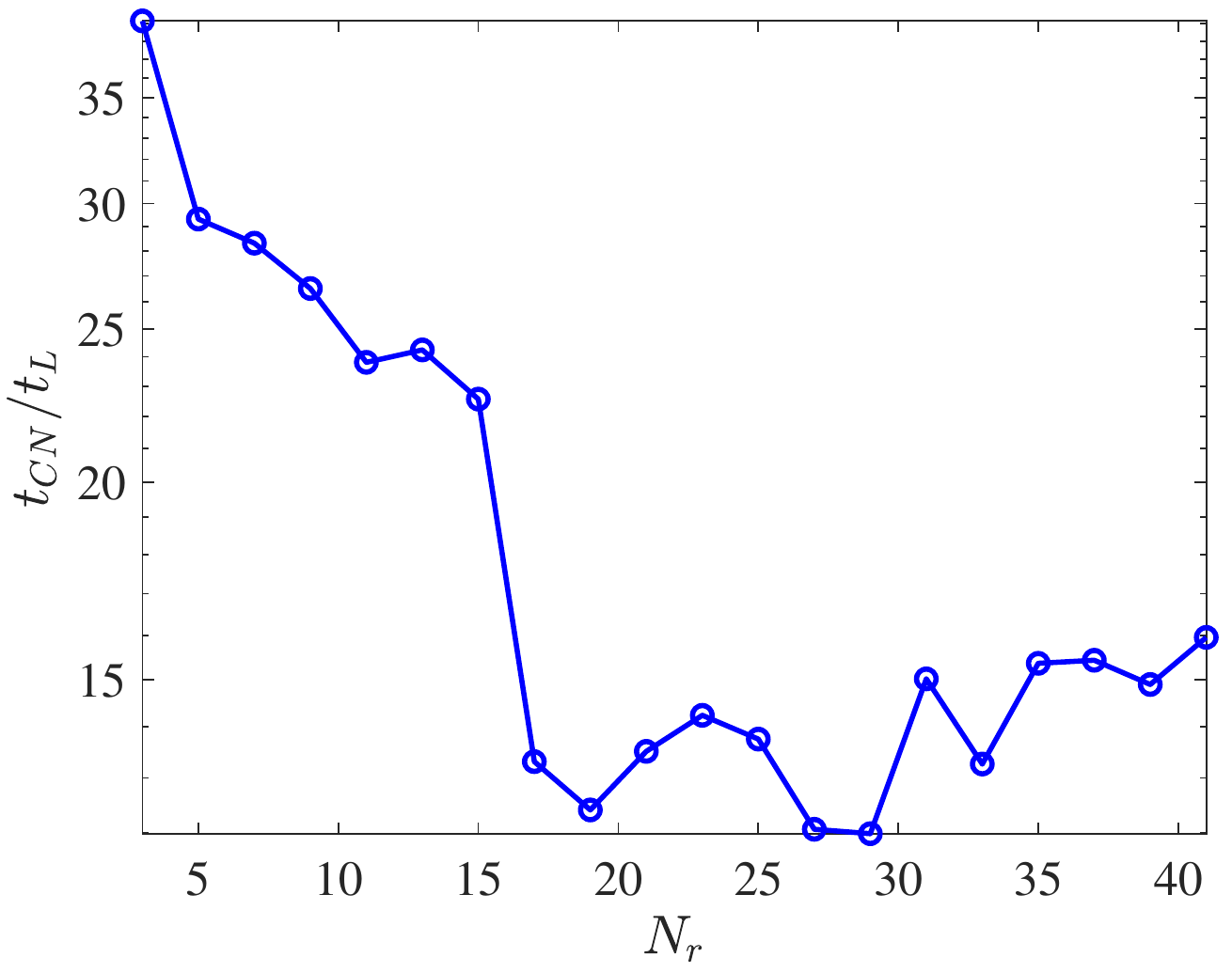}}
		}	
		\caption{Heston test problem. Behaviour of the computational time versus the reduced spaces size for Crank-Nicolson and Laplace methods. Direct comparison (left) and ratio between the averages CPU time (right).}
		\label{fig6}
	\end{figure}

	Also for the Heston problem we compare the Laplace reduced order method with the classical one based on the Crank-Nicolson time scheme (with $\Delta t=10^{-4}$) over the time window $[0.5,1]$. The error in the full time discretization is approximately $10^{-2}$ for both Laplace and Crank-Nicolson. Figure \ref{fig6} shows the average CPU time for the solution of the reduced problem over a test set of $100$ uniformly randomly distributed values of $\bmu \in \mathcal{D}_2$. The contour integral method is between $11$ and $39$ times faster than the classical Crank-Nicolson scheme in the online phase.

\begin{figure}[t]
		\centering{
			\subfigure{
				\includegraphics[width=0.45\textwidth]{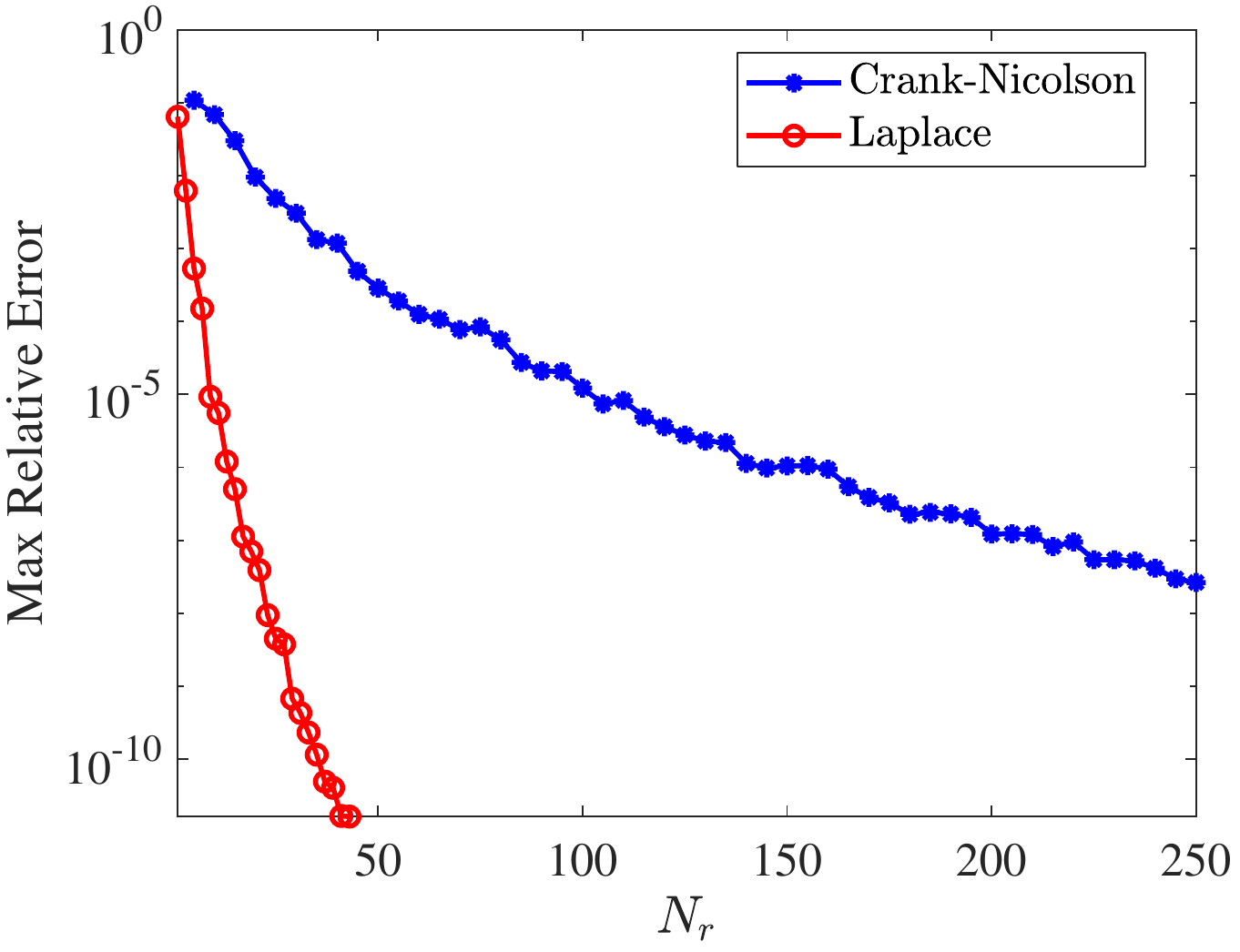}}
			\hspace{0.05cm}
			\subfigure{
				\includegraphics[width=0.45\textwidth]{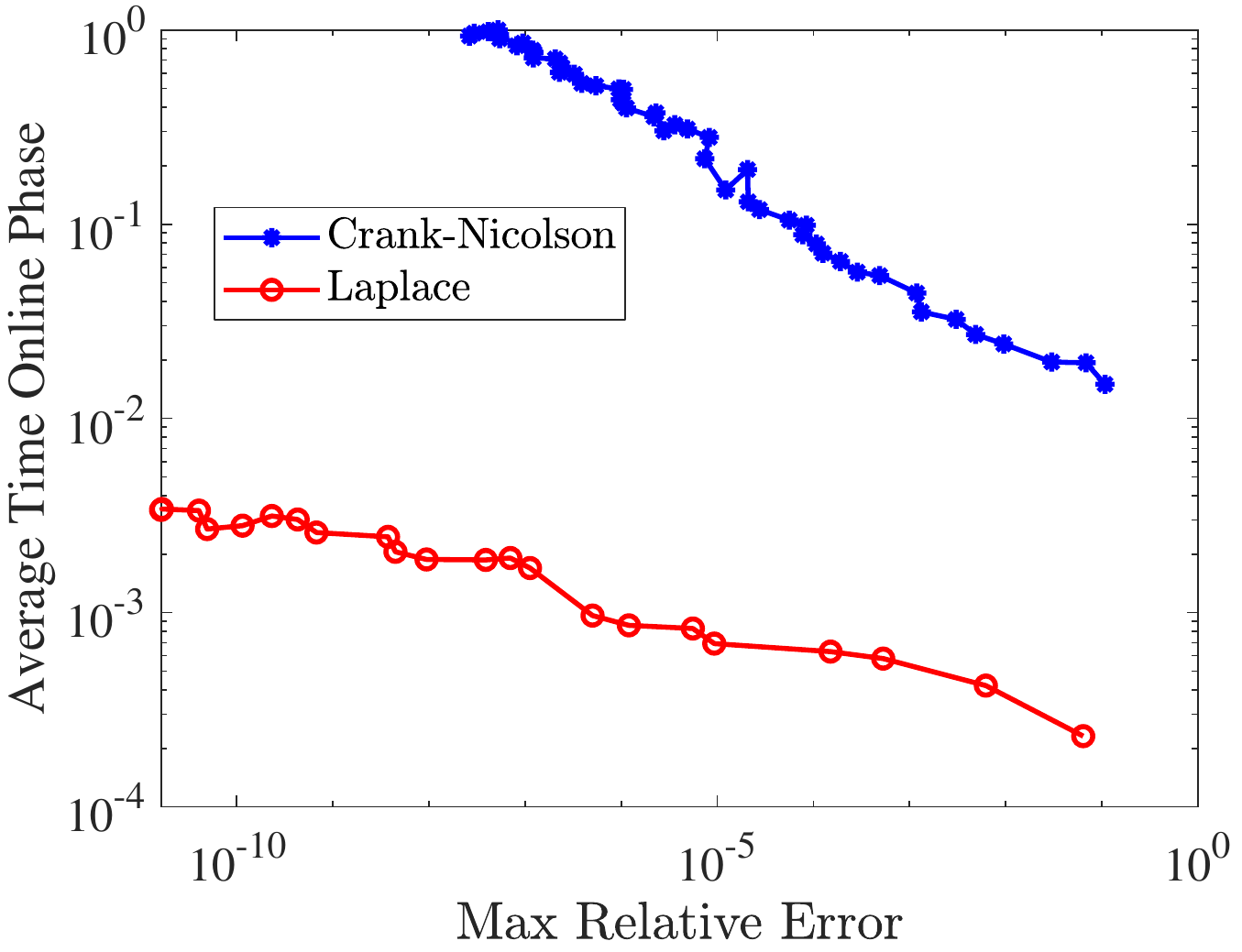}}
		}	
		\caption{Heston test problem. Decay of the relative reduction error (\normalfont{\ref{eq:errBS}}) with respect to the reduced spaces size (left) and computational time with respect to the reduction error (right) for Crank-Nicolson and contour integral methods. Time window $[0.5, 1]$ and $|\Xi|=225$ uniformly distributed sample.}
		\label{fig7}
	\end{figure}

	Figure \ref{fig7} (left) shows the behaviour of the relative error with respect to the size of the reduced spaces for the two reduction techniques. The error originated by the Laplace reduced order method is always significantly smaller than the one associated to the classical reduced order method. Figure \ref{fig7} (right) compares the two reduction techniques in terms of computational time, with respect to the associated relative error. It can be seen that the Laplace reduced order method is two orders of magnitude faster than the classical reduced order method.
	
	Finally in Table \ref{T2} are reported the results of the gradient based strategy for the computation of ${\sigma}_{LB}$ for the Heston problem with $\bmu \in \mathbb{R}^5$. The number of selected starting point is $32$, while for the estimation of ${\sigma}_{LB}$ we made use of the discrete domain $\Xi$ previously defined. It can been seen that with a maximum of $331$ eigenvalue computations the method is able to approximate, with high accuracy, the lower bound on the parametric domain (and in some cases it recover smaller values than the ones associated to the discrete set $\Xi$). The number of solved eigenvalue problems considerably increases with respect to the Black-Scholes test problem (see Table \ref{T1}); despite this fact we observe that the dimension of the parametric domain is much larger and that the cost of our procedure determines a strong speed up with respect to the direct evaluation on the discrete grid $\Xi$. 
	\begin{table}
	\begin{small}
		\centering
		\begin{tabular}{|c|cccc|}
			\hline
			\multirow{2}{*}{$z$}&\multirow{2}{*}{${\sigma}_{GS}$}&\multirow{2}{*}{${\sigma}_{LB}$}&\multirow{2}{*}{$\frac{{\sigma}_{GS}-{\sigma}_{LB}}{{\sigma}_{LB}}$}&\multirow{2}{*}{\#$EP$}\\
			&&&&\\
			\hline
			\hline
			$6.3125 + 0.6398i$   & $0.0958$    &$0.0958$& $0$&$121$ \\
			
			
			$4.0569 + 8.2377i$   & $3.5251$    &$3.5234$& $0.0005$&$212$ \\
			
			
			
			
			$-11.4305 +21.8591i$   & $0.2101$    &$0.2101$& $0$&$131$ \\
			
			$-17.3254 +24.6875i$    & $0.0633$ &$0.0633$& $0$&$131$\\
			
			$-23.8980 +27.2103i$   & $0.0164$    &$0.0165$& $-0.0081$&$160$ \\
			
			$-31.0669 +29.3965i$    & $0.0026$ &$0.0032$& $-0.1955$&$331$\\
			\hline
		\end{tabular}
		\caption{Heston test problem. Approximations of ${\sigma}_{LB}$ through the gradient type method computed in $6$ quadrature points. The table reports the relative error in the approximation and the number of solved eigenvalue problems ($EP$) to compute ${\sigma}_{GS}$.}
		\label{T2}
		\end{small}
		\vspace*{-1.0cm}
	\end{table}
	\vspace{-0.3cm}
	\subsection{The advection equation} As last illustrative example we apply our me\-tho\-do\-lo\-gy to the linear one dimensional advection equation
	\begin{equation}\label{eq:cov}
		u_t+\bmu u_x=0,\quad t\ge0,\quad x\in[0,1],
	\end{equation}
	with initial and boundary conditions
	\begin{equation}\label{eq:covdata}
		u(x,0)=H(x-0.2)=\begin{cases}0 &  x < 0.2\\
			1 & 0.2\le x\le 1\end{cases}, \quad u(0, t)=0,\quad u(1,t)=1.
	\end{equation}
	It has been shown (see \cite{Ohlberger2015}) that for such a problem (\ref{eq:cov})-(\ref{eq:covdata}) the Kolmogorov $N$-width associated to the solution set originated by the map
	\begin{equation*}
		\Phi:\quad \bmu\longrightarrow H(x-0.2-\bmu t),
	\end{equation*}
	has only a polynomial type decay, which suggests that reduced order methods might be inappropriate. We discretize in space (\ref{eq:cov}) with an upwind scheme, setting $\Delta x=10^{-3}$. For the time treatment and successive reduction procedure we compare our approach and the classical reduction based on a backward Euler scheme implemented with constant stepsize $\Delta t=10^{-3}$, and final time $T=0.5$. For the parametric domain we consider $\bmu\in\mathcal{D}=[0.1,1]$ and we consider $\Xi\subset\mathcal{D}$ constructed by taking $100$ uniformly randomly distributed values of the velocity $\bmu$. For this problem we do not employ the greedy strategy, we directly evaluate the full solution taking $20$ values of the velocity $\bmu$ uniformly distributed in $\mathcal{D}$ and we construct the reduced space through the POD decomposition of the snapshots collected. 
	\begin{figure}
		\centering{
			\subfigure{
				\includegraphics[width=0.45\textwidth]{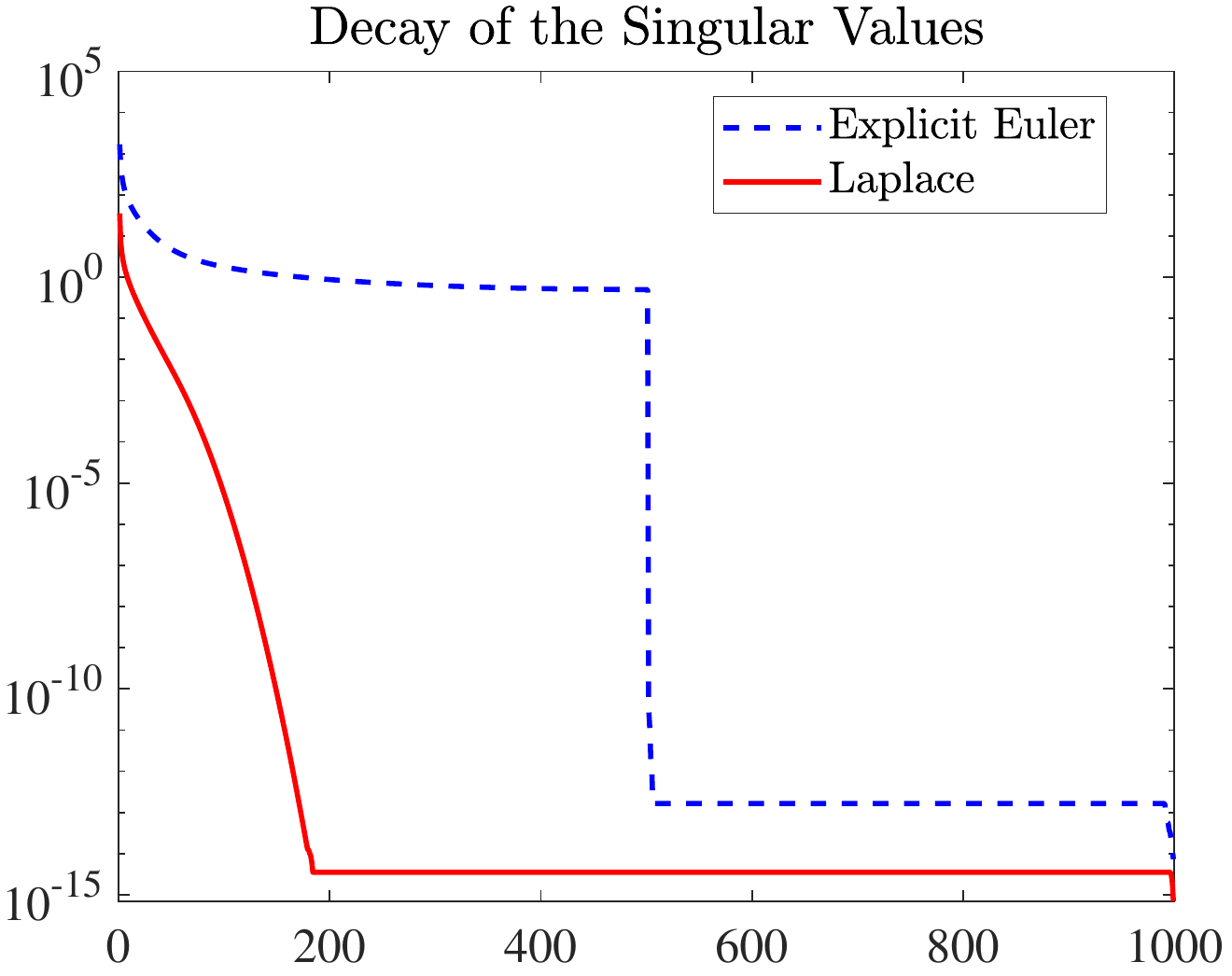}}
			\hspace{0.05cm}
			\subfigure{
				\includegraphics[width=0.45\textwidth]{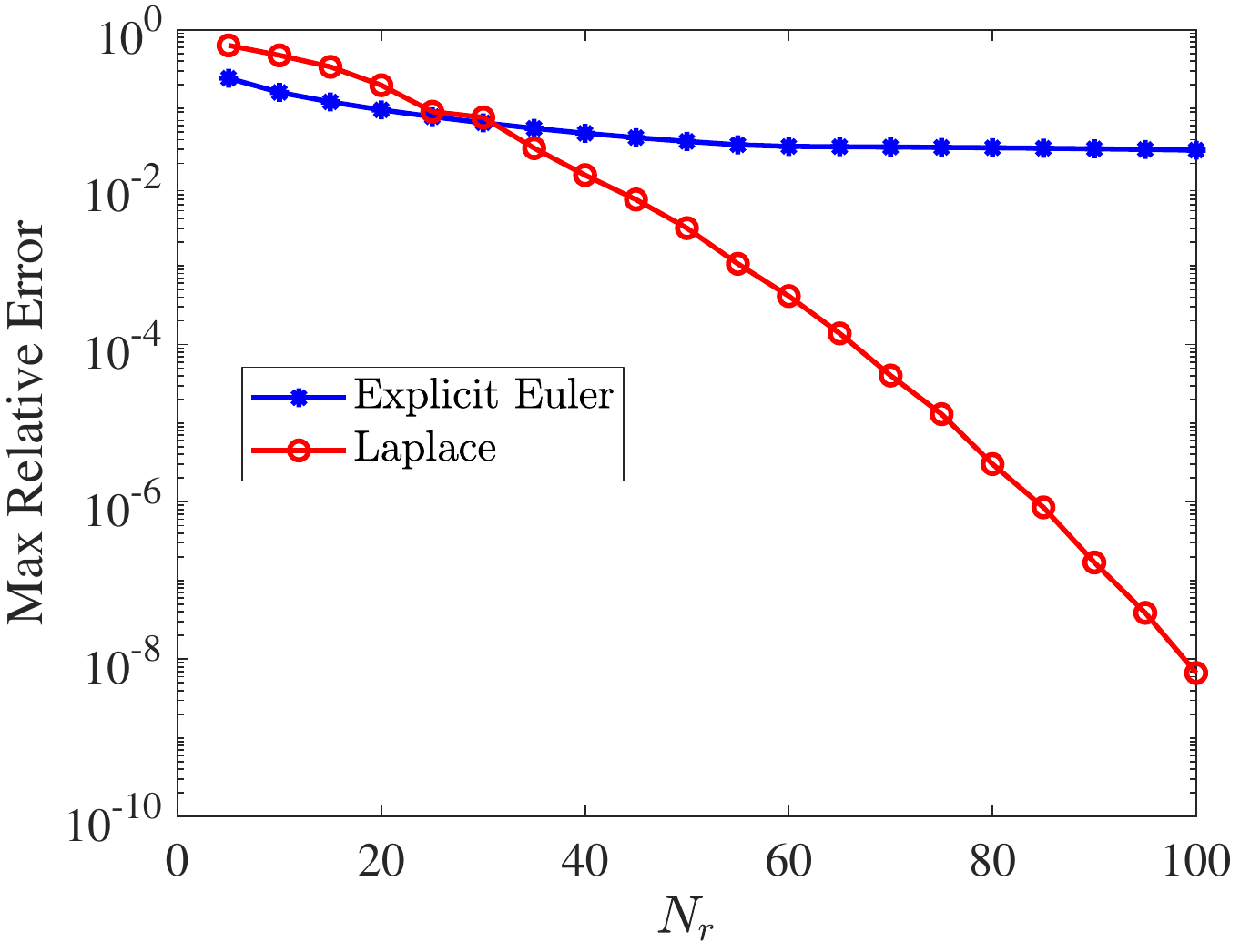}}
		}	
		\caption{Linear 1D advection problem with $\bmu\in[0.1,1]$. Decay of the singular values associated to the snapshots matrix (left) and of the relative reduction error (\normalfont{\ref{eq:errAd}}) (right) for the two reduction techniques.}
		\label{fig1}
	\end{figure}

	Figure \ref{fig1} shows that our reduction strategy is quite effective. The decay of the singular values, which, in some sense, mimics the behaviour of the Kolmogorov $n$-width \cite{Djouadi2008}, is faster. In order to validate our strategy we consider the relative reduction error
	\begin{equation}\label{eq:errAd}
		E_r=\max_{\bmu\in\Xi}\frac{\|\bff u(T;\bmu)-\bff u_r(T;\bmu)\|}{\|\bff u(T;\bmu)\|},
	\end{equation}
	that shows a decay behaviour similar to the one of the singular values.
	
\begin{figure}
		\centering{
			
			\includegraphics[width=0.45\textwidth]{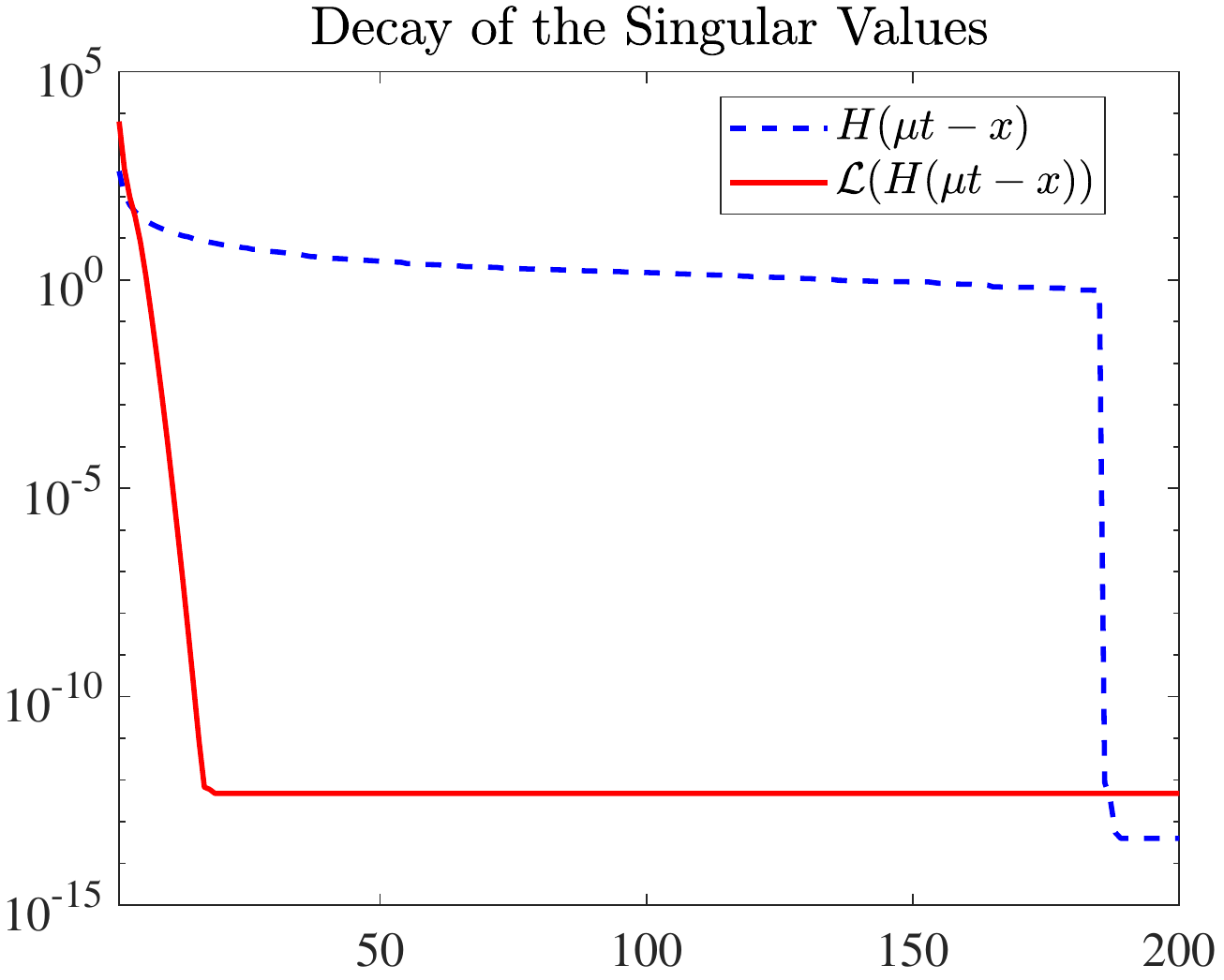}
			
		}	
		\caption{Decay of the first $200$ singular values associated to the snapshots matrices of $H(\bmu t-x)$ and $L(H(\bmu t-a))$ for $\bmu\in[0.5, 5]$ and $x\in [0,1]$. The interval $[0,1]$ is discretized with $1000$ points.}
		\label{fig13}
	\end{figure}

	Note the jump occurring between the $500$th and $501$th singular value for the classical reduced order method, which is due to the presence of $500$ step functions with a different location of the discontinuity point, in the discrete solution set. Indeed, when taking the velocity $\bmu=1$ it holds
	$\bmu({\Delta x}/{\Delta t})=1,$ which implies that if at time $t_i$ the jump is located at $x_i$ then at time $t_{i+1}=t_i+\Delta t$ the jump will be located at $x_{i+1}=x_i+\Delta x$. Since the final time is $T=0.5$ and $\Delta t=10^{-3}$, we have exactly $500$ step functions. Our motivation to employ the simple backward Euler scheme was exactly to bring out this feature. The set of step functions is isomorphic to an orthogonal set of the same size, for this reason the $n$-width has only an algebraic decay (see \cite{Greif2019}). This explains the jump in the blue line of Figure \ref{fig1} (left); the decay of the magnitude of the singular values is slow until all the information related to the $500$ step functions is included in the reduced space. Then the jump in the decay of singular values suggests that all the other elements in the solution set are well approximated by the subset generated by these $500$ step functions. The reason why our strategy is effective in dealing with the reduction phase for the advection problem is that the solution set (\ref{eq:manLap}) is defined on the Laplace transform domain where the step function becomes
	\begin{equation*}
		\mathcal{L}(H(\bmu t-a))=\frac{1}{|\bmu|}\frac{e^{-a\frac{s}{\bmu}}}{s},\quad a\in\mathbb{R},
	\end{equation*}
	which does not lead to a slow decay of the Kolmogorov $n$-width. We show this in Figure \ref{fig13}, where we report the decay of the singular values associated to matrices built taking snapshots of $H(\bmu t-x)$ and $\mathcal{L}(H(\bmu t-a))$ for $\bmu\in[0.5, 5]$, $x\in [0,1]$ making use of a size $N_h=1000$ for the discrete space domain.
	\section{Conclusions}
	In this article we have presented a new approach, based on the coupling of contour integral methods and {projection reduced order methods} to solve time dependent parametric PDEs. The main features of the proposed method, which distinguish it from classical approaches are the following: 
	\begin{itemize}
	    \item[(1) ] it allows the computation of the solution for all times $t$ in a suitable time window, a feature which is preserved by the reduced problem;
		\item[(2) ] it makes possible to parallelize the computation of $\hat{\bff u}(z)$ on the quadrature nodes and the construction of the reduced spaces in the local greedy strategy; 
		\item[(3) ] in the considered test problem (Black-Scholes and Heston) it determines a significant speed up in computational time of the online phase with respect to the classical time stepping scheme;
		\item[(4) ] it allows for efficient treatment in a projection reduced order setting of linear hyperbolic problems with discontinuous initial data. 
	\end{itemize}
	For this reasons we consider the Laplace transform based approach very promising.

	\subsection*{Acknowledgments}

{The authors wish to thank Karsten Urban (Universit\"at Ulm) for his insightful comments and suggestions regarding a preliminary version of this article. Also the} authors acknowledge support by the INdAM Research group GNCS (Gruppo Nazionale di Calcolo Scientifico).
Nicola Guglielmi acknowledges that his research was supported by funds from the Italian MUR (Ministero dell'Universit\`a e della Ricerca) within the PRIN 2017 Project ``Discontinuous dynamical systems: theory, numerics and applications''. 
    \subsubsection*{Code availability} The code is available upon request to the authors.
	\bibliographystyle{plain}
	\bibliography{references}
\end{document}